\documentclass{article}
\usepackage{caption}
\usepackage{rotating}
\usepackage[maxfloats=100]{morefloats}
\usepackage{listings}
\usepackage{booktabs}
\usepackage{xcolor}
\usepackage[title]{appendix}
\usepackage{makecell}
\usepackage{multirow}
\usepackage[symbol]{footmisc}
\usepackage{longtable}
\usepackage{subcaption}
\usepackage[numbers]{natbib}
\usepackage{geometry}
\usepackage{verbatim}
\usepackage{graphicx}
\usepackage{enumitem}
\usepackage{comment}
\usepackage{subfiles}
\usepackage{dsfont}
\usepackage{todonotes}
\usepackage{tikz-qtree,tikz-qtree-compat}
\usepackage{scrextend}
\usepackage{framed}
\usepackage{placeins}
\usepackage{setspace}
\usepackage[ruled,vlined]{algorithm2e}
\usepackage{mathtools,collcell,eqparbox}
\usepackage{multicol}
\usepackage{pgfplots}
\usepackage{tikz}
\usepackage{courier}
\usepackage{url}
\usepackage{color}
\usepackage{authblk}
\usepackage{amssymb}
\usepackage{amsthm}
\usepackage{amsmath}
\usepackage[pagebackref = true]{hyperref}

\hypersetup{
colorlinks = true,
linkcolor = teal,
anchorcolor = teal,
citecolor = teal,
filecolor = teal,
urlcolor = teal
}

\pgfplotsset{compat=1.18}
\numberwithin{equation}{section}

\makeatletter
\newcounter{BMatrix}

\newcommand{\setmaxwd}[1]{%
\eqmakebox[BM-\theBMatrix][\BMalign]{$#1$}%
}
\MHInternalSyntaxOn

\MHInternalSyntaxOff
\makeatother

\newtheorem{theorem}{Theorem}
\newtheorem{assumption}{Assumption}

\newtheorem{proposition}{Proposition}[section]

\newtheorem{corollary}{Corollary}[theorem]
\newtheorem{lemma}{Lemma}
\theoremstyle{definition}

\theoremstyle{remark}
\newtheorem{remark}{Remark}
\newtheorem*{remark*}{Remark}

\newcommand{\inv}{^{-1}}

\newcommand{\1}{\mathds{1}}

\newcommand{\A}{\mathbb{A}}

\newcommand{\R}{\mathbb{R}}

\newcommand{\X}{\mathbb{X}}

\newcommand{\Z}{\mathbb{Z}}
\newcommand{\N}{\mathbb{N}}
\newcommand{\E}{\mathbb{E}}

\newcommand{\eqd}{\overset{d}{=}}

\newcommand{\var}{\mathrm{Var}}

\newcommand{\cov}{\mathrm{Cov}}

\newcommand{\del}{\partial}
\newcommand{\da}{\downarrow}
\newcommand{\ra}{\rightarrow}
\newcommand{\Ra}{\Rightarrow}

\newcommand{\cd}{\cdot}
\newcommand{\ds}{\dots}

\newcommand{\diag}{\mathrm{diag}}

\newcommand{\mrm}[1]{\mathrm{#1}}

\newcommand{\mfk}[1]{\mathfrak{#1}}

\newcommand{\image}{\mathrm{image}}

\newcommand{\cF}{\mathcal{F}}
\newcommand{\cG}{\mathcal{G}}

\newcommand{\cK}{\mathcal{K}}
\newcommand{\cL}{\mathcal{L}}

\newcommand{\cR}{\mathcal{R}}

\newcommand{\cT}{\mathcal{T}}

\newcommand{\set}[1]{\left\{{#1}\right\}}

\newcommand{\floor}[1]{\left\lfloor{#1}\right\rfloor}

\newcommand{\norm}[1]{\left\|#1\right\|}

\newcommand{\abs}[1]{\left|#1\right|}
\newcommand{\sqbk}[1]{\left[ #1 \right]}
\newcommand{\sqbkcond}[2]{\left[ #1 \middle| #2 \right]}

\newcommand{\crbk}[1]{\left( #1 \right)}

\newcommand{\pmx}[1]{\begin{pmatrix} #1 \end{pmatrix}}

\title{Understanding the Asymptotic Transience of Q-Learning\\ via Distributional Expansion}
\author{Shengbo Wang}
\affil{Daniel J. Epstein Department of Industrial and Systems Engineering, University of Southern California}
\date{November 2026}

\begin{document}
\newcommand{\norminf}[1]{\norm{#1}_\infty}
\maketitle
\begin{abstract}
Q-learning has been empirically observed to exhibit substantial transient error. This transient remains uncaptured by the leading Gaussian limit because the central-limit normalization removes the smaller-scale displacement generated by the initial nonlinear evolution. On the other hand, worst-case finite-sample guarantees, while accounting for transient effects, can obscure the mechanisms responsible for the error on a particular problem instance. We develop a finer asymptotic analysis that gives an asymptotically exact characterization of this transient for a fixed instance. Specifically, for synchronous Q-learning in finite discounted Markov decision processes, we separate the last-iterate error into an optimal-policy temporal-difference (TD) component and a nonnegative displacement caused by maximization. The centered TD component has Gaussian fluctuations at scale $n^{-1/2}$, whereas the maximization transient is of order $n^{-1}$ and converges almost surely, after rescaling by $n$, to a random limit. This limit represents the propagated effect of maximization errors generated before the greedy policy stabilizes. We combine this characterization with an Edgeworth expansion of the TD component to identify the corrections to its Gaussian approximation, separating the effects of initialization, non-Gaussian TD fluctuations, and maximization. The results reveal how transient errors arise for individual Q-learning instances and how their effects persist below the central-limit scale.
\end{abstract}

\section{Introduction}
\label{sec:introduction}

The convergence of a learning algorithm is commonly studied through finite-sample bounds or asymptotic limits. These approaches can provide different information about the same stochastic dynamics. Finite-sample analysis can quantify how initialization, noise, and step sizes affect convergence before the algorithm reaches its asymptotic regime. However, a guarantee that holds uniformly over a broad class of problems must accommodate its most difficult instances. Thus, even a sharp worst-case bound need not explain the behavior of a particular instance, and bounding the accumulated error can obscure the mechanism, direction, and sample-path dependence of its transient component. This distinction is particularly relevant to Q-learning, for which both sharp worst-case guarantees and substantial variation across problem instances have been established \citep{wainwright2019cone,li2024minimax,wang2026optimalsamplecomplexitystable,xia2026instance}.

Classical asymptotic analysis provides a complementary description. Under suitable stability and regularity conditions, stochastic approximation (SA) algorithms usually admit Gaussian limits determined by the local drift and noise near the limiting point \citep{fabian1968normal,kushner2003stochastic,borkar2025asymptotic}. This description identifies the exact leading fluctuation scale and its covariance. The normalization that produces an asymptotic limiting law can nevertheless remove the effect of the earlier nonlinear evolution. Consequently, for a learning algorithm, the central limit theorem (CLT) alone does not determine whether the initial transient leaves a systematic displacement, how that displacement depends on the sample path, or at what smaller scale it remains visible.

These observations motivate the question studied in this paper: \textit{can a finer asymptotic analysis identify the leading effect of the transient phase of a learning algorithm?}

We address this question for synchronous Q-learning in a finite discounted Markov decision process. The maximization in the Bellman update makes the recursion nonlinear. However, if every state has a unique optimal action, the Bellman operator is affine in a neighborhood of the optimal action-value function $Q^*$. Almost-sure convergence of Q-learning then implies that the greedy decisions agree with the optimal policy after a finite random time. From that point onward, Q-learning uses the same update rule as policy evaluation under the optimal policy. Even though the maximization will not induce an error compared to the optimal policy evaluation recursion from that point onward, the displacement induced by the earlier maximizations have already changed the iterate, and the subsequent updates continue to propagate that displacement. Our analysis identifies the asymptotic magnitude and structure of this accumulated effect.

Specifically, to describe the results, let $\gamma\in(0,1)$ be the discount factor, let $P$ be the transition operator from state values to state-action values, and write $\Pi^*q(x)=q(x,\pi^*(x))$ for evaluation at the unique optimal action. At iteration $k$, the generative model supplies a bounded reward vector $R_k$ and a sampled transition operator $\widehat P_k$, with independent and identically distributed observations across iterations and $\mathbb E\widehat P_k=P$. Consider stepsize $\alpha_k \sim 1/((1-\gamma)k)$. With $V_k(x)=\max_a Q_k(x,a)$, we couple Q-learning to a policy-correct recursion using the same initialization and data: $\widetilde Q_0=Q_0$, 
\[
\begin{aligned}
Q_{k+1}&=(1-\alpha_{k+1})Q_k+\alpha_{k+1}\bigl(R_{k+1}+\gamma\widehat P_{k+1}V_k\bigr),\\
\widetilde Q_{k+1}&=(1-\alpha_{k+1})\widetilde Q_k+\alpha_{k+1}\bigl(R_{k+1}+\gamma\widehat P_{k+1}\Pi^*\widetilde Q_k\bigr).
\end{aligned}
\]
The second recursion is the synchronous temporal-difference (TD) learning for the optimal policy, which we call the policy-correct recursion. It provides a benchmark that preserves the observation noise and isolates the effect of replacing policy evaluation with maximization.

We decompose the iteration-$n$ error as $Q_n-Q^*=Y_n+m_n+D_n$, where
\[
Y_n=\widetilde Q_n-\mathbb E\widetilde Q_n,\quad m_n=\mathbb E\widetilde Q_n-Q^*,\quad D_n=Q_n-\widetilde Q_n.
\]
Here $Y_n$ is the centered policy-correct fluctuation, $m_n$ is the deterministic bias of the policy-correct recursion, and $D_n$ is the displacement caused by maximization. The latter can be written as
\[
D_n=\gamma\sum_{k=0}^{n-1}\alpha_{k+1}C_nC_{n-1}\cdots C_{k+2}\widehat P_{k+1}(V_k-\Pi^*Q_k),
\]
where $C_k$ is a nonnegative random matrix. Each summand represents a maximization residual generated at iteration $k$ and propagated through subsequent updates. Since every factor is nonnegative, $D_n\ge0$ on every sample path. We call $D_n$ the \emph{maximization transient}.

Our first result characterizes this displacement at order $n^{-1}$. We show that $nD_n\ra D$ almost surely and in $L^p$ for every $1\le p<\infty$ to some random variable $D$. To interpret this limit, recall that $Q_n\ra Q^*$ almost surely and the optimal action at each state is unique. These properties imply that there exists an almost surely finite random time after which the greedy decisions agree with the optimal policy. Consequently, $V_k-\Pi^*Q_k=0$ for all subsequent iterations, so only finitely many indices contribute to the sum as $n\ra\infty$. Although no new maximization residuals are generated after this random time, later transition samples continue to modulate the accumulated displacement through the matrix products $C_n\ds C_{k+2}$. The limit $D$ therefore records the lasting effect of the maximization errors generated during this finite initial phase, together with their subsequent propagation.

The magnitude of this $n^{-1}$ transient can nevertheless be large. In a one-state, two-action example with bounded rewards, we show that $\mathbb ED$ is bounded below by a positive constant times the inverse optimality gap. This provides an asymptotic explanation consistent with observations that maximization bias can substantially affect Q-learning performance \citep{vanhasselt2010double}. In particular, even though the maximization transient is at the smaller scale $n^{-1}$, it can still generate substantial displacement at reasonably large iteration counts. The inverse-gap dependence clarifies one source of this difficulty: when competing actions are hard to distinguish, maximization errors can accumulate during the initial phase and leave a large displacement under subsequent propagation. Thus, our fixed-instance asymptotic analysis identifies this displacement pathwise and characterizes its limiting effect. In comparison, a finite-sample minimax guarantee that is uniform over the action gap must cover such hard instances, hence the resulting worst-case bound can obscure how their large errors arise. By isolating the accumulated maximization residuals and their propagation, our analysis makes this mechanism explicit.

Identifying $D$ does not by itself give a complete distributional approximation at its scale. By the CLT for SA algorithms, the centered policy-correct iterate $Y_n$ is of order $n^{-1/2}$, whereas the maximization transient is of order $n^{-1}$. To examine the contribution of $Y_n$ at the $n^{-1}$ scale, we recall that the first order correction in the classical Edgeworth expansion directly corresponds to this finer scale. Therefore, we develop a multivariate Edgeworth expansion for the centered policy-correct iterate $Y_n$ and characterize the relevant Gaussian corrections. At this finer scale, however, dependence between the policy-correct fluctuation and the random maximization transient could affect the resulting distributional correction. To clarify this potential interaction, we establish the joint convergence $(\sqrt nY_n,nD_n)\Rightarrow(Z,D)$, where $Z\sim N(0,\Sigma)$ and $Z$ is independent of $D$. This asymptotic independence is crucial for combining the two contributions.

Specifically, when $\Sigma$ is positive definite, we show the following distributional expansion for the $\sqrt{n}$-normalized error; formally, 
\[
\sqrt n(Q_n-Q^*)\sim \phi_\Sigma(z)\left[1+n^{-1/2}\crbk{P_1(z)+\mu^\top\Sigma^{-1}z}+o\crbk{n^{-1/2}}\right]dz.
\]
Here, $\phi_\Sigma$ denotes the Gaussian density with covariance $\Sigma$, $P_1$ is the first centered Edgeworth polynomial for $\sqrt nY_n$ specified in Theorem \ref{thm:edgeworth}, and $\mu=m+\mathbb ED$, where $m=\lim_{n\to\infty}nm_n$ is the limiting coefficient of the deterministic initialization bias. The rigorous version of this result is stated in Theorem \ref{thm:combined}. In particular, although $D$ is random, its contribution to the first correction is additive and depends only on $\mathbb ED$. The expansion therefore captures the effect of the early maximization phase together with the Gaussian correction polynomial that are relevant at the $n^{-1}$ scale.

Beyond our contribution to understanding how the maximization transient arises, propagates, and interacts with the Gaussian correction, the paper develops a framework and explicit formulas for the Edgeworth expansion of multivariate affine SA with multiplicative and martingale-difference noise. This extends the existing work on univariate SA by \citet{dippon2008edgeworth}. We also note that the affine structure is not essential to our approach. Under suitable smoothness and stability conditions, we expect the analysis to extend to nonlinear SA when \(b\) is sufficiently large, with higher-order derivatives of the drift at the limiting point incorporated into the Gaussian correction coefficients.

This work finds inspiration from two lines of related literature. \citet{srikant2025rates} uses Stein's method and Poisson equations to establish quantitative central limit theorems for martingales and Markov-chain functionals, with an application to Polyak--Ruppert averaged TD learning. These results provide bounds on the distributional error from the Gaussian limit, whereas our expansion identifies the explicit corrections up to the scale of the maximization transient; see Theorem \ref{thm:edgeworth}. 

A second line develops expansions for constant-step stochastic approximation. \citet{huo2026bias} expand the stationary bias of Markovian linear stochastic approximation in powers of the step size, and \citet{zhang2024constant} derive a first-order stationary-bias expansion for asynchronous Q-learning under a unique optimal policy. Their limiting object is the stationary law obtained at a fixed step size, followed by an expansion as that step size decreases. Our limiting object retains the effect of the initial nonlinear evolution along a diminishing-step recursion converging to $Q^*$: The almost-sure limit $nD_n\ra D$ identifies the random displacement along individual trajectories, and its joint convergence with an independent Gaussian fluctuation determines how this displacement enters the distributional expansion.

\subsection{Review of Related Literature}
\label{subsec:related-literature}

\paragraph{Stochastic approximation and Q-learning asymptotics.}
Q-learning \citep{watkins1992q} can be viewed as a SA algorithm \citep{robbins1951stochastic}. Its convergence was developed through asynchronous stochastic approximation and stochastic dynamic programming by \citet{tsitsiklis1994asynchronous} and \citet{jaakkola1994convergence}. The ODE approach of \citet{borkar2000ode} provides a general route from dynamical stability to convergence; \citet{kushner2003stochastic} and \citet{borkar2008stochastic} give systematic treatments. Classical asymptotic normality results for SA algorithms include \citet{fabian1968normal}, while \citet{borkar2025asymptotic} develop functional central limit theorems and asymptotic covariance results under Markovian noise. For Q-learning, \citet{devraj2021qlearningvariance} analyze asymptotic covariance, step-size design, and relative Q-learning. 

\paragraph{Finite-sample bounds for Q-learning.}
\citet{evendar2003learning} establish the sensitivity of Q-learning convergence rates to the learning-rate schedule. \citet{wainwright2019cone} develops sample-complexity bounds through cone-contractive SA, including instance-dependent refinements. \citet{chen2022finite} establish finite-sample convergence bounds for contractive SA using Lyapunov approaches and apply them to Q-learning. \citet{li2024minimax} establish matching upper and lower bounds for Q-learning, showing that Q-learning is not minimax optimal. To address this gap, \citet{wainwright2019variancereducedqlearningminimaxoptimal} proposes variance-reduced Q-learning and establish its minimax optimality. \citet{xia2026instance} identify instance-dependent functionals governing the statistical risks of optimal value estimation through local minimax lower bounds. They also establish a matching instance-dependent sample complexity upper bound for variance-reduced Q-learning.

\paragraph{Maximization bias.}
\citet{vanhasselt2010double} identifies the overestimation produced by maximizing noisy value estimates and introduces Double Q-learning to reduce it. \citet{vanhasselt2016deep} generalize this approach to deep reinforcement learning. The worst-case analysis of \citet{li2024minimax} further connects overestimation with the statistical inefficiency of Q-learning. In this paper, we identify the maximization transient as a pathwise nonnegative displacement relative to the optimal-policy evaluation recursion dynamics and characterize its asymptotic limit.

\paragraph{Gaussian approximation bounds and stationary expansions.}
Beyond \citet{srikant2025rates}, \citet{samsonov2024gaussian} establish multivariate Berry--Esseen bounds and multiplier-bootstrap guarantees for Polyak--Ruppert averaged linear stochastic approximation, with applications to TD learning. In the constant-step setting, \citet{huo2024collusion} separate the leading stationary bias into contributions from Markovian dependence, nonlinearity, and their interaction. \citet{zhang2026prelimit} study steady-state limits of nonsmooth contractive stochastic approximation, including square-root-step-size bias. \citet{zhang2026steadystate} obtain quantitative steady-state Gaussian approximation bounds under Markovian multiplicative noise and local quadratic linearization assumptions. These results motivate structural analyses beyond error rates, considering finer approximation error bounds.

\paragraph{Edgeworth expansions.}
The classical theory of Edgeworth expansions uses cumulants to construct polynomial corrections to Gaussian laws; see \citet{bhattacharya2010normal,hall1992bootstrap}. Edgeworth expansions for SA algorithms have also been explored. \citet{dippon2004representations} derives higher-order representations of \textit{scalar} Robbins--Monro iterates using linear, quadratic, and cubic forms in the observation errors. \citet{dippon2008edgeworth} establishes Edgeworth expansions for the resulting approximating sequences, and \citet{dippon2008asymptotic} obtains expansions for the Robbins--Monro process itself. In this paper, we extend and develop a multivariate Edgeworth expansion for synchronous TD-learning, where the techniques generalize to nonlinear SA with martingale difference noise. More recently, \citet{konakov2026local} develop local limit theorems and strong approximations for Robbins--Monro procedures.

\section{Q-Learning and Error Decomposition}
\label{section:setup}

\subsection{Markov Decision Processes}

We work in the finite-state, finite-action synchronous generative-model setting. Let $\X$ and $\set{\A_x:x\in\X}$ be finite and set $\Z=\set{z=(x,a):x\in\X,a\in\A_x}$. Fix $\gamma\in(0,1)$ and a controlled transition kernel $P(\cdot| \cd)$. We view this kernel as a linear operator from $\R^\X\ra\R^{\Z}$; i.e. for $v\in\R^{\X}$, $Pv(z)=\sum_{x'\in\X}P(x'|z)v(x'),$ for all $z\in\Z$.

Let $F_k:\Z\times\X\to\R$ be i.i.d. random reward functions, independent of the transition samples, such that
\[
\max_{z\in\Z,x'\in \X}|F_k(z,x')|\le R_{\max}.
\]
Then, for a controlled Markov chain $\set{(X_k,A_k)\in\Z:k\ge 0}$ generated by any admissible policy $\pi$ and initial state $x\in\X$, the expected discounted cumulative reward is \[V^\pi(x) =E^{\pi}\sqbkcond{\sum_{k=0}^\infty \gamma^k F_k(X_k,A_k,X_{k+1})}{X_0=x},\]
where the expectation is taken over the joint law of the controlled Markov chain under $(P,\pi)$ and the i.i.d. sequence $\set{F_k:k\ge 0}$. This is also known as the value function, seen as a function in $x$, associated with policy $\pi$.

Given this setup, we seek an optimal admissible policy that maximizes $V^\pi(x)$. Classical dynamic programming theory for MDPs provides a value-based approach. Specifically, we define the mean reward as
\[
r(z) = \sum_{x'\in\X} P(x'|z)E F_1(z,x'),\quad\forall z\in\Z
\]
where the expectation is taken with respect to the law of $F$. Consider the unique solution $Q^*$ to the following Bellman equation:
\[
Q^*(x,a) = r(x,a)+\gamma \sum_{x'\in\X} P(x'|x,a)\max_{b\in\A_x}Q^*(x',b).
\]
Then, for all $x\in\X$, $V^\pi(x)$ is simultaneously maximized by any deterministic Markovian policy that selects $\pi(x)\in \arg\max_{a\in\A_x} Q^*(x,a)$ at state $x\in\X$. The function $Q^*$ is called the optimal action-value function, or simply the optimal Q-function. For a comprehensive exposition of dynamic programming theory, see \citet{puterman2014MDP}.

This property of $Q^*$ motivates many modern Markovian control and RL algorithms, which construct an estimate of $Q^*$ from data and use its argmax as the deployed policy. Among these approaches, Q-learning \citep{watkins1992q} has been receiving significant attention from both practitioners and theorists. It enjoys algorithmic simplicity, making it readily compatible with function approximation, as well as favorable theoretical properties, such as asymptotic normality. In this paper, we seek a deeper understanding of the asymptotic behavior of a simple variant known as synchronous Q-learning.

\subsection{Synchronous Q-Learning and its Central Limit Theorem}

In this paper, we consider synchronous Q-learning under a generative model that produces transition and reward data. Let $(\Omega,\cF,\Pr)$ be the supporting probability space for the generative model, and let $\E$ denote the expectation under $\Pr$. At each iteration $k\ge1$, the generative model draws a vector $\set{X_k(z)\in\X:z\in\Z}$, where $X_k(z)\sim P(\cdot|z)$ independently across $z\in\Z$ and across iterations. At the same time, we also observe the entire reward vector
\[
R_k(z)=F_k(z,X_k(z)),\quad \forall z\in\Z.
\]
However, the reward function $F_k$ itself is not assumed to be observable. 

Then, by the definition of $r$, we have $r(z)=\E R_k(z)$. We also define the sampled transition kernel by $\widehat P_kv(z)=v(X_k(z))$ for $v:\X\ra \R$, which can be represented as a matrix in $\R^{|\Z|\times|\X|}$ with entries taking values in zero and one. It follows that $\E\widehat P_k=P$. Let $\cF_k = \sigma((\widehat P_j,R_j):j\leq k)$ denote the filtration generated by the observed data.

Next, we introduce the synchronous Q-learning algorithm. Write $V_k(x)=\max_{a\in\A_x}Q_k(x,a)$. The synchronous Q-learning maintains an estimate of $Q^*$ by iterating the following (vector) update
\begin{equation}
Q_{k+1}=(1-\alpha_{k+1})Q_k+
\alpha_{k+1}\sqbk{R_{k+1}+\gamma\widehat P_{k+1}V_k}
\label{eqn:q-update}
\end{equation}
starting from a deterministic initialization $Q_0$. 

\begin{assumption}[Unique optimal action]
\label{assump:unique_optimal_action} Assume that every state has a unique optimal action $\pi^*(x) := \arg\max_{a\in\A_x}Q^*(x,a)$.
\end{assumption}
Under Assumption \ref{assump:unique_optimal_action}, the finiteness of the state and action spaces implies a strictly positive action gap
\[
\delta^*=\min_{x\in\X}\min_{\pi^*(x)\ne a\in\A_x}
\sqbk{Q^*(x,\pi^*(x))-Q^*(x,a)}>0.
\]
Define the optimal-policy application map
\[
\Pi^* q(x)=q(x,\pi^*(x)), \quad\forall x\in\X
\]
which is a linear operator, seen as a $\R^{|\X|\times|\Z|}$ matrix of zeros and ones. Moreover, for any controlled transition kernel $P\in\R^{|\Z|\times|\X|}$, $P\Pi^*$ is a square row-stochastic matrix.

\begin{assumption}[Rescaled-linear step sizes]
\label{assump:stepsize}
We use the rescaled-linear step size $\alpha_k:=\alpha/(k+k_0)$ and assume that $k_0+1\ge\alpha$ and $b:=\alpha(1-\gamma)>1/2$.
\end{assumption}

The offset condition ensures that $0<\alpha_k\le1$ and $k_0+1>b$. The step sizes satisfy the Robbins-Monro conditions $\sum_{k\ge1}\alpha_k=\infty$ and $\sum_{k\ge1}\alpha_k^2<\infty$. Consequently, by \citet[Proposition~4.1]{abounadi2002stochastic}, the iterates in \eqref{eqn:q-update} converge almost surely to $Q^*$. This almost-sure convergence does not require a unique optimal action or the restriction $b>1/2$.

Under Assumption~\ref{assump:unique_optimal_action}, the drift $q\mapsto\cT q-q$ is affine in a neighborhood of $Q^*$, with Jacobian $\gamma P\Pi^*-I$. For rescaled-linear step sizes, the standard spectral stability condition for the $\sqrt n$-normalized stochastic approximation recursion requires every eigenvalue of $\alpha(\gamma P\Pi^*-I)+I/2$ to have strictly negative real part; see \citet{kushner2003stochastic}. Since $P\Pi^*$ is stochastic, the largest real part of these eigenvalues is $\alpha(\gamma-1)+1/2=1/2-b$. Thus $b>1/2$ is precisely this stability condition, under which the iterates attain the canonical rate $Q_n-Q^*=O_{\Pr}(n^{-1/2})$ and satisfy the CLT
\begin{equation}\label{eqn:QL-CLT}
\sqrt n(Q_n-Q^*)\Ra N(0,\Sigma),
\end{equation}
where the covariance $\Sigma$ is specified in \eqref{eqn:covariance}; see \citet{devraj2021qlearningvariance}. The analysis in this paper goes beyond this first-order limit to identify the contribution of maximization and the Edgeworth expansion of this CLT.

Throughout the rest of the paper, Assumptions~\ref{assump:unique_optimal_action} and~\ref{assump:stepsize} are in force.
\subsection{Decomposition of the Last Iterate Error}

To isolate the displacement generated by maximization, initialize $\widetilde Q_0=Q_0$ and define the policy-correct iterate by
\begin{equation}
\widetilde Q_{k+1}=(1-\alpha_{k+1})\widetilde Q_k+
\alpha_{k+1}\sqbk{R_{k+1}+\gamma\widehat P_{k+1}\Pi^*\widetilde Q_k}.
\label{eqn:policy-update}
\end{equation}
This is precisely the synchronous TD-learning recursion for the optimal policy $\pi^*$.

We separate the policy-correct error into its centered fluctuation $Y_k$ and a deterministic mean $m_n$, and use $D_k$ to record the error created by maximization instead of directly implementing $\pi^*$:
\[
Y_k:=\widetilde Q_k-\E\widetilde Q_k,\qquad
m_k:=\E\widetilde Q_k-Q^*,\qquad
D_k:=Q_k-\widetilde Q_k.
\]

To study the asymptotics of the last iterate error, we first write down the recursion that $(Y,m,D)$ satisfies. We center the one-step Bellman innovation at $Q^*$ and define
\[
\xi_k:=R_k-r+\gamma(\widehat P_k-P)\Pi^* Q^*.
\]
Then, the pairs $\set{(\widehat P_k,\xi_k):k\ge1}$ are i.i.d. and bounded, and $\E\xi_k=0$. Also, define the random matrix
\[
C_k:=(1-\alpha_k)I+\gamma\alpha_k\widehat P_k\Pi^*.
\]

Note that each $C_k$ is nonnegative with the deterministic row sum $c_k= 1-(1-\gamma)\alpha_k$. Thus, 
\begin{equation}
\begin{aligned}
\norm{ C_n\ds C_{k+1}}_\infty
&\le\prod_{j=k+1}^n c_j = \prod_{j=k+1}^n \crbk{1-\frac{b}{j+k_0}}\\
&=\exp\set{\sum_{j=k+1}^n\log\crbk{1-\frac{b}{j+k_0}}}\\
&\leq  \exp\set{-\int_{k+1}^n\frac{b}{t+k_0}dt}\\
&=\crbk{\frac{k+1+k_0}{n+k_0}}^b
\end{aligned}
\label{eqn:product-bound}
\end{equation}
where we used the inequality $\log(1-x)\le-x$. Thus, for a fixed starting index $k$, $\norm{ C_n\ds C_{k+1}}_\infty = O(n^{-b})$.

Using these notations and subtracting the fixed-point equation from \eqref{eqn:policy-update} gives
\[
\widetilde Q_{k+1}-Q^*=C_{k+1}(\widetilde Q_k-Q^*)+\alpha_{k+1}\xi_{k+1}.
\]
Then, we center the policy-correct recursion, take expectations, and subtracting \eqref{eqn:policy-update} from \eqref{eqn:q-update} to get
\begin{equation}
\begin{aligned}
Y_{k+1}&=C_{k+1}Y_k+\alpha_{k+1}\sqbk{\xi_{k+1}+\gamma(\widehat P_{k+1}-P)\Pi^* m_k},\\
m_{k+1}&=\sqbk{I-\alpha_{k+1}(I-\gamma P\Pi^*)}m_k,\\
D_{k+1}&=C_{k+1}D_k+\gamma\alpha_{k+1}\widehat P_{k+1}\crbk{V_k-\Pi^* Q_k}.
\end{aligned}
\label{eqn:YmD-recursions}
\end{equation}
We observe that the matrix $I-\alpha_{k+1}(I-\gamma P\Pi^*)$ also has row sum $c_k$; hence applying \eqref{eqn:product-bound} to the deterministic recursion for $m_k$ gives $\norm{m_k}_\infty\le C(k+1)^{-b}$, for some constant $C$ that is independent of $k$.

Importantly, the last recursion for $D_k$ is driven exactly by the maximization residual $V_k-\Pi^* Q_k$, which is nonnegative. Since $D_0=0$ and every $C_k$ is nonnegative, $D_k\ge0$ pathwise. Consequently, $D_k$ is the upward displacement produced by maximizing rather than evaluating the optimal policy. Thus, we call this term the maximization transient.  Also note that $Y_k + m_k = \widetilde Q_k-Q^*$, which is the error of the policy-correct recursion. Since $\E Y_k=0$, we refer to $Y_k$ as the centered policy-correct error and $m_n$ as policy-correct bias.

These three components decompose the normalized last-iterate error as follows
\begin{equation}
\sqrt n(Q_n-Q^*)=\sqrt nY_n+\sqrt n m_n+\sqrt n D_n.
\label{eqn:decomposition}
\end{equation}
The analysis follows this decomposition: we first analyze the maximization transient, then derive the Edgeworth corrections to the CLT for $\sqrt nY_n$, and finally identify the deterministic mean and combine all three terms.

\section{The Maximization Transient}
\label{section:transient}

We first focus on the maximization transient term. From the recursion in \eqref{eqn:YmD-recursions}, the maximization transient $D_n$ is generated by the residuals $\set{V_k-\Pi^* Q_k:k\le n-1}$ and subsequently propagated by the matrices $\set{C_j: k+1\leq j\le n}$. In this section, we show that, under appropriate scaling, the maximization transient has an almost sure limit, as the residuals converge to zero almost surely.

To analyze the asymptotics of $D_n$, it is useful to separate the contraction magnitude and the redistribution induced by the propagation. Specifically, recall that $C_k$ is non-negative and has row sum $c_k$. Define the stochastic matrix
\begin{equation}
\overline C_k:=\frac{C_k}{c_k}
=I+\frac{\gamma\alpha_k}{c_k}(\widehat P_k\Pi^*-I).
\label{eqn:normalized-matrix}
\end{equation}
Also let
\[
\mrm K^*=\lim_{n\ra\infty}\frac1 n\sum_{j=0}^{n-1}(P\Pi^*)^j.
\]
As $P\Pi^*$ is a finite stochastic matrix, this limit always exists, and it is a projection onto $\ker(I-P\Pi^*)$; see \citep[Appendix A.4]{puterman2014MDP}.

\begin{theorem}[Maximization transient]
\label{thm:maximization_transient}
The random matrix product limits
$W_r:=\lim_{n\ra\infty}\overline C_n\ds\overline C_r$
exist almost surely for each $r\ge1$. Define
\begin{equation}
D:=\alpha\gamma\sum_{k=0}^\infty
\frac{\Gamma(k+k_0+1)}{\Gamma(k+k_0+2-b)}
W_{k+2}\widehat P_{k+1}\crbk{V_k-\Pi^* Q_k},
\label{eqn:D}
\end{equation}
where the infinite sum in \eqref{eqn:D} has only finitely many nonzero terms almost surely. Then, $n^bD_n\ra D$ almost surely and in $L^p$  for all $1\leq p<\infty$.
Moreover, $D$ is nonnegative, belongs to $\ker(I-P\Pi^*)$, and satisfies
\begin{equation}
\E D=\alpha\gamma\sum_{k=0}^\infty
\frac{\Gamma(k+k_0+1)}{\Gamma(k+k_0+2-b)}
\mrm K^*P\,\E\sqbk{V_k-\Pi^* Q_k}.
\label{eqn:D-mean}
\end{equation}
\end{theorem}
Here, $\Gamma$ denotes the Gamma function, a natural extension of the factorial when $b$ is not an integer. The condition $k_0+1>b$ ensures that every Gamma function argument above is positive.

Therefore, Theorem \ref{thm:maximization_transient} shows that the maximization transient $D_n$ is $O_{\Pr}(n^{-b})$. Moreover, under the common choice $b = 1$, we have
$$
\E D = \alpha\gamma \mrm K^* P \sum_{k=0}^\infty
\E\sqbk{V_k-\Pi^* Q_k},$$
which accumulates the expected error arising from taking the maximum rather than using the correct optimal policy $\pi^*$. As we will present in the next proposition, even for an MDP instance with one state and two actions (a two-armed bandit setting), $\E D$ can be inversely related to the optimality gap between the two actions, and hence can be unbounded over the class of all $\gamma$-discounted MDPs. Thus, the maximization transient can still be large for moderately large $n$.
\begin{proposition}
\label{prop:large-maximization-transient}
Consider an MDP with one state, two actions $\A=\set{1,2}$.  Let $\gamma=1/2$, $\alpha=2$, $k_0=1$, and $Q_0=0$, so that $\alpha_k=2/(k+1)$ and $b=\alpha(1-\gamma)=1$. Let $\set{\varepsilon_k:k\ge1}$ be i.i.d. Rademacher random variables and define the rewards $R_k(1)=0$, $R_k(2)=\varepsilon_k-\delta,$
where $0<\delta\le1/20$. Then, $\E D=d_\delta e$ where $d_\delta\ge\frac{1}{10\delta}$ and $e$ is the vector of all ones. 
\end{proposition}

The proof of Proposition \ref{prop:large-maximization-transient} is deferred to Appendix \ref{section:proof:prop:large-maximization-transient}. In this MDP instance, we use randomized rewards to make the connection with the classical two-armed bandit setting transparent. The same inverse-gap phenomenon can also arise with deterministic, state-dependent rewards in a three-state, two-action MDP. Let the states be $\set{0,+,-}$. From state $0$, action $1$ transitions to $+$ or $-$ with equal probability, whereas action $2$ transitions to $+$ w.p.$(1-\delta)/2$ and to $-$ w.p.$(1+\delta)/2$. Then, it is not hard to show that the same $O(1/\delta)$ scaling of the mean maximization transient follows.

\subsection{Proof of Theorem \ref{thm:maximization_transient}}
The proof relies on two important lemmas. First, we show that the last time at which Q-learning identifies the wrong action has a sub-exponential tail. Second, to propagate the error, we show that the almost sure limit matrix product $W_r$ exists for every $r\geq 1$. Moreover, the range of $W_r$ is contained in $\ker(I-P\Pi^*)$.

We first establish the following bound on the wrong-action identification time. This implies that the sum in \eqref{eqn:D} is finite almost surely and also provides a moment bound that will be used later.

\begin{lemma}
\label{lemma:policy_identification_time}
Let $\tau=\sup\{k\ge0:V_k-\Pi^* Q_k\ne0\}$. There are constants $C,c>0$ such that $\Pr(\tau\ge m)\le Ce^{-cm}$ for all $m\ge1$.
\end{lemma}
Here we use the convention that the supremum of an empty set is $0$. By integrating the tail bound, Lemma \ref{lemma:policy_identification_time} implies that $\E\tau^p<\infty$ for every $1\leq p<\infty$. 
\begin{proof}
Let $M=\max\{\norm{ Q_0}_\infty,R_{\max}/(1-\gamma)\}$. We show by induction that $\norm{ Q_k}_\infty\le M$. The base case is trivial; for the induction step, 
\begin{equation}\label{eqn:Q_k_M_bd}\begin{aligned}
\norminf{ Q_{k+1}}&\leq(1-\alpha_{k+1})\norminf{ Q_k}+
\alpha_{k+1}\sqbk{\norminf{R_{k+1}}+\gamma\norminf{ V_k}}\\
&\leq (1-\alpha_{k+1})M + \alpha_{k+1}\sqbk{R_{\max}+\gamma M} \\
&\leq M;
\end{aligned}
\end{equation}
where the last step follows from $R_{\max} + \gamma M -M = R_{\max}-\max\set{(1-\gamma) \norminf{Q_0},R_{\max}}\le 0$.

We next establish the concentration of $Q_n$. Regard $Q_n$ as a measurable function of the independent random tuples $\set{(R_k,\widehat P_k):1\le k\le n}$: $Q_n=f_n\crbk{(R_k,\widehat P_k):1\le k\le n}.$ We verify the bounded-difference property of each coordinate of $f_n$.

Fix $1\le j\le n$, and consider two deterministic input sequences $\zeta=\set{(r_k,p_k):1\le k\le n}$ and $\zeta'=\set{(r_k,p_k):1\leq k\ne j\le n}\cup\set{(r_j',p_j')},$
which differ only at time $j$. Let $q_k$ and $q_k'$ be the iterates generated from the same initial condition $Q_0$ by $\zeta$ and $\zeta'$, respectively. Thus,
\[
q_k=(1-\alpha_k)q_{k-1}+\alpha_k\sqbk{r_k+\gamma p_kv( q_{k-1})},
\]
with the analogous recursion for $q_k'$. By the same argument as in \eqref{eqn:Q_k_M_bd}, $\norminf{q_k}\le M$, $\norminf{q_k'}\le M$ for every admissible input sequence. Here the admissible inputs range over reward vectors in $[-R_{\max},R_{\max}]^{\Z}$ and sampled transition kernels.

Since the two input sequences agree before time $j$, $q_{j-1}=q_{j-1}'.$ So, the initial difference created by changing the input at time $j$ is
\[\begin{aligned}
\norminf{q_j-q_j'}
&=\alpha_j\norminf{r_j-r_j'
+\gamma\crbk{p_j-p_j'}v(q_{j-1})}\\
&\le\alpha_j\sqbk{\norminf{r_j-r_j'}
+\gamma\norminf{p_j v(q_{j-1})}
+\gamma\norminf{p_j'v(q_{j-1})}}\\
&\le2\alpha_j\crbk{R_{\max}+\gamma M}.
\end{aligned}\]
For every $k>j$, the two recursions use the same input $(r_k,p_k)$. Since $q\ra v(q)$ is $1$-Lipschitz in the sup norm and $p_k$ is stochastic,
\[\begin{aligned}
\norminf{q_k-q_k'}
&\le(1-\alpha_k)\norminf{q_{k-1}-q_{k-1}'}
+\gamma\alpha_k\norminf{v( q_{k-1})-v( q_{k-1}')}\\
&\le c_k\norminf{q_{k-1}-q_{k-1}'},
\end{aligned}\]
where $c_k=1-(1-\gamma)\alpha_k$ is the row sum. Iterating this inequality gives
\[
\norminf{f_n(\zeta)-f_n(\zeta')}
=\norminf{q_n-q_n'} \leq \sigma_{n,j}
:=2\crbk{R_{\max}+\gamma M}\alpha_j\prod_{k=j+1}^n c_k,
\]

Thus, for every state-action pair $z$,
\[
\sup_{\substack{\zeta,\zeta'}}
\abs{f_n(\zeta)(z)-f_n(\zeta')(z)}
\le \sigma_{n,j}
\]
where the sup is taken over admissible sequences $\zeta$ and $\zeta'$ that differ at only time $j$. This shows the bounded-difference property needed for McDiarmid's inequality.

Using \eqref{eqn:product-bound}, we have that
\[
\sigma_{n,j}
\le \frac{c}{j+k_0}\crbk{\frac{j+1+k_0}{n+k_0}}^b
\le c(n+k_0)^{-b}(j+k_0)^{b-1},
\]
where $c$ absorbs the constants that are independent of $n$ and $j$, and can change value between inequalities. Since $b>1/2$, we have $2b-2>-1$, and therefore
\[
\sum_{j=1}^n\sigma_{n,j}^2
\le c(n+k_0)^{-2b}\sum_{j=1}^n(j+k_0)^{2b-2}
\le \frac{c}{n}.
\]
Applying McDiarmid's inequality to the scalar function $f_n(\cdot)(z)$ gives, for every $u>0$,
\[
\Pr\crbk{\abs{Q_n(z)-\E Q_n(z)}>u}
\le 2\exp\set{-\frac{2u^2}{\sum_{j=1}^n\sigma_{n,j}^2}}
\le 2e^{-cnu^2}.
\]
Taking a union bound over the state-action pairs yields
\begin{equation}\label{eqn:Q_n_MD_bd}
\Pr\crbk{\norminf{Q_n-\E Q_n}>u}
\le 2|\Z|e^{-c nu^2}.
\end{equation}

It remains to translate this concentration bound into optimal policy identification. Since $Q_n\ra Q^*$ almost surely and $\norminf{Q_n}\le M$, by bounded convergence $\E Q_n\ra Q^*.$ Hence there exists $n_0$ such that for all $ n\ge n_0$
$\norminf{\E Q_n-Q^*}<\delta^*/{4}$. We claim that for all $n\ge n_0$, 
\begin{equation}\label{eqn:tu_wrong_policy_event}
    \set{V_n-\Pi^*Q_n\neq 0}\subseteq \set{\|Q_n-\E Q_n\|_\infty\ge \frac{\delta^*}{4}}. 
\end{equation}

Indeed, on $\set{V_n-\Pi^*Q_n\neq 0}$, some wrong action must be selected at some state due to overestimation; i.e., there exists a state-action pair $(x',a')\in\Z$ such that $Q_n(x',a')> Q_n(x',\pi^*(x')).$ Since $\pi^*$ is unique, by the definition of the action gap,
\[\begin{aligned}
\delta^*
&\le Q^*(x',\pi^*(x'))-Q^*(x',a')\\
&=\sqbk{Q^*(x',\pi^*(x'))-Q_n(x',\pi^*(x'))}
+\sqbk{Q_n(x',\pi^*(x'))-Q_n(x',a')}+\sqbk{Q_n(x',a')-Q^*(x',a')}\\
&< 2\norminf{Q_n-Q^*}\\
&\leq 2(\norminf{Q_n-\E Q_n}+\norminf{\E Q_n-Q^*})\\
&\leq 2\norminf{Q_n-\E Q_n} +\frac{\delta^*}{2}.
\end{aligned}\]
Rearranging, we see that $\|Q_n-\E Q_n\|_\infty\ge \delta^*/4$.

Consequently, for any $n\ge n_0$, 
\[
\Pr\crbk{V_n-\Pi^*Q_n\ne0}
\le
\Pr\crbk{\norminf{Q_n-\E Q_n}\ge\frac{\delta^*}{4}}
\le Ce^{-c(\delta^*/4)^2n},
\]
where $C$ is another $n$-independent constant that can change from line to line. Consolidating the constants and applying the union bound, for $m\ge n_0$, we get 
\[\begin{aligned}
\Pr(\tau\ge m) 
\le\sum_{n=m}^\infty\Pr\crbk{V_n-\Pi^*Q_n\ne0}\le C\sum_{n=m}^\infty e^{-cn}\le Ce^{-cm}. 
\end{aligned}\]
After increasing $C$ to include $e^{cn_0}$, the same bound holds for every $m\ge1$. 
\end{proof}

Next, we establish the second lemma that clarifies the existence of the almost sure limit $W_r$. 

\begin{lemma}
\label{lemma:C_products_limit}
For each fixed $r\ge1$, there is a random stochastic matrix $W_r$ such that $\overline C_n\ds\overline C_r\ra W_r$
almost surely and in $L^p$ for every $1\leq p < \infty$. Moreover, $\image(W_r)\subseteq\ker(I-P\Pi^*)$ and $\E W_r=\mrm K^*.$
\end{lemma}

\begin{proof}
To prove this lemma, we consider the action of $\overline C_{n}\ds \overline C_r$ on $\R^{|\Z|}$. Fix $v\in\R^{\abs \Z}$, set $v_{r-1}=v$ and $v_n=\overline C_nv_{n-1}$. Thus, $v_n = \overline C_{n}\ds \overline C_rv$. Since $\overline C_k$ are stochastic, $\norm{ v_n}_\infty\le\norm{ v}_\infty$. Put $\beta_n=\gamma\alpha_n/c_n$. Then, by the definition of $\overline C_k$ in \eqref{eqn:normalized-matrix},
\[
v_n=v_{n-1}-\beta_n(I-P\Pi^*)v_{n-1}
+\beta_n(\widehat P_n-P)\Pi^* v_{n-1}.
\]

We decompose $v_n = \mrm K^* v_n + (I-\mrm K^*)v_n$. By \citet[Appendix A.4]{puterman2014MDP}, $(\mrm K^*)^2 = \mrm K^*$ and $\mrm K^* = P\Pi^*\mrm K^*=\mrm K^* P\Pi^*$; hence multiplying by $\mrm K^*$, the drift vanishes:
\[
\mrm K^*v_n=\mrm K^*v_{n-1}
+\beta_n\mrm K^*(\widehat P_n-P)\Pi^* v_{n-1}.
\]
The final term has conditional mean zero. Since $\mrm K^*$ is also stochastic, $\mrm K^*v_n$ is a bounded martingale and converges almost surely and in every finite $L^p$ to a random vector. This implies that $\mrm K^* \overline C_n\ds \overline C_r$ converges almost surely and in every finite $L^p$ to a limiting matrix, denoted by $W_r$. Since $\E\mrm K^*v_n = \E \mrm K^*v_{n-1}$ for all $n$, we have that by bounded convergence $\E W_r =\lim_{n\ra\infty} \E \mrm K^*  \overline C_n\ds \overline C_r = \mrm K^*$ .

We now show that the component in $\image(I-\mrm K^*)$ vanishes. Let $w_n=(I-\mrm K^*)v_n$. Since $\mrm K^* =\mrm K^* P\Pi^*$, we have $(I-\mrm K^*)(I-P\Pi^*) = I-P\Pi^*$ and
$$(I-P\Pi^*) w_n =(I-P\Pi^*)(I-\mrm K^*) v_n = (I-P\Pi^*) v_n.$$
Then, 
$$\begin{aligned}
w_n &= w_{n-1} - \beta_n(I-\mrm K^*)(I-P\Pi^*)v_{n-1} + \beta_n (I-\mrm K^*)(\widehat P_n-P)\Pi^* v_{n-1}\\
&= w_{n-1}  + \beta_n(P\Pi^*-I)w_{n-1} + \beta_n \epsilon_n
\end{aligned}$$
where $\epsilon_n\in \image(I-\mrm K^*)$ is a bounded martingale difference sequence. In addition, $(P\Pi^*-I)w_{n-1} = (I-\mrm K^*)(P\Pi^*-I)w_{n-1}$. Thus, this is the stochastic approximation recursion for the fixed point of $P\Pi^*$ in the image of $I-\mrm K^*$. 

We note that the step size   $$\gamma \alpha_n\leq \beta_n = \gamma \frac{\alpha_n}{1-(1-\gamma)\alpha_n}\le \alpha_n,$$
satisfying the Robbins-Monro condition. Since the noise sequence is bounded, by \citet{abounadi2002stochastic}, to show that $w_n\ra 0$ almost surely, it suffices to show that the ODE $\dot w = (P\Pi^*-I)w$ has $0$ as its globally asymptotically stable fixed point in $\image(I-\mrm K^*)$. 

Since $P\Pi^*$ is a stochastic matrix, the eigenvalues of $P\Pi^*-I$ have non-positive real parts, and the spectral radius of $P\Pi^*$ is 1. Let $H$ be the restriction of $P\Pi^*-I$ on $\image(I-\mrm K^*)$. Then we claim that the eigenvalues of $H$, i.e. $\set{\lambda:\lambda v = (P\Pi^*-I) v,\ 0\ne v\in\image(I-\mrm K^*)}$ must have strictly negative real parts. 

Indeed, the only eigenvalue of $P\Pi^*$ with real part equal to 1 is 1. Moreover, if $v = P\Pi^* v$ and $v\in \image(I-\mrm K^*)$ then it must be that $v = 0$. To see this, note that as $v\in  \image(I-\mrm K^*)$, $v = (I-\mrm K^*)u$ for some $u\in\R^{\abs{\Z}}$ and hence $\mrm K^* v = \mrm K^*(I-\mrm K^*)u = 0$, where we used $(\mrm K^*)^2 = \mrm K^*$. On the other hand, by the definition of $\mrm K^*$
\[
0 = \mrm K^*v=\lim_{n\ra\infty}\frac1 n\sum_{j=0}^{n-1}(P\Pi^*)^j v = v,
\]
i.e. $v=0$. So, $0$ cannot be an eigenvalue of $H$. Hence, the real parts of the eigenvalues of $H$ must be negative. 

Therefore, the ODE $\dot w = (P\Pi^*-I)w = Hw$ on $\image (I-\mrm K^*)$ has 0 as its unique global exponentially stable fixed point. Hence $w_n\ra 0$ almost surely. Applying this identity with $w_{r-1}= (I-\mrm K^*)e_j $ for $j=1,\ds,|\Z|$, we see that $(I- \mrm K^*)\overline C_n\ds\overline C_r\ra 0$ almost surely and in $L^p$ by the boundedness of the sequence. 

To conclude Lemma \ref{lemma:C_products_limit}, we have shown that
\[
\overline C_n\ds\overline C_r =   \mrm K^*\overline C_n\ds\overline C_r + (I- \mrm K^*)\overline C_n\ds\overline C_r \ra W_r \]
almost surely and in $L^p$ for any $1\le p< \infty$ with $W_r u\in \image(\mrm K^*),$ $\forall u\in\R^{|\Z|}$ and $\E W_r = \mrm K^*$. In addition, for $v\in\image(\mrm K^*)$, $v = \mrm K^*u$ for some $u\in \R^{|\Z|}$. Thus, $(I-P\Pi^*) v = (I-P\Pi^*) \mrm K^*u = 0$. This implies that $\image(W_r)\subseteq  \ker (I-P\Pi^*)$, completing the proof. 
\end{proof}

Equipped with Lemmas \ref{lemma:policy_identification_time} and \ref{lemma:C_products_limit}, we present the proof of Theorem \ref{thm:maximization_transient}. 

\begin{proof}[Proof of Theorem~\ref{thm:maximization_transient}]
By Lemma \ref{lemma:policy_identification_time}, $\tau:=\sup\{k\ge0:V_k-\Pi^* Q_k\ne0\}<\infty$ almost surely. Thus, the infinite sum in \eqref{eqn:D} has only finitely many nonzero terms almost surely. Moreover, by Lemma \ref{lemma:C_products_limit}, $W_r$ exists, and hence $D$ is well defined.

We first prove that $n^b D_n\ra D$ almost surely. Iterating the recursion for $D_k$ in \eqref{eqn:YmD-recursions} from $D_0=0$, we see that
\begin{equation}\label{eqn:nb_Dn}
n^bD_n=\gamma\sum_{k=0}^{n-1}n^b\alpha_{k+1}
\crbk{\prod_{j=k+2}^n c_j}
\overline C_n\ds\overline C_{k+2}\widehat P_{k+1}\crbk{V_k-\Pi^* Q_k}
\end{equation}
where we use the convention that an empty matrix product is the identity. The scalar product can be evaluated exactly. Recall the factorial identity for the Gamma function $\Gamma(n) = (n-1)!$, $n\in\N$. Extending to real-valued arguments, we have 
\begin{equation}
\label{eqn:cj_prod_using_gamma}
\prod_{j=k+2}^n c_j=\prod_{j=k+2}^n\frac{j+k_0-b}{j+k_0}=\frac{\Gamma(n+k_0+1-b)}{\Gamma(k+k_0+2-b)}
\frac{\Gamma(k+k_0+2)}{\Gamma(n+k_0+1)}.
\end{equation}

By Binet's second formula \citep[Section~12.33, Example~12.3.8]{Whittaker_Watson_2021}, for positive real $z$ tending to infinity,
\[
\log\Gamma(z)
=\crbk{z-\frac12}\log z-z+\frac12\log(2\pi)+O(z^{-1}).
\]
Set $x=n+k_0+1$. Applying this formula at $x-b$ and $x$ gives
\[\begin{aligned}
\log\frac{\Gamma(x-b)}{\Gamma(x)}
&=\crbk{x-b-\frac12}\log(x-b)-(x-b)
-\crbk{x-\frac12}\log x+x+O(x^{-1})\\
&=\crbk{x-b-\frac12}\sqbk{\log x+\log\crbk{1-\frac bx}}
-\crbk{x-\frac12}\log x+b+O(x^{-1})\\
&=-b\log x+\crbk{x-b-\frac12}\log\crbk{1-\frac bx}+b+O(x^{-1})\\
&=-b\log x+\crbk{x-b-\frac12}\crbk{-\frac bx+O(x^{-2})}+b+O(x^{-1})\\
&=-b\log x+O(x^{-1})\\
&=-b\log n+O(n^{-1}).
\end{aligned}\]
Exponentiating yields
\begin{equation}\label{eqn:tu_nb_gamma_ratio}
n^b\frac{\Gamma(n+k_0+1-b)}{\Gamma(n+k_0+1)}
=1+O(n^{-1}).
\end{equation}
Moreover, the identity $\Gamma(x+1)=x\Gamma(x)$ gives
\[
\alpha_{k+1}\Gamma(k+k_0+2)
=\frac{\alpha}{k+k_0+1}(k+k_0+1)\Gamma(k+k_0+1)
=\alpha\Gamma(k+k_0+1).
\]
Combining this with \eqref{eqn:tu_nb_gamma_ratio} and \eqref{eqn:cj_prod_using_gamma}, for each fixed $k$,
\begin{equation}
n^b\alpha_{k+1}\prod_{j=k+2}^n c_j
\ra
\alpha\frac{\Gamma(k+k_0+1)}{\Gamma(k+k_0+2-b)}.\label{eqn:scalar_limit}
\end{equation}

Thus, Lemma~\ref{lemma:C_products_limit} and \eqref{eqn:scalar_limit} imply 
\[
n^b\alpha_{k+1}
\crbk{\prod_{j=k+2}^n c_j}
\overline C_n\ds\overline C_{k+2}\widehat P_{k+1}\crbk{V_k-\Pi^* Q_k} \ra \alpha\frac{\Gamma(k+k_0+1)}{\Gamma(k+k_0+2-b)}W_{k+2}\widehat P_{k+1}\crbk{V_k-\Pi^* Q_k}.
\]
On the other hand, on $\set{\tau<\infty}$, which happens with probability 1, only the finitely many indices $k\le\tau$ contribute to the sum defining $D$ and $D_n$. Thus, $n^bD_n\ra D$ almost surely, as the nonzero terms in the sum converge almost surely.

For convergence in $L^p$, \eqref{eqn:product-bound} and boundedness of $V_k-\Pi^* Q_k$ give
\begin{equation}
\norm{ n^bD_n}_\infty
\le C\sum_{k=0}^{\tau}(k+1)^{b-1}
\le C(1+\tau)^b.
\label{eqn:d-domination}
\end{equation}
Lemma~\ref{lemma:policy_identification_time} thus implies that $n^bD_n$ is $L^{p+\delta}$ bounded for every $1\leq p<\infty$ and some $\delta>0$. Thus, dominated convergence implies the convergence in $L^p$. 

Nonnegativity of $D$ follows from the nonnegativity of all factors in \eqref{eqn:D}, and $D\in \ker(I-P\Pi^*)$ follows from $\image(W_{k+2})\subseteq\ker(I-P\Pi^*)$.

Finally, we compute the expectation $\E D$. The matrix $W_{k+2}$ is in $\sigma(C_r:r\ge k+2)$; hence it is independent of $\cF_{k+1}$. Also, $V_k-\Pi^* Q_k$ is $\cF_k$-measurable. Thus
\[\begin{aligned}
\E[W_{k+2}\widehat P_{k+1}(V_k-\Pi^* Q_k)]
&=\E\sqbk{\E[W_{k+2}|\cF_{k+1}]
\widehat P_{k+1}(V_k-\Pi^* Q_k)}\\
&=\mrm K^*\E[\widehat P_{k+1}(V_k-\Pi^* Q_k)]\\
&=\mrm K^*\E\sqbk{\E[\widehat P_{k+1}|\cF_k](V_k-\Pi^* Q_k)}\\
&=\mrm K^*P\,\E(V_k-\Pi^* Q_k).
\end{aligned}\]
where $\E[W_{k+2}|\cF_{k+1}] = \E W_{k+2} = \mrm K^*$ follows from independence and Lemma \ref{lemma:C_products_limit}. Therefore, by Tonelli's theorem
\[\begin{aligned}
\E D&= \E\sqbk{\alpha\gamma\sum_{k=0}^\infty
\frac{\Gamma(k+k_0+1)}{\Gamma(k+k_0+2-b)}
W_{k+2}\widehat P_{k+1}\crbk{V_k-\Pi^* Q_k}}\\
&=\alpha\gamma\sum_{k=0}^\infty
\frac{\Gamma(k+k_0+1)}{\Gamma(k+k_0+2-b)}
\mrm K^*P\,\E\sqbk{V_k-\Pi^* Q_k}.
\end{aligned}\]
This completes the proof. 
\end{proof}

\subsection{Maximization Transient under a Single Recurrent Class}
\label{section:single-recurrent-class}

In many operations applications, good policies are expected to stabilize the controlled system so that the long-run system dynamics are insensitive to the initial state; see \citet{wang2026optimalsamplecomplexitystable}. This motivates considering an optimal policy whose induced chain has a single recurrent class, while allowing transient states. The following corollary shows that, under this condition, the limiting maximization transient takes the same random value at every state-action pair. 

\begin{corollary}
\label{cor:single-recurrent-class}
Suppose that the stochastic matrix $P\Pi^*$ on $\Z$ has a single recurrent class. Then, there is a nonnegative random variable $\overline D$ such that $n^bD_n\ra D= \overline De$
almost surely and in $L^p$, for all $1\le p<\infty$.
\end{corollary}

\begin{proof}
We recall the definition of $\mrm K^*$. Since the limit converges absolutely, telescoping the sums yields
\[
\mrm K^*P\Pi^*-\mrm K^* = \lim_{n\ra\infty} \frac{1}{n}\sum_{j=0}^{n-1} (P\Pi^*)^{j+1} -(P\Pi^*)^j=\lim_{n\ra\infty}\frac{(P\Pi^*)^n-I}{n}=0,
\]
where the last equality follows from boundedness of the powers of a stochastic matrix. In particular, each row $\mrm K^*_j$ of $\mrm K^*$ satisfies $\mrm K^*_j P\Pi^* = \mrm K^*_j$. 

Since $P\Pi^*$ has only one recurrent class, there is a unique stationary distribution $\eta$, seen as a probability row vector, so that $\eta P\Pi^* = \eta$. By uniqueness, each row of $\mrm K^*$ must be $\eta$; i.e $\mrm K^*=e\eta$. 

By Theorem~\ref{thm:maximization_transient}, $D\in\ker(I-P\Pi^*)$ almost surely. Therefore, $(P\Pi^*)^jD=D$ for every $j\ge0$, and taking Ces\`aro averages yields $D=\mrm K^*D=(\eta D)e$ almost surely. This can hold if and only if $D = \overline D e$ for some $\overline D\in\R$ a.s. Finally, the almost sure and $L^p$ convergence of $n^bD_n$ to $ \overline D e$ follow directly from Theorem~\ref{thm:maximization_transient}.
\end{proof}

\section{The Centered Policy-Correct Recursion}
\label{section:coefficients}

After isolating the maximization transient, the centered policy-correct recursion $\set{Y_n:n\ge 0}$ is an affine SA. Classical SA theory suggests the CLT $\sqrt{n}Y_n\Ra N(0,\Sigma)$ under the rescaled-linear step size and unique optimal policy assumptions. However, since we have identified the maximization transient to be of order $O_P(n^{-b})$ with $b>1/2$, clarifying the interaction between the maximization transient and the policy-correct recursion requires analyzing the asymptotics of $Y_n$ beyond the canonical $n^{-1/2}$ scale. This motivates us to consider an Edgeworth expansion for the affine SA recursion $Y_n$, which provides higher-order corrections to the CLT for $\sqrt{n}Y_n$.

In this paper, we first provide a concise introduction to the Edgeworth expansion for a normalized sum of i.i.d. mean-zero random vectors $S_n := n^{-1/2}\sum_{i=1}^n\varepsilon_i$. In the i.i.d. case, obtaining an Edgeworth expansion proceeds in three steps:
\begin{enumerate}
\item Expand the log characteristic function using the cumulants $\set{\kappa_r(S_n):r\ge1}.$
\item Exponentiate the expanded log characteristic function to obtain an expansion of $\E e^{i u^\top S_n}$. The resulting polynomial expansion is expressed through Bell polynomials.
\item Convert the polynomial expansion of the characteristic function into an expansion of the density using the inverse Fourier transform. This inversion can be carried out systematically using Hermite polynomials.
\end{enumerate}
Here, the first step involves a computation specific to the distribution of $S_n$, whereas the second and third steps are standard machinery.

Building on the intuition and construction from the i.i.d. case, the key to obtaining an Edgeworth expansion for $\sqrt{n}Y_n$ is to compute its cumulants. In the affine SA setting, the cumulants of $\sqrt n Y_n$ are complicated, but they can themselves be further expanded as polynomials. The coefficients in these expansions can then be computed recursively, and hence algorithmically. We then show that plugging these coefficients into the standard Edgeworth expansion machinery yields the desired expansion of the density. Finally, we provide concrete computations for the second- and third-order cumulants.

\subsection{Cumulants and the I.I.D. Edgeworth Expansion}
\label{section:iid-edgeworth}

We begin with the cumulant notation used throughout. For scalar random variables $W_1,\ds,W_r$ with finite mixed moments through order $r$, define their joint cumulant by
\begin{equation}
\kappa(W_1,\ds,W_r)
:=
\sum_{\pi\in\mfk{P}_r}
(-1)^{|\pi|-1}(|\pi|-1)!
\prod_{B\in\pi}
\E\prod_{i\in B}W_i .
\label{eqn:cumulant_def}
\end{equation}
where $\mfk{P}_r$ is the set of all partitions of $\set{1,\ds,r}$. Equivalently, using the characteristic function,
\[
\kappa(W_1,\ds,W_r)
=
i^{-r}
\frac{\partial^r}{\partial t_1\ds\partial t_r}
\sqbk{\log \E\exp\crbk{i\sum_{j=1}^r t_jW_j}
}_{t_1=\ds=t_r=0},
\]
where the logarithm denotes the branch continuous at the origin and equal to zero there.
For a random vector $Y$, define its symmetric rank-$r$ cumulant tensor by
\[
\kappa_r(Y)[u_1,\ds,u_r]
:=
\kappa(u_1^\top Y,\ds,u_r^\top Y).
\]
In particular, $\kappa_1(Y)=\E Y$, $\kappa_1(Y)[u] = \E u^\top Y$, $\kappa_2(Y)=\cov(Y)$, and $\kappa_2(Y)[u_1,u_2] = u_1^\top \cov (Y)u_2$. If $Y$ is a scalar, then we specialize $\kappa_r(Y) = \kappa(Y,\ds  ,Y)$ where $Y$ is repeated $r$ times.  

For a rank-$r$ tensor $T$, write $T[v^{[r]}]=T[v,\ds,v]$ for evaluation on $r$ copies of $v$, and use $T[v^{[a]},w^{[b]}]$ analogously when $a+b=r$. The notation $T:S$ denotes the entrywise contraction of tensors of the same rank; i.e., for the rank $r$ case,
\[
T:S = \sum_{i_1,\ds,i_r} T_{i_1,\ds,i_r} S_{i_1,\ds,i_r}. 
\]
Throughout the rest of the discussion, we consider $v\in i\R^{|\Z|}$; i.e. $v=iu$ for some $u\in\R^{|\Z|}$. For bounded random vectors, the Taylor expansion of the log characteristic function is well defined in a neighborhood of the origin with
\[
\log\E e^{v^\top Y}
=\sum_{r\ge1}\frac{1}{r!}\kappa_r(Y)[v^{[r]}].
\]

Let $\varepsilon_1,\varepsilon_2,\ds$ be bounded, centered i.i.d. vectors with covariance $\Sigma\succ0$, and put $S_n=n^{-1/2}\sum_{i=1}^n\varepsilon_i$. Then by independence,
\[
\log\E e^{v^\top S_n}
=n\log\E e^{v^\top\varepsilon_1/\sqrt n}=\sum_{r\ge2}\frac{n^{1-r/2}}{r!}
\kappa_r(\varepsilon_1)[v^{[r]}].
\]
Thus the covariance at $r=2$ is $O(1)$ after normalization, while the rank-$r$ cumulant first enters at order $n^{-(r-2)/2}$ for $r\ge3$. The quadratic term produces the Gaussian limit, and the higher cumulants determine the successive shape corrections.

To see how these logarithmic corrections combine to form an expansion of the characteristic function, we consider the following expansion that retains terms through order $n^{-1}$. For each fixed $u\in\R^{|\Z|}$,
\[
\log\E e^{iu^\top S_n}
=-\frac12u^\top\Sigma u
+\frac{\kappa_3(\varepsilon_1)[(iu)^{[3]}]}{3!\sqrt n}
+\frac{\kappa_4(\varepsilon_1)[(iu)^{[4]}]}{4!n}
+O(n^{-3/2}).
\]
Using $e^x=1+x+x^2/2+O(x^3)$ for the part after the quadratic term yields
\[\begin{aligned}
\E e^{iu^\top S_n}
=e^{-\frac{1}{2}u^\top\Sigma u}\sqbk{1+\frac{\kappa_3(\varepsilon_1)[(iu)^{[3]}]}{3!\sqrt n}+\frac1n\crbk{
\frac{\kappa_4(\varepsilon_1)[(iu)^{[4]}]}{4!}
+\frac{\crbk{\kappa_3(\varepsilon_1)[(iu)^{[3]}]}^2}{2(3!)^2}
}
+O(n^{-3/2})}.
\end{aligned}\]
The order-$n^{-1}$ correction therefore contains both the fourth cumulant and the square of the third cumulant, where the latter arises from exponentiating the logarithmic expansion, not from a rank-four cumulant. At higher orders, the same exponentiation produces the Bell polynomial coefficients.

\subsubsection{From Characteristic Expansion to Gaussian Correction}
Let $\phi_\Sigma$ be the centered Gaussian density on $\R^{|\Z|}$, and define the rank-$r$ Hermite tensors by
\[
H_r(z;\Sigma)=(-1)^r\phi_\Sigma(z)^{-1}D^r\phi_\Sigma(z).
\]
Here, for coordinates $1\le i_1,\ds,i_r\le |\Z|$, $D_{i_1,\ds,i_r}^r f(z) := \del_{i_1}\del_{i_2}\ds \del_{i_r}f(z)$. 

We consider the following identity. Integrating by parts, we have
\[\begin{aligned}
\int_{\R^{|\Z|}} e^{iu^\top z}\phi_\Sigma(z)H_{r,i_1\ds i_r}(z;\Sigma)\,dz
&=(-1)^r\int_{\R^{|\Z|}} e^{iu^\top z}
\partial_{i_1}\ds\partial_{i_r}\phi_\Sigma(z)\,dz\\
&=(-1)^r(-iu_{i_1})\ds(-iu_{i_r})
\int_{\R^{|\Z|}} e^{iu^\top z}\phi_\Sigma(z)\,dz\\
&=(iu_{i_1})\ds(iu_{i_r})e^{-u^\top\Sigma u/2}.
\end{aligned}\]
The boundary terms vanish because Gaussian derivatives are polynomials times the Gaussian, and hence decay fast enough. Thus, for rank-$r$ tensor $T$, summing over the coordinates gives
\begin{equation}
\int_{\R^{|\Z|}} e^{iu^\top z}\phi_\Sigma(z)(T:H_r(z;\Sigma))\,dz
=T[(iu)^{[r]}]e^{-u^\top\Sigma u/2}.
\label{eqn:hermite-fourier}
\end{equation}

This identity implies that any degree-$r$ term $T_r[(iu)^{[r]}]$ in the polynomial expansion of the characteristic function will correspond to a density correction $\phi_\Sigma(z)(T_r:H_r(z;\Sigma))$. For example, define the density correction terms
\begin{equation}\label{eqn:Piid_1_2}\begin{aligned}
P_1^{\mathrm{iid}}(z)
&=\frac1{3!}\kappa_3(\varepsilon_1):H_3(z;\Sigma),\\
P_2^{\mathrm{iid}}(z)
&=\frac1{4!}\kappa_4(\varepsilon_1):H_4(z;\Sigma)
+\frac1{2(3!)^2}
\crbk{\kappa_3(\varepsilon_1)\otimes\kappa_3(\varepsilon_1)}:H_6(z;\Sigma),
\end{aligned}
\end{equation}
where we note that the squared third-cumulant term has degree 6 with $$\crbk{\kappa_3(\varepsilon_1)\otimes\kappa_3(\varepsilon_1)}[(iu)^{[6]}]
=\crbk{\kappa_3(\varepsilon_1)[(iu)^{[3]}]}^2.$$ 
Then by linearity and \eqref{eqn:hermite-fourier}, the Fourier transform of the signed density
\[
\phi_\Sigma(z)\sqbk{1+n^{-1/2}P_1^{\mathrm{iid}}(z)
+n^{-1}P_2^{\mathrm{iid}}(z)}
\]
is given by
\[\begin{aligned}
e^{-\frac{1}{2}u^\top\Sigma u}\sqbk{1+\frac{\kappa_3(\varepsilon_1)[(iu)^{[3]}]}{3!\sqrt n}+\frac1n\crbk{
\frac{\kappa_4(\varepsilon_1)[(iu)^{[4]}]}{4!}
+\frac{\crbk{\kappa_3(\varepsilon_1)[(iu)^{[3]}]}^2}{2(3!)^2}
}}.
\end{aligned}\]
This expression agrees with the coefficients of the Fourier transform $S_n$ through order $n^{-1}$. Thus, $P_1^{\mathrm{iid}}$ and $P_2^{\mathrm{iid}}$ should be the first and second order correction to the Gaussian CLT. 

We remark that, at this stage however, we have identified only the formal coefficients with no verification of uniformity, and additional issues can appear for ``rough" distributions. The rigorous statement of the Edgeworth expansion for the i.i.d. case is not presented, and we refer the interested reader to \citet[Chapters~4--5]{bhattacharya2010normal}. 

For the policy-correct recursion, the exponentiation and conversion using Hermite tensors remain unchanged. The difference lies in the cumulant calculation. Whereas the rank-$r$ cumulant of an i.i.d. normalized sum has the single scale $\kappa_r(S_n)=n^{1-r/2}\kappa_r(\varepsilon_1)$, a fixed cumulant rank of the policy-correct recursion generally contributes at several successive orders:
\[
\kappa_r(\sqrt nY_n)
=n^{1-r/2}\sum_{q=0}^p n^{-q}K_{r,q}+o\crbk{n^{1-r/2-p}}
\]
up to some power $p$. Thus the rank-$r$ cumulant contributes at the powers $n^{1-r/2-q}$ for $0\le q\le p$. Therefore, to obtain the density correction, we need to group all coefficients $K_{r,p}$ that contribute to each specific order of the characteristic function expansion and then apply the Hermite conversion.

\subsection{Computing the Cumulant Coefficients}
\label{section:conditional-cumulants}
\label{section:general-recursion}

In this section, we compute the expansion coefficients of the cumulants of $Y_n$. To achieve this, we first establish the asymptotic scale of the cumulants of $Y_n$ as $n\ra\infty$. Using the scale as an ansatz, we propose a polynomial expansion for the cumulants as a fractional power series and determine the successive coefficients by an inductive construction.

\subsubsection{The Asymptotic Scale of the Cumulants}
Since the centered policy-correct recursion $Y_n$ converges at the canonical $n^{-1/2}$ rate, this would imply a rank-$r$ moments bound of $O(n^{-r/2})$. This bound is overly conservative: we show that the rank-$r$ cumulant scales like $n^{-(r-1)}$ for all $r\ge 2$ up to some threshold $r_0$. 

We equip rank-$r$ tensors with the entrywise infinity norm $\norminf{T}=\max_{i_1,\ds,i_r}|T_{i_1,\ds, i_r}|.$ Then, we have
\begin{equation}
\label{eqn:tensor_norm_identity}|T[u_1,\ds,u_r]|\le\norminf{T}\prod_{i=1}^r\norm{u_i}_1,\qquad
\norminf{T}=\sup_{\norm{u_i}_1\le1}|T[u_1,\ds,u_r]|,
\end{equation}
where $\norm{u}_1=\sum_i|u_i|$. 

\begin{lemma}[Scale of the Cumulants]
\label{lemma:cumulant_scale}
For every integer $r\ge2$ satisfying $r(1-b)<1$, there is a constant $C<\infty$ independent of $n$ such that $\norminf{\kappa_r(Y_n)}\le C n^{-(r-1)}$ for all $n\ge1$.
\end{lemma}

\begin{remark}\label{rmk:large_stepsize}
Following the same proof structure as Lemma \ref{lemma:cumulant_scale}, one can show that when $r(1-b)>1$, the cumulants satisfy $\norminf{\kappa_r(Y_n)}\le C n^{-rb}$, which is larger than the scale in Lemma \ref{lemma:cumulant_scale}. This can occur when $1>b>1/2$ and $r$ is sufficiently large. Thus, the cumulants $\kappa_r(\sqrt nY_n) = \widetilde O(n^{r/2-rb})$ may fall outside the canonical half-integer Edgeworth expansion powers $\set{n^{-k/2}:k\ge 0}$, where the $\widetilde O$ suppresses log order terms when $r/2-rb$ is a half-integer. Similarly, an additional log term can appear when $r(1-b)=1$. While these are interesting phenomena to explore mathematically, from an algorithmic perspective, the higher-order cumulants satisfying $r(1-b)\ge 1$ are asymptotically negligible relative to the maximization transient. Indeed, when $r(1-b)\ge 1$, $\sqrt n D_n = O_{\Pr}(n^{1/2-b})$, which dominates $\kappa_r(\sqrt nY_n) = \widetilde O(n^{r/2-rb})$ for $r\ge 2$. Thus, in this paper, we restrict attention to the more algorithmically relevant regime $r(1-b)<1$.
\end{remark}

\begin{proof}[Proof of Lemma \ref{lemma:cumulant_scale}]
We prove by strong induction on the rank $r$ that $\norminf{\kappa_r(Y_n)}\le C(n+1)^{-(r-1)}$ for every $n\ge0$, which implies Lemma \ref{lemma:cumulant_scale}. Throughout the proof, the constant $C$ may depend on $r$ and increase from line to line but is independent of time and of the directions $u_1,\ds,u_r$ used below.

We first consider the base case $r=2$. The random elements $(C_{k+1},\widehat P_{k+1},\xi_{k+1})$ are independent of $Y_k$, and $\E Y_k=0$. Thus the law of total covariance applied to the recursion \eqref{eqn:YmD-recursions} gives
\[\begin{aligned}
\kappa_2(Y_{k+1}) = \cov (Y_{k+1})=\E\sqbk{C_{k+1}\kappa_2(Y_k)C_{k+1}^\top}+\alpha_{k+1}^2\cov\crbk{\xi_{k+1}+\gamma(\widehat P_{k+1}-P)\Pi^*m_k}.
\end{aligned}\]
Each entry of $C_{k+1}\kappa_2(Y_k)C_{k+1}^\top$ is a weighted sum of entries of $\kappa_2(Y_k)$, with nonnegative weights summing to $c_{k+1}^2$. The covariance in the second term is uniformly bounded in $k$ because $\widehat P_{k+1}$, $\xi_{k+1}$, and $m_k$ are uniformly bounded. Consequently, for every $k\ge0$,
\[
\norminf{\kappa_2(Y_{k+1})}\le c_{k+1}^2\norminf{\kappa_2(Y_k)}+C\alpha_{k+1}^2\le c_{k+1}^2\norminf{\kappa_2(Y_k)}+C(k+1)^{-2}.
\]
Since $Y_0=0$, we have $\norminf{\kappa_2(Y_0)}=0$. Iterating this scalar recursion up to time $n$ and applying \eqref{eqn:product-bound} yields
\[\begin{aligned}
\norminf{\kappa_2(Y_n)}
&\le C\sum_{k=1}^n k^{-2}\prod_{\ell=k+1}^n c_\ell^2\\
&\le C(n+1)^{-2b}\sum_{k=1}^n(k+1)^{2b-2}\\
&\le C(n+1)^{-1},
\end{aligned}\]
where the last inequality follows from $b > 1/2$ and the integral comparison bound \begin{equation}\label{eqn:int_comparison_bd}\sum_{k=1}^n (k+1)^p \le (n+1)^p +\int_1^{n+1}x^pdx \leq C(n+1)^{p+1}
\end{equation} for $p>-1$. This proves the base case.

For the induction step, fix an integer $r\ge3$ satisfying $r(1-b) < 1$. Then, we have $s(1-b) <1$ for every integer $2\le s<r$. Assume the strong induction hypothesis
\[
\norminf{\kappa_s(Y_k)}\le C(k+1)^{-(s-1)}
\]
for all $2\le s<r$ and $ k\ge0$. We will derive a scalar recursion for $\set{\norminf{\kappa_r(Y_k)}:k\ge 0}$ and iterate it to obtain the rank-$r$ bound for all time.

Fix a time $k\ge0$ and directions $u_1,\ds,u_r$ with $\norm{u_i}_1\le1$, and set
\[
W_i=u_i^\top Y_{k+1},\qquad \cG=\sigma(\widehat P_{k+1},\xi_{k+1}).
\]
Since $Y_k$ is independent of $\cG$ and has mean zero, for each nonempty block $B\subseteq\{1,\ds,r\}$,
\begin{equation}
c_B:=\kappa(W_i:i\in B\mid\cG)
=\begin{cases}
\alpha_{k+1}u_i^\top\sqbk{\xi_{k+1}+\gamma(\widehat P_{k+1}-P)\Pi^*m_k},&B=\{i\},\\
\kappa_{|B|}(Y_k)[C_{k+1}^\top u_i:i\in B],&|B|\ge2.
\end{cases}\label{eqn:tu_cB_def}
\end{equation}
The singleton case is the conditional mean. For $|B|\ge2$, the entire additive term, including the mean coupling, is fixed under conditioning and therefore does not affect the conditional cumulant. The law of total cumulance in \eqref{eqn:mixed-total-cumulance} gives
\begin{equation}
\kappa_r(Y_{k+1})[u_1,\ds,u_r]=\sum_{\pi\in\mfk{P}_r}\kappa(c_B:B\in\pi).
\label{eqn:tu_tot_cumulance}
\end{equation}
We separately bound the multi-block partitions, which involve only lower-rank cumulants, and the one-block partition, which propagates the rank-$r$ cumulant from time $k$.

\paragraph{One-block contribution.}
For $\pi=\{\{1,\ds,r\}\}$, the outer cumulant has one argument and is therefore its expectation. Thus, its contribution to the sum in \eqref{eqn:tu_tot_cumulance} is $$\kappa(c_{\set{1,\ds,r}}) = \E c_{\set{1,\ds,r}} = \E\kappa_r(Y_k)[C_{k+1}^\top u_1,\ds,C_{k+1}^\top u_r].$$

Recall that $C_{k+1}$ is nonnegative with row sum $c_{k+1}$. So, for every vector $u$, $\norm{C_{k+1}^\top u}_1\leq c_{k+1}\norm{u}_1.$
Therefore, recalling $\norm{u_i}_1\le 1$ for all $i$ and \eqref{eqn:tensor_norm_identity}, we apply this bound to get
\begin{equation}\label{eqn:tu_one_block_bd}
    \begin{aligned}
\kappa(c_{\set{1,\ds,r}})
&\le\E\abs{\kappa_r(Y_k)[C_{k+1}^\top u_1,\ds,C_{k+1}^\top u_r]}\\
&\le \norminf{\kappa_r(Y_k)}\E\prod_{i=1}^r\norm{C_{k+1}^\top u_i}_1\\
&\le c_{k+1}^r \norminf{\kappa_r(Y_k)}.
\end{aligned}
\end{equation}

\paragraph{Multi-block contributions.}
Fix a partition $\pi$ with $h=|\pi|\ge2$. For $B\in \pi$, we bound the conditional cumulant by bounding $\norm{c_B - \E c_B}_{L^\infty}$. 

We first consider singleton blocks $B=\{i\}$. By the formula in \eqref{eqn:tu_cB_def}, we have $\E c_B=0$. Since $(\xi_{k+1},\widehat P_{k+1},m_k)$ are bounded, we have
\[
\norm{c_B-\E c_B}_{L^\infty}=\norm{c_B}_{L^\infty}\le C\alpha_{k+1}(1+\norminf{m_k})\le C(k+1)^{-1}.
\]

For a block $B$ of size $2\le s < r$, we use the induction hypothesis at rank $s$. We write
\[
G_{k+1}=\gamma\widehat P_{k+1}\Pi^*-I \quad \text{and hence}\quad C_{k+1}=I+\alpha_{k+1}G_{k+1}.
\]
By multilinearity of the cumulant, 
\[\begin{aligned}
c_B&=\kappa_s(Y_k)[u_i+\alpha_{k+1} G_{k+1}^\top u_i:i\in B]\\
&=\kappa_s(Y_k)[u_i:i\in B] + \sum_{\varnothing\ne A\subseteq B}\alpha_{k+1}^{|A|}\kappa_s(Y_k)[G_{k+1}^\top u_i:i\in A;\ u_i:i\in B\setminus A]
\end{aligned}\]
Notice that the first term is deterministic. Subtracting the expectation cancels out this term.  Then, taking the $L^\infty$ norm, we have
\[\begin{aligned}
\norm{c_B-\E c_B}_{L^\infty} &\leq  2\sum_{\varnothing\ne A\subseteq B}\alpha_{k+1}^{|A|}\norm{\kappa_s(Y_k)[G_{k+1}^\top u_i:i\in A;\ u_i:i\in B\setminus A]}_{L^\infty}\\
&\stackrel{(i)}{\le} C \norm{\kappa_s(Y_k)}_\infty\sum_{\varnothing\ne A\subseteq B}\alpha_{k+1}^{|A|}\norm{\prod_{i\in A}\norm{G_{k+1}^\top u_i}_1}_{L^\infty}\\
&\stackrel{(ii)}{\le} (1+\gamma)^s C\norm{\kappa_s(Y_k)}_\infty \sum_{j=1}^s\binom{s}{j}\alpha_{k+1}^j\\
&\stackrel{(iii)}{\le} C(k+1)^{-(s-1)}\cd  s \cd s!\cd  \alpha_{k+1}\\
&\leq C(k+1)^{-s},
\end{aligned}\]
where $(i)$ follows from \eqref{eqn:tensor_norm_identity} and $\norm{u_i}\le 1$, $(ii)$ follows from the absolute row sum of $G_{k+1}$ is at most $1+\gamma$ almost surely, and $(iii)$ applies the induction hypothesis.

Combining the $|B| = 1$ and $|B| \ge 2$ cases, we have that $\norm{c_B-\E c_B}_{L^\infty} \le C(k+1)^{-|B|}.$ Thus, by shift invariance and the product bound \eqref{eqn:cumulant-product-bound} in Lemma~\ref{lemma:cumulant_identities}, we conclude that for the multi-block term, 
\begin{equation}\label{eqn:tu_cond_cumulant_bd}
|\kappa(c_B:B\in\pi)|=|\kappa(c_B-\E c_B:B\in\pi)|\le C \prod_{B\in\pi}(k+1)^{-|B|}=C(k+1)^{-r}.
\end{equation}

\paragraph{Propagation to time $n$.}
Combining \eqref{eqn:tu_cond_cumulant_bd} and \eqref{eqn:tu_one_block_bd} with \eqref{eqn:tu_tot_cumulance}, and applying \eqref{eqn:tensor_norm_identity}, we see that
\[
\norminf{\kappa_r(Y_{k+1})}= \sup_{\norm {u_i}_1\le 1} \kappa_r(Y_{k+1})[u_1,\ds,u_r]\leq c_{k+1}^r \norminf{\kappa_r(Y_k)} + C(k+1)^{-r}
\]
for all $k\ge 0$, where we recall that $C$ is independent of $u_1,\ds,u_r$ and $k$. This provides a scalar recursive relationship of $\norminf{\kappa_r(Y_{k})}$ in $k$.

Since $\norminf{\kappa_r(Y_0)}=0$, iterating up to time $n$ and applying \eqref{eqn:product-bound} gives for all $n\ge 0$
\[\begin{aligned}
\norminf{\kappa_r(Y_n)}
&\le C\sum_{k=1}^n k^{-r}\prod_{\ell=k+1}^n c_\ell^r\\
&\le C\sum_{k=1}^n k^{-r}\crbk{\frac{k+1+k_0}{n+k_0}}^{rb}\\
&\stackrel{(i)}{\le} C(n+1)^{-rb}\sum_{k=1}^n(k+1)^{rb-r}\\
&\le C(n+1)^{-rb}(n+1)^{rb-r+1}\\
&=C(n+1)^{-(r-1)}
\end{aligned}\]
where $(i)$ applies the rough estimates $k^{-r}\leq 2^r (k+1)^{-r}$ and $(k+1+k_0)^{rb}\leq (k+1)^{rb}(2+k_0)^{rb}$, and the bound on the power sum uses the assumption that $r(1-b)<1$ and integral comparison bound \eqref{eqn:int_comparison_bd}. This completes the induction on $r$. 
\end{proof}

\subsubsection{Recursive Coefficient Computation}
Using the scale in Lemma \ref{lemma:cumulant_scale} as an ansatz, for each rank $r\ge2$, we formally express $\kappa_r(Y_n)$ as an expansion involving symmetric tensors $K_{r,p}$ of the form
\begin{equation}
\kappa_r(Y_n)=n^{-(r-1)}\sum_{q=0}^p n^{-q}K_{r,q}+o\crbk{n^{-(r-1+p)}}.
\label{eqn:desired-cumulant-expansion}
\end{equation}
Here, we do not claim that this ansatz is fully justified, as the following remark points out some caveats. Nevertheless, it provides useful guidance for the subsequent construction of the Edgeworth expansion, which we will rigorously establish.
\begin{remark}\label{rmk:ansatz_p_range}
One cannot expect \eqref{eqn:desired-cumulant-expansion} to hold for arbitrary $p$, since behavior similar to that in Remark \ref{rmk:large_stepsize} can occur, and the expansion is limited by the deterministic term $m_n = O(n^{-b})$. Thus, we require $p<b$. A second consideration is that the construction and computation of $K_{r,p}$ require $r-1+p <rb$ to ensure stability. In this paper, we restrict attention to $r-1+2p\le2b$, which implies $p\leq \min\set{b,rb-r+1}$, and construct the coefficients inductively. The main reason for considering this reduced range is that the distributional expansion in Theorem \ref{thm:combined}, which combines the maximization transient with the Edgeworth expansion of the policy-correct recursion, is itself limited by the maximization transient term of order $O(n^{1/2-b})$. Expanding only to this scale leads precisely to the condition $r-1+2p\le2b$.
\end{remark}

\newcommand{\coeff}{\mathrm{Coeff}}
In this section, we provide a formal construction of the coefficient recursion, and justify it rigorously in Proposition \ref{prop:coefficient-recursion}. For each integer $j\ge0$, define the collection of cumulant coefficients through order $n^{-j}$ by
\[
\coeff(j):=\set{K_{s,q}:s\ge2,\ q\ge0,\ s-1+q\le j,\ s-1+2q\le2b}.
\]
We recursively compute the coefficients in $\coeff(j+1)$ using those in $\coeff(j)$. Specifically, assuming that the coefficients in $\coeff(j)$ satisfy the induction hypothesis stated below, we uniquely identify every new coefficient $K_{r,p}$ with $r-1+p=j+1$ so that the induction hypothesis is maintained. 

For the induction hypothesis, we look at the log  characteristic function of $Y_n$ in a neighborhood of 0. Specifically, the characteristic function $\E e^{v^\top Y_k}=1$ at $v=0$. Thus, for all $v\in i\R^{|\Z|}$ sufficiently close to the origin and using the branch of the complex logarithm for which $\log 1=0$, $\Psi_k(v)=\log\E e^{v^\top Y_k}$ is well defined. To simplify notation, write $\varepsilon_{n+1}:=\gamma(\widehat P_{n+1}-P)\Pi^*m_n$. Since $(\widehat P_{n+1},\xi_{n+1})$ is independent of $Y_n$ and $m_n$ is deterministic, the first recursion in \eqref{eqn:YmD-recursions} gives
\[\begin{aligned}
\E\sqbkcond{e^{v^\top Y_{n+1}}}{\widehat P_{n+1},\xi_{n+1}}&=\exp\set{\alpha_{n+1}v^\top\crbk{\xi_{n+1}+\varepsilon_{n+1}}}\E e^{(C_{n+1}^\top v)^\top Y_n}\\
&=\exp\set{\alpha_{n+1}v^\top\crbk{\xi_{n+1}+\varepsilon_{n+1}}+\Psi_n(C_{n+1}^\top v)}.
\end{aligned}\]
Taking the expectation and logarithm yields
\begin{equation}
\Psi_{n+1}(v)=\log\E\exp\set{\alpha_{n+1}v^\top\crbk{\xi_{n+1}+\varepsilon_{n+1}}+\Psi_n(C_{n+1}^\top v)}.
\label{eqn:cgf-update}
\end{equation}

Since cumulants form the power series $\Psi_n(v)=\sum_{s\ge1}\kappa_s(Y_n)[v^{[s]}]/s!$, we use the available coefficients in $\coeff(j)$ to define the polynomial
\begin{equation}
\Psi_n^j(v):=\sum_{\substack{s\ge2,\ q\ge0\\s-1+q\le j\\s-1+2q\le2b}}\frac{n^{-(s-1+q)}}{s!}K_{s,q}[v^{[s]}].
\label{eqn:known-cumulant-polynomial}
\end{equation}
The superscript $j$ is the highest order of the terms in the power series approximation. Inserting it into the one-step recursion \eqref{eqn:cgf-update} and subtracting the same expansion evaluated at time $n+1$ gives the following discrepancy:
\begin{equation}
\Delta_n^j(v):=\log\E\exp\set{\alpha_{n+1}v^\top\crbk{\xi_{n+1}+\varepsilon_{n+1}}+\Psi_n^j(C_{n+1}^\top v)}-\Psi_{n+1}^j(v).\label{eqn:def_Deltanj}
\end{equation}
This is the gap at order $j$ between approximating $\Psi_{n+1}$ using the dynamic update \eqref{eqn:cgf-update} and the version simply advancing the time $n\ra n+1$. 

\paragraph{Induction hypothesis.}
Write $[t^s]f(tv)$ to extract the coefficient of $t^s$ term in the power series of $f(tv)$ at $t=0$, which selects the degree-$s$ part in $v$. The induction hypothesis we maintain requires that at stage $j$, for each fixed $v\in i\R^{|\Z|}$ in a neighborhood of the origin,
\begin{equation}
[t^s]\Delta_n^j(tv)=o\crbk{n^{-(s+q)}}
\label{eqn:coefficient-induction}
\end{equation}
for every pair $(s,q)$ such that $K_{s,q}\in\coeff(j)$. We note that this hypothesis holds trivially when $j=0$, since $\coeff(0)=\varnothing$ and $\Psi_n^0=0$.

\paragraph{Recursive coefficient matching.}

Suppose \eqref{eqn:coefficient-induction} holds at stage $j$, and fix a new pair $(r,p)$ with $r-1+p=j+1$. If $p\ge1$, the coefficient $K_{r,p-1}$ belongs to $\coeff(j)$, and the induction hypothesis gives $[t^r]\Delta_n^j(tv)=o\crbk{n^{-(j+1)}}.$
This forces all the preceding integer-power coefficients to vanish, leaving $n^{-(j+2)}$ as the next possible order. If $p=0$, the degree-$r$ term starts with order  $n^{-(j+2)}$ as $r = j+2$. Therefore, in both cases, for each fixed $v$, we should expect a symmetric tensor $F_{r,p}$ so that
\begin{equation}
[t^r]\Delta_n^j(tv)=\frac{n^{-(j+2)}}{r!}F_{r,p}[v^{[r]}]+o\crbk{n^{-(j+2)}};
\label{eqn:discrepancy-coefficient}
\end{equation}
i.e. $F_{r,p}[v^{[r]}]/r!=[t^r n^{-(j+2)}]\Delta_n^j(tv)$ using the coefficient selection notation. In addition, we note that $F_{r,p}$ is computable from $\coeff(j)$, because $\Delta_n^j(v)$ is completely determined by coefficients in $\coeff(j)$. 

Thus, our goal is to choose the coefficient $K_{r,p}$ so that the discrepancy $\Delta_n^{j+1}$ satisfies the induction hypothesis. Define the rank-$r$ drift operator by
\[
(\cL_rT)[u_1,\ds,u_r]:=\alpha\sum_{i=1}^r T[u_1,\ds,(I-\gamma P\Pi^*)^\top u_i,\ds,u_r].
\]
Then, we show that by expanding \eqref{eqn:def_Deltanj}, in order to cancel out the order $n^{-(j+2)}$ term $[t^r]\Delta_n^{j+1}(tv)$, we should choose the coefficient $K_{r,p}$ so that
\begin{equation}
\crbk{\cL_r-(j+1)I}K_{r,p}=F_{r,p}.
\label{eqn:coefficient-equation}
\end{equation}

The recursive coefficient computation is summarized below. We start from $\coeff(0)=\varnothing$, then recursively enlarge 
\begin{equation}\label{eqn:coeff_recursion}
\coeff(j)
\xrightarrow{\text{form }\Psi_n^j }
\Delta_n^j
\xrightarrow{\text{extract coefficient}}
F_{r,p}
\xrightarrow{\text{apply }(\cL_r-(j+1)I)^{-1}}
K_{r,p}
\ra
\coeff(j+1).
\end{equation}
In the following Proposition \ref{prop:coefficient-recursion}, we show that the process above can be executed rigorously. The inductive hypothesis is maintained, the coefficients $F_{r,p}$ are well defined, and a unique solution to \eqref{eqn:coefficient-equation} exists. The proof is given in Appendix~\ref{section:proof:prop:coefficient-recursion}.

To simplify notation, we consider the exponential of $t\cL_r$ as a linear operator. Observe that
\begin{equation}\label{eqn:cL_exponential}
e^{-t\cL_r}T[u_1,\ds,u_r]=T[e^{-\alpha t(I-\gamma P\Pi^*)^\top}u_1,\ds,e^{-\alpha t(I-\gamma P\Pi^*)^\top}u_r].    
\end{equation}

\begin{proposition}[Cumulant Coefficients]
\label{prop:coefficient-recursion}
The coefficient construction in \eqref{eqn:coeff_recursion} uniquely determines symmetric tensors $F_{r,p}$ and $K_{r,p}$ for all pairs $r\ge2$, $p\ge0$, and $r-1+2p\le2b$. The unique solution $K_{r,p}$ to \eqref{eqn:coefficient-equation} is given by
\begin{equation}
K_{r,p}=\int_0^\infty e^{(r-1+p)t}e^{-t\cL_r}F_{r,p}\,dt.
\label{eqn:coefficient-integral}
\end{equation}
Then, the resulting polynomials $\Psi_n^j$ in \eqref{eqn:known-cumulant-polynomial} satisfy the induction hypothesis \eqref{eqn:coefficient-induction} for every integer $j\ge0$.
\end{proposition}

\subsection{Edgeworth Expansion for the Centered Policy-Correct Recursion}
\label{section:policy_correct_edgeworth}

The coefficient matching in Proposition~\ref{prop:coefficient-recursion} constructs a cumulant expansion polynomial $\Psi_{n}^j$ for the log characteristic function of $Y_n$ that is consistent with the one-step dynamics. To obtain the Edgeworth expansion, we follow the same procedure as in the i.i.d. case: normalize the argument, exponentiate the resulting polynomial, and convert each term to a Gaussian density correction using the Hermite tensors.

Set $\Sigma=K_{2,0}$ and $J=\floor{2b-1}$. For $1\le j\le J$ and $v \in i\R^{|\Z|}$, define
\begin{equation}
\label{eqn:def_ell_j}
\ell_j(v):=\sum_{\substack{r\ge2,\ p\ge0\\r-2+2p=j}}\frac{K_{r,p}[v^{[r]}]}{r!}.
\end{equation}

Since $\Psi^j_n(v)$ approximates the log characteristic function of $Y_n$, we scale $v$ by $\sqrt{n}$ to approximate that of $\sqrt{n}Y_n$. Consider any pair $(r,p)$ such that $r-1+2p\le2b$. Since $r-2+2p$ is an integer, we have $r-2+2p\le J$ and hence $r-1+p\le J+1$. So, the coefficient $K_{r,p}$ belongs to $\coeff(J+1)$ and the corresponding contribution to $\Psi^{J+1}_n(\sqrt{ n}v)$ is 
\[
\frac{n^{-(r-1+p)}}{r!}K_{r,p}[(\sqrt n\,v)^{[r]}]=n^{-(r-2+2p)/2}\frac{K_{r,p}[v^{[r]}]}{r!}. 
\]
Therefore, grouping the contributions by their asymptotic order in $n$, we have that
\begin{equation}
\Psi_n^{J+1}(\sqrt n\,v)=\frac12\Sigma[v^{[2]}]+\sum_{j=1}^J n^{-j/2}\ell_j(v).
\label{eqn:Psi_J_1}
\end{equation}

Exponentiating \eqref{eqn:Psi_J_1}, we first factor out the quadratic Gaussian term:
\[
\begin{aligned}
\exp\set{\Psi_n^{J+1}(\sqrt n\,v)}
&=
e^{\frac{1}{2}\Sigma[v^{[2]}]}
\exp\set{\sum_{j=1}^J n^{-j/2}\ell_j(v)}\\
&= e^{\frac{1}{2}\Sigma[v^{[2]}]}\crbk{1+\sum_{j=1}^J n^{-j/2}q_j(v)+O\crbk{n^{-(J+1)/2}}},
\end{aligned}
\]
where, for $1\le j\le J$, 
\begin{equation}
q_j(v):=\sum_{\substack{a_1,\ds,a_j\ge0\\\sum_{k=1}^jka_k=j}}
\prod_{k=1}^j\frac{\ell_k(v)^{a_k}}{a_k!}.
\label{eqn:bell}
\end{equation}
The formula for the coefficients $q_j$ is from the Bell polynomials associated with $\{\ell_k(v):k\ge 1\}$. See, e.g., \citet[Theorem A, Section~3.3]{comtet2012advanced}.

Assume $\Sigma\succ0$. From the i.i.d. setting,  the density correction term $P_j$ of order $n^{-j/2}$, should satisfy
\begin{equation}
\int_{\R^{|\Z|}}e^{v^\top z}\phi_\Sigma(z)P_j(z)\,dz=e^{\frac{1}{2}\Sigma[v^{[2]}]}q_j(v)
\label{eqn:correction-transform}
\end{equation}
for any $v\in i\R^{|\Z|}$. This can be done using the same Hermite tensor as in \eqref{eqn:hermite-fourier}. We give the example computing $P_1$ and $P_2$ as follows. 

\subsubsection{Computing the First and Second Order Corrections}
For $j=1$, from \eqref{eqn:bell} we have $q_1(v)=\ell_1(v)$. Provided that $J\ge1$, the only pair satisfying $r-2+2p=1$ is $(r,p)=(3,0)$. Thus, $q_1(v)=\ell_1(v)=K_{3,0}[v^{[3]}]/3!$. 
Applying \eqref{eqn:hermite-fourier} to 
\begin{equation}
P_1(z):=\frac1{3!}K_{3,0}:H_3(z;\Sigma).
\label{eqn:first-correction}
\end{equation}
implies \eqref{eqn:correction-transform} with $j=1$. 

Next, \eqref{eqn:bell} gives $q_2=\ell_2+\ell_1^2/2$. When $J\ge2$, the pairs contributing to $\ell_2$ are $(2,1)$ and $(4,0)$, so
\[
\ell_2(v)=\frac{K_{2,1}[v^{[2]}]}{2!}+\frac{K_{4,0}[v^{[4]}]}{4!}
\]
and hence
\[
q_2(v)=\frac{K_{2,1}[v^{[2]}]}{2!}+\frac{K_{4,0}[v^{[4]}]}{4!}+\frac{\bigl(K_{3,0}[v^{[3]}]\bigr)^2}{2(3!)^2}.
\]
Note that the last squared term has degree six in $v$. As in the i.i.d. case, we write $(K_{3,0}[v^{[3]}])^2=(K_{3,0}\otimes K_{3,0})[v^{[6]}]$. Then, letting
\[
\begin{aligned}
P_2(z)
:=\frac1{2!}K_{2,1}:H_2(z;\Sigma)+\frac1{4!}K_{4,0}:H_4(z;\Sigma)+\frac1{2(3!)^2}(K_{3,0}\otimes K_{3,0}):H_6(z;\Sigma)
\end{aligned}
\]
and applying \eqref{eqn:hermite-fourier} separately to the three terms give \eqref{eqn:correction-transform} when $j=2$. Thus the order-$n^{-1}$ correction uses Hermite tensors of ranks 2, 4, and 6, which has an additional rank 2 contribution compared to the i.i.d case \eqref{eqn:Piid_1_2}.

\subsubsection{Edgeworth Expansion under Smooth Test Functions}
Let $C_b^2(\R^{|\Z|})$ be the class of real-valued, twice continuously differentiable functions for which $f$, the gradient $\nabla f$ (seen as a row vector), and the Hessian matrix $\nabla^2 f$ are bounded. Equip this class with the norm
\[\norm{f}_{C_b^2}:=\norm{f}_\infty+\norm{\nabla f}_\infty+\norm{\nabla^2 f}_\infty\]
where the infinity norm for the gradient and the Hessian is defined first by taking the vector and matrix 2-norm and then sup over $z\in\R^{|\Z|}$.

In Theorem \ref{thm:edgeworth}, we show that the centered policy-correct recursion admits an Edgeworth expansion up to an $o(n^{1/2-b})$ uniform error over the test function class $\set{f\in C_b^2(\R^{|\Z|}):\norm{f}_{C_b^2}\le 1}$. The proof of this theorem is deferred to Appendix~\ref{section:edgeworth-proof}.

\begin{theorem}[Edgeworth expansion for the centered policy-correct recursion]
\label{thm:edgeworth}
Let $\Sigma=K_{2,0}\succ0$, $J = \floor{2b-1}$, and $K_{r,p}$ be the tensors constructed in Proposition~\ref{prop:coefficient-recursion}. Define the signed density
\[
\psi_n(z):= \phi_\Sigma(z)\sqbk{1+\sum_{j=1}^J n^{-j/2}P_j(z)}
\]
Then, $\int_{\R^{|\Z|}}\psi_n(z)dz =1$ and $\int_{\R^{|\Z|}}z\psi_n(z)dz =0\in\R^{|\Z|}$. Moreover, 
\begin{equation}
\sup_{\norm{f}_{C_b^2}\le1}\abs{\E f(\sqrt nY_n)-\int_{\R^{|\Z|}}f(z)\psi_n(z)dz}=o(n^{1/2-b}),
\label{eqn:policy-edgeworth}
\end{equation}
where $P_j$ are constructed from \eqref{eqn:bell} and Hermite tensors \eqref{eqn:hermite-fourier} such that \eqref{eqn:correction-transform} holds. 
\end{theorem}

In particular, this test function class is rich enough to imply $\sqrt nY_n\Ra N(0,\Sigma)$. Indeed, applying \eqref{eqn:policy-edgeworth} to the real and imaginary parts of $z\ra e^{iu^\top z}$ for each fixed $u\in\R^{|\Z|}$, then $\E e^{iu^\top\sqrt nY_n}\ra e^{-u^\top\Sigma u/2}.$
Thus, L\'evy's continuity theorem implies $\sqrt nY_n\Ra N(0,\Sigma)$. In addition, the expansion in Theorem \ref{thm:edgeworth} is uniform over the $C^2_b$ class and hence stronger than the convergence needed for the CLT. We also note that Theorem \ref{thm:edgeworth} doesn't require the iterate $Y_n$ to have a density, which is usually the case for MDPs with non-random reward.

If $2b-1$ is an integer, the sum in \eqref{eqn:policy-edgeworth} includes the correction of exact order $n^{1/2-b}$. Otherwise, all corrections in Theorem \ref{thm:edgeworth} have larger orders than $n^{1/2-b}$. In both cases, the remainder is smaller than $O(n^{1/2-b})$, the order of the normalized maximization transient.

\subsection[The Gaussian Correction When b=1]{The Gaussian Correction When $b=1$}
\label{section:worked-coefficients}

We now specialize to $b=1$, which is widely used in theoretical analysis of Q-learning \citep{wainwright2019cone,devraj2021qlearningvariance}. In particular, \citet[Proposition 5.2]{devraj2021qlearningvariance} show that the variance on the slowest eigen-direction is minimized by $b=1$. 

Under this choice, $J=\floor{2b-1} = 1$, and both the first centered correction and the normalized maximization transient $\sqrt n D_n$ are of order $O(n^{-1/2})$. The relevant coefficient tensors are $K_{2,0} = \Sigma$, which we will show is the usual covariance, and $K_{3,0}$, which enters $P_1$ as in \eqref{eqn:first-correction}.

\subsubsection[The covariance K20]{The Covariance $K_{2,0}= \Sigma$}
Note that $\coeff(0)=\varnothing$ and $\coeff(1)=\{K_{2,0}\}$. By \eqref{eqn:coeff_recursion}, we first form the discrepancy $\Delta_n^0$ and extract the rank-2 coefficient $F_{2,0}$. 

Recall the definition of $\Delta_n^0$ in \eqref{eqn:def_Deltanj} and $\Psi_n^0 = 0$. The rank-2 term for the expansion of $\Delta_n^0(v)$ is $1/2!$ times the second cumulant, a.k.a. the covariance, of $\alpha_{n+1}(\xi_{n+1}+\varepsilon_{n+1})$. Since $\E \xi_{n+1} = \E \varepsilon_{n+1}=0$ and $\varepsilon_{n+1}= O(n^{-b})=O(n^{-1})$, writing $\Xi=\E[\xi_1\xi_1^\top]$, we have
\[\begin{aligned}
[t^2]\Delta_n^0(tv)
&=\frac{1}{2!}\alpha_{n+1}^2\E\sqbk{\bigl(v^\top(\xi_{n+1}+\varepsilon_{n+1})\bigr)^2}=\frac{\alpha^2}{2n^2}v^\top\Xi v+O(n^{-3}).
\end{aligned}\]
Thus, $F_{2,0}=\alpha^2\Xi$ by \eqref{eqn:discrepancy-coefficient}. 

Viewing rank-two tensors as matrices, we have
\[
\cL_2K_{2,0} =\alpha\sqbk{(I-\gamma P\Pi^*)K_{2,0}+ K_{2,0}(I-\gamma P\Pi^*)^\top}.
\]
Therefore, writing $\Sigma = K_{2,0}$, the coefficient matching equation \eqref{eqn:coefficient-equation} and Proposition \ref{prop:coefficient-recursion} implies that $\Sigma$ should uniquely solve
\[
\alpha\sqbk{(I-\gamma P\Pi^*)\Sigma+\Sigma(I-\gamma P\Pi^*)^\top}-\Sigma=\alpha^2\Xi.
\]
This is precisely the Lyapunov equation that determines the asymptotic covariance matrix in \eqref{eqn:QL-CLT}, see \citet{devraj2021qlearningvariance}. Then, the integral formula in Proposition~\ref{prop:coefficient-recursion} becomes
\begin{equation}
\Sigma=\alpha^2\int_0^\infty e^t e^{-t\alpha(I-\gamma P\Pi^*)}\Xi e^{-t\alpha(I-\gamma P\Pi^*)^\top }\,dt.
\label{eqn:covariance}
\end{equation}

\subsubsection[The Skewness Correction K03]{The Skewness Correction $K_{3,0}$}
At stage $j=1$, the known polynomial is $\Psi_n^1(v)=\Sigma[v^{[2]}]/(2n)$, and the new coefficient is $K_{3,0}$. 

With $p = 0$ and $r = 3$, the decomposition \eqref{eqn:expanded-discrepancy} in the proof of Proposition \ref{prop:coefficient-recursion} shows that
\[
[t^3]\Delta_n^1(tv)=\frac{2}{2!\cd 2!}\kappa(c_1^1,c_2^1)+\frac{1}{3!}\kappa(c_1^1,c_1^1,c_1^1) = \frac{1}{2}\kappa(c_1^1,c_2^1)+\frac{1}{3!}\kappa_3(c_1^1),
\]
where
\[
c_1^1=\alpha_{n+1}v^\top(\xi_{n+1}+\varepsilon_{n+1})\quad \text{and}\quad c_2^1=\frac1n v^\top C_{n+1}\Sigma C_{n+1}^\top v.
\]

Since cumulants are shift-invariant, we expand $c_2^1$ and subtract its expectation to get
\[\begin{aligned}
c_2^1
&=\frac{v^\top\Sigma v}{n}+\frac{2\alpha_{n+1}}n v^\top(\gamma\widehat P_{n+1}\Pi^*-I)\Sigma v+\frac{\alpha_{n+1}^2}n v^\top(\gamma\widehat P_{n+1}\Pi^*-I)\Sigma(\gamma\widehat P_{n+1}\Pi^*-I)^\top v,\\
c_2^1-\E c_2^1
&=\frac{2\alpha\gamma}{n^2}v^\top(\widehat P_{n+1}-P)\Pi^*\Sigma v+O(n^{-3}).
\end{aligned}\]
Since $\E c_1^1=0$ and $c_1^1=\alpha n^{-1}v^\top\xi_{n+1}+O(n^{-2})$, the two contributions are
\[\begin{aligned}
\frac{1}{2}\kappa(c_1^1,c_2^1)
&=\frac{\alpha^2\gamma}{n^3}\E\sqbk{(v^\top\xi_{n+1})\bigl(v^\top(\widehat P_{n+1}-P)\Pi^*\Sigma v\bigr)}+O(n^{-4}),\\
\frac{1}{3!}\kappa_3(c_1^1)
&=\frac{\alpha^3}{6n^3}\kappa_3(v^\top\xi_{n+1})+O(n^{-4}).
\end{aligned}\]

By the definition of $F_{3,0}$ in \eqref{eqn:discrepancy-coefficient} as the degree-3 coefficient at order $n^{-3}$, we identify
\begin{equation}
\begin{aligned}
F_{3,0}[v^{[3]}] = \alpha^3\kappa_3(v^\top\xi_1)
+6\alpha^2\gamma\E\sqbk{(v^\top\xi_1)\bigl(v^\top(\widehat P_1-P)\Pi^*\Sigma v\bigr)}
\label{eqn:F30}
\end{aligned}
\end{equation}
where we have relabeled the sample indices as $1$ using i.i.d. Thus, $K_{3,0}$ is the unique solution to 
\[
{(\cL_3-2I)K_{3,0}}=F_{3,0},
\]
which is seen as a set of linear equations.

Therefore, by Theorem \ref{thm:edgeworth}, the first order Gaussian correction is given by the polynomial $P_1$ in \eqref{eqn:first-correction}, which leads to the signed density
\[
\psi_n(z) = \phi_\Sigma(z)\sqbk{1+\frac{1}{6\sqrt n}K_{3,0}:H_3(z;\Sigma)}.
\]
This density provides a refinement of the CLT at order $n^{-1/2}$. Recall that when $b=1$, the normalized centering and maximization transient are also of this order: $\sqrt{ n}(m_n+ D_n)= O(n^{-b+1/2}) = O(n^{-1/2})$. Since $D_n$ is random, there could, in principle, be interference between the corrections at order $n^{-1/2}$. However, in the next section, we show that $(\sqrt{n}Y_n,n^{b}D_n)$ are asymptotically independent, and hence such interference does not manifest at this scale.

\section{Distributional Expansion for Q-Learning}
\label{section:combination}

In this section, we establish the full distributional expansion for synchronous Q-learning. By Theorem \ref{thm:edgeworth}, the Edgeworth expansion of the centered policy-correct recursion provides the higher-order distributional correction to a Gaussian CLT. For the entire Q-learning recursion, we must also account for the deterministic mean $m_n$ and the maximization transient $D_n$ in \eqref{eqn:decomposition}. Both can contribute at order $n^{-b}$ before the $\sqrt{n}$ normalization. 

Crucially, we show that the joint distribution of $\sqrt nY_n$ and $n^bD_n$ converges to that of an independent pair $(Z,D)$, where $Z$ is Gaussian. We then use this result to establish that the maximization transient contributes to the Q-learning distributional expansion at order $n^{-b}$ only through its mean $\E D$.

\begin{theorem}[Distributional expansion]
\label{thm:combined}
Assume the setup in Theorem \ref{thm:edgeworth}. Then the following statements hold.

\begin{enumerate}[label=(\roman*)]
\item The bias for the policy-correct recursion satisfies
\begin{equation}
n^bm_n\ra m:=\frac{\Gamma(k_0+1)}{\Gamma(k_0+1-b)}\mrm K^*(Q_0-Q^*).
\label{eqn:mean-limit}
\end{equation}
\item The centered policy-correct error and the maximization transient jointly satisfy
\begin{equation}
(\sqrt nY_n,n^bD_n)\Ra(Z,D),\qquad 
\label{eqn:joint-limit}
\end{equation}
where $Z\sim N(0,\Sigma)$ and $Z$ is  independent of $D$.

\item Define the total bias parameter
\[
\mu:=\frac{\Gamma(k_0+1)}{\Gamma(k_0+1-b)}\mrm K^*(Q_0-Q^*)+\E D.
\]
Then
\begin{equation}
\sup_{\norm{f}_{C_b^2}\le1}\abs{\E f\crbk{\sqrt n(Q_n-Q^*)}-\int_{\R^{|\Z|}}f(z)\sqbk{\psi_n(z)+ n^{1/2-b}\mu^\top \Sigma^{-1}z\phi_\Sigma(z)}dz}=o(n^{1/2-b}).
\label{eqn:combined}
\end{equation}
\end{enumerate}
\end{theorem}

Note that $\Sigma\inv z = H_1(z;\Sigma)$ is the rank-one Hermite tensor. The corresponding term shifts the mean of the Gaussian correction. This suggests that the first-order effect of the maximization transient, viewed as a perturbation to the policy-correct dynamics, manifests only through its limiting mean when we consider the test function class $C_b^2$. The remaining distributional information of $D$ can only appear in higher-order expansions. This observation echoes the extensive literature on the role of ``maximization bias'' \citep{vanhasselt2010double,vanhasselt2016deep} in causing initial instabilities. On the other hand, by Theorem \ref{thm:edgeworth}, $\psi_n$ has mean zero.

\begin{remark}
Since we retain the expansion only through order $n^{1/2-b}$, the distribution of $D$ enters \eqref{eqn:combined} only through its mean. Indeed, for the independent pair $(Z,D)$ in \eqref{eqn:joint-limit}, Taylor's formula gives, for every $f\in C_b^2(\R^{|\Z|})$ and $\epsilon\ge0$,
\[
\abs{\E f(Z+\epsilon D)-\E f(Z+\varepsilon\E D)}\le\frac{\epsilon^2}{2}\norm{\nabla^2f}_\infty\E|D-\E D|^2,
\]
as the expected linear term vanishes by independence. Thus, setting $\epsilon =n^{1/2-b}$, we expect the higher-order contribution of $D$ to be of order $\epsilon^2=n^{1-2b}=o(n^{1/2-b})$ because $b>1/2$. 
\end{remark}

\subsection{Proof of Theorem~\ref{thm:combined}}

\begin{proof}
We prove the three statements in order.

\paragraph{Part (i).}
We recall that $\E C_{k+1} = I-\alpha_{k+1}(I-\gamma P\Pi^*)$. Thus, the recursion for $m_k$ in \eqref{eqn:YmD-recursions} is equivalent to $m_{k+1}=\E C_{k+1}\,m_k$. Independence of the matrices $C_k$ and the identity $C_k=c_k\overline C_k$ imply
\[
m_n=\E\sqbk{C_n\ds C_1}(Q_0-Q^*)=\crbk{\prod_{k=1}^n c_k}\E\sqbk{\overline C_n\ds\overline C_1}(Q_0-Q^*).
\]
By Lemma~\ref{lemma:C_products_limit}, the product $\overline C_n\ds\overline C_1$ converges in $L^1$ to $W_1$ and $\E W_1=\mrm K^*$. Its expectation therefore converges to $\mrm K^*$. By \eqref{eqn:scalar_limit} with $k=0$, the scalar factor satisfies
\[
\begin{aligned}
n^b\prod_{j=1}^n c_j &= \frac{c_1}{\alpha_{1}}\crbk{n^b\alpha_1\prod_{j=2}^n c_j }\\
&\ra \frac{(k_0+1)(k_0+1-b)}{\alpha(k_0+1)}\frac{\alpha\Gamma(k_0+1)}{\Gamma(k_0+2-b)}\\
&=\frac{\Gamma(k_0+1)}{\Gamma(k_0+1-b)}.
\end{aligned}
\]
Combining these limits proves \eqref{eqn:mean-limit}.

\paragraph{Part (ii).}
The key idea behind the independence of the limiting pair stems from the following observation. Although the maximization transient $D_n$ depends on the entire sample history $\cF_n$, it is primarily influenced by earlier samples, since $V_n-\Pi^*Q_n$ converges to zero almost surely and in $L^p$. On the other hand, the CLT for $\sqrt nY_n$ is an asymptotic result that is insensitive to the behavior of the recursion at early iterations. 

Concretely, we fix $N\ge1$ and consider a restarted centered policy-correct recursion from the 0 vector at time $N$:
\[
\begin{aligned}
Y_N^{(N)}&=0,\\
Y_{k+1}^{(N)}&=C_{k+1}Y_k^{(N)}+\alpha_{k+1}\sqbk{\xi_{k+1}+\gamma(\widehat P_{k+1}-P)\Pi^*m_k},\qquad k\ge N.
\end{aligned}
\]
Here $m_k$ is the original deterministic centering sequence. Thus $Y_n^{(N)}$ depends only on the pairs $(\widehat P_k,\xi_k)$ with $k>N$, is independent of $\cF_N$, and has mean zero. 

Subtracting its recursion from that of $Y_n$, we see that
\[
Y_n-Y_n^{(N)}=C_n\ds C_{N+1}Y_N.
\]
Thus, by \eqref{eqn:product-bound},
\begin{equation}
\E\abs{\sqrt nY_n-\sqrt nY_n^{(N)}}\le C\sqrt n\crbk{\frac{N+1+k_0}{n+k_0}}^b\E|Y_N|\le C_N n^{1/2-b}\ra0.
\label{eqn:restart-error}
\end{equation}
In particular, $\sqrt nY_n-\sqrt nY_n^{(N)}\ra 0$ in probability. On the other hand, Theorem \ref{thm:edgeworth} implies that $\sqrt nY_n\Ra Z$ for some $Z$ independent of $\cF_\infty$. Hence \eqref{eqn:restart-error} and Slutsky's lemma show that $\sqrt nY_n^{(N)}\Ra Z$ for each fixed $N$.

We note that since  $N^bD_N$ is $\cF_N$ measureable, $Y_n^{(N)}$ and $N^bD_N$ are independent. Thus, we have that $\crbk{\sqrt nY_n-\sqrt nY_n^{(N)}, \sqrt nY_n^{(N)}, N^bD_N}\Ra (0,Z, N^bD_N)$. By L\'evy's continuity theorem, for $v,w\in i\R^{|\Z|}$, 
\[
\begin{aligned}
\E \exp\set{v^\top\sqrt nY_n+w^\top N^bD_N}&= \E \exp\set{v^\top\sqbk{\sqrt nY_n-\sqrt nY_n^{(N)}} +v^\top\sqrt nY_n^{(N)} +w^\top N^bD_N}\\
&\ra e^{\Sigma[v^{[2]}]/2}\E e^{w^\top N^bD_N}.
\end{aligned}\]

Since for imaginary $v$ and real $x,y$, $|e^{v^\top x}-e^{v^\top y}|\leq \abs v \abs{x-y}$, we have that for all $k\ge 1$
\[
\abs{\exp\set{v^\top\sqrt nY_n+w^\top D} -  \exp\set{v^\top\sqrt nY_n+w^\top k^bD_k}}\leq |w|  \abs{D-k^bD_k}
\]
where we also use $\abs{e^{v^\top\sqrt nY_n} }=1$. Consequently,
\[
\begin{aligned}
&\limsup_{n\to\infty}\abs{\E e^{v^\top\sqrt nY_n+w^\top n^bD_n}-e^{\Sigma[v^{[2]}]/2}\E e^{w^\top D}} \\
&\quad \leq\limsup_{n\to\infty}\abs{\E e^{v^\top\sqrt nY_n+w^\top n^bD_n}\pm \E e^{v^\top\sqrt nY_n+w^\top D}- \E e^{v^\top\sqrt nY_n+w^\top N^bD_N}}\\
&\qquad + \limsup_{n\to\infty}
\abs{\E e^{v^\top\sqrt nY_n+w^\top N^bD_N}\pm e^{\Sigma[v^{[2]}]/2}\E e^{w^\top N^bD_N}-e^{\Sigma[v^{[2]}]/2}\E e^{w^\top D}}\\
&\quad\stackrel{(i)}{\le} \limsup_{n\ra\infty}\abs {w}\E\abs{n^bD_n-D} + 2|w|\E |N^bD_N-D|\\
&\quad =2|w|\E |N^bD_N-D|
\end{aligned}
\]
where $(i)$ uses $\Sigma[v^{[2]}] = v^\top \Sigma v \leq 0$ for imaginary $v$, and the last equality recalls the $L^1$ convergence $n^bD_n\ra D$ in Theorem~\ref{thm:maximization_transient}.  

Since $N$ is arbitrary, letting $N\to\infty$ proves that the characteristic function converges to that of $(Z,D)$. Thus, L\'evy's continuity theorem implies \eqref{eqn:joint-limit}.

\paragraph{Part (iii).}
First, we note that by Lemma~\ref{lemma:cumulant_scale} and the $L^2$ convergence in Theorem~\ref{thm:maximization_transient},
\begin{equation}
\sup_{n\ge 0}\E\sqbk{|\sqrt nY_n|^2+|n^b(m_n+D_n)|^2}<\infty.\label{eqn:tu_L2_bd_iterates}
\end{equation}
Thus, for $f\in C_b^2(\R^{|\Z|})$, the decomposition \eqref{eqn:decomposition} and the mean value theorem imply that 
\begin{equation}
\E f\crbk{\sqrt n(Q_n-Q^*)}
=\E f(\sqrt nY_n)+n^{1/2-b}\E\sqbk{\nabla f(\sqrt nY_n)n^b(m_n+D_n)}+O(n^{1-2b})\norm{\nabla^2f}_\infty.
\label{eqn:tu_mean_val_thm_fQ}   
\end{equation}

We claim that
\begin{equation}
\sup_{\norm{f}_{C_b^2}\le1}\abs{\E\sqbk{\nabla f(\sqrt nY_n)n^b(m_n+D_n)}-\E\nabla f(Z)\mu}\longrightarrow0.
\label{eqn:shift-expectation}
\end{equation}
By parts (i) and (ii),
\[
(\sqrt nY_n,n^b(m_n+D_n))\Ra(Z,\mu-\E D+D).
\]
In addition, by \eqref{eqn:tu_L2_bd_iterates}, the joint laws have uniformly integrable first moments, so the law converges in Wasserstein-1 distance; i.e. 
\[
d_n := W_1(\cL(\sqrt nY_n,n^b(m_n+D_n)), \cL(Z,\mu-\E D+D))\ra 0.
\]
For simplicity, we denote the two laws by $\cL_n$ and $\cL$, respectively. 

Moreover, for $\norm{f}_{C_b^2}\le1$, the mean value theorem gives
\[
\begin{aligned}
\abs{\nabla f(z)w-\nabla f(z')w'}&= \abs{\nabla f(z)(w-w')+(\nabla f(z)-\nabla f(z'))w'}\\
&\le |w-w'|+|w'|\min\set{|z-z'|,2}.
\end{aligned}\]

Let $\pi_n$ be an arbitrary coupling between $\cL_n$ and $\cL$. Then, by the independence of $Z$ and $D$, 
\[
\begin{aligned}
&\abs{\E\sqbk{\nabla f(\sqrt nY_n)n^b(m_n+D_n)}-\E\nabla f(Z)\mu}\\
&\quad=\abs{\E\sqbk{\nabla f(\sqrt nY_n)n^b(m_n+D_n)}-\E\sqbk{\nabla f(Z)(\mu-\E D+D\bigr)}}\\
&\quad = \abs{\int \nabla f(z)w \cL_n(d z,dw)- \int \nabla f(z)w \cL(d z,dw)}\\
&\quad = \abs{\int [\nabla f(z)w- \nabla f(z')w']d\pi_n}\\
&\quad\le\int |w-w'|\,d\pi_n+\int |w'|\min\set{|z-z'|,2}\,d\pi_n\\
&\quad\stackrel{(i)}{\le}\int |w-w'|\,d\pi_n+\crbk{\int|w'|^2\,d\pi_n}^{1/2}\crbk{\int\min\set{|z-z'|,2}^2\,d\pi_n}^{1/2}\\
&\quad\le\int|(z,w)-(z',w')|\,d\pi_n+\crbk{2\E|\mu-\E D+D|^2}^{1/2}\crbk{\int|(z,w)-(z',w')|\,d\pi_n}^{1/2},
\end{aligned}
\]
where $(i)$ applies the Cauchy-Schwarz inequality, and the last inequality uses the estimate $\min\set{c,2}^2 \leq 2c$ for $c\ge 0$. 

Therefore, taking the infimum over all couplings $\pi_n$ yields
\[
\abs{\E\sqbk{\nabla f(\sqrt nY_n)n^b(m_n+D_n)}-\E\nabla f(Z)\mu} \leq d_n+\crbk{2\E|\mu-\E D+D|^2}^{1/2}d_n^{1/2}
\]
which goes to $0$ as $n\ra\infty$. Since the right-hand side is uniform in $f$, this shows \eqref{eqn:shift-expectation}.

Finally, we apply  integration by parts on the Gaussian density to obtain
\[
\E\nabla f(Z)\mu=-\int_{\R^{|\Z|}}f(z)\nabla\phi_\Sigma(z)\mu\,dz=\int_{\R^{|\Z|}}f(z)\phi_\Sigma(z)z^\top\Sigma^{-1}\mu\,dz.
\]
Here, the boundary terms vanish because $f$ and its first derivatives are bounded and the Gaussian density decays fast enough at infinity. Combining this identity, \eqref{eqn:shift-expectation}, and \eqref{eqn:tu_mean_val_thm_fQ}, and applying Theorem~\ref{thm:edgeworth} to the $\E f(\sqrt{ n}Y_n)$ term proves \eqref{eqn:combined}.
\end{proof}

\section{Computation and Visualization: A Two-State Example}
\label{section:two-state-example}

In this section, we illustrate the distributional expansion in Theorem \ref{thm:combined} using a simple family of two-state MDPs parameterized by $0 < \delta<1/2$. For every instance in this family, the unique optimal policy induces a Markov chain with a single recurrent class. Moreover, as $\delta\da0$, the value gap between the optimal and suboptimal actions also approaches zero. Thus, distinguishing the optimal action becomes increasingly difficult when $\delta$ is small, a phenomenon similar to that observed in the example in Proposition \ref{prop:large-maximization-transient}.

Consider the following MDP instance with two states $\X=\{0,1\}$ and deterministic rewards $r(0)=0$, $r(1)=1$. At state $0$, there are two actions $a_1$ and $a_2$, with state-action pairs denoted by $z_1 = (0,a_1)$ and $z_2 = (0,a_2)$. They move to state $1$ with probabilities $1/2$ and $1/2+\delta$, respectively, where $0<\delta<1/2$; otherwise, the process remains at state $0$. The only action at state $1$, with state-action pair denoted by $z_3$, returns to state $0$ with probability one. 

We order the coordinates as $(z_1,z_2,z_3)$, and take $\gamma=1/2$, $\alpha = 2$ (so that $b=1$), $\alpha_k=2/(k+1)$ (so that $k_0=1$), and $Q_0=0$. It is not hard to see that for this MDP, the optimal action at state $0$ is $a_2$, 
\[
Q^*=\frac{2}{5+2\delta}\pmx{1+\delta\\1+2\delta\\3+2\delta}, \quad \text{and}\quad P\Pi^*=\pmx{0&1/2&1/2\\0&1/2-\delta&1/2+\delta\\0&1&0}.
\]

Let $B_1\sim\mathrm{Bernoulli}(1/2)$ and $B_2\sim\mathrm{Bernoulli}(1/2+\delta)$ be independent. Then the sampled transition matrix has the same distribution as
\[
\widehat P_k\eqd\pmx{1-B_1&B_1\\1-B_2&B_2\\1&0}.
\]
Since the rewards are deterministic and $Q^*(z_3)-Q^*(z_2)=4/(5+2\delta)$, $\xi_k$ simplifies to
\[
\xi_k=\gamma(\widehat P_k-P)\Pi^*Q^*\eqd \frac{2}{5+2\delta}\pmx{B_1-1/2\\B_2-1/2-\delta\\0}.
\]
In addition, we note that the one-sample transition matrix error can be expressed using $\xi_k$ as
\begin{equation}
\label{eqn:P-P_using_xi}
(\widehat P_k-P)\Pi^*=\frac{5+2\delta}{2}\xi_k\pmx{0&-1&1}.
\end{equation}

We first verify that $\Sigma\succ0$. Let $\Xi=\E[\xi_1\xi_1^\top]$. Since the two transition samples are independent,
\[
\Xi=\frac{4}{(5+2\delta)^2}\diag\crbk{\frac14,\frac14-\delta^2,0}.
\]
By the integral representation of the solution of the Lyapunov equation \eqref{eqn:covariance},
\[
\Sigma=4\int_0^\infty e^te^{-t(2I-P\Pi^*)}\Xi e^{-t(2 I-P\Pi^*)^\top }= 4\int_0^\infty e^{-t(\frac{3}{2}I-P\Pi^*)}\Xi e^{-t(\frac{3}{2} I-P\Pi^*)^\top }.
\]
Let $A = \frac{3}{2} I-P\Pi^*$. If $u^\top\Sigma u=0$ and $\abs u = 1$, then by continuity, we must have $\Xi^{1/2}e^{-A^\top t}u=0$ for all $t\ge0$. If so, taking $t=0$ implies that $u=e_3=(0,0,1)^\top$. However, $A^\top e_3=(0,-1,3/2)^\top$ is not proportional to $e_3$. Noting that $e^{-tA^\top} = I-tA^\top+O(t^2)$, we have $e^{-tA^\top}e_3$ cannot be proportional to $e_3$ for all sufficiently small but positive $t$. Hence, it is not possible to have $\Xi^{1/2}e^{-A^\top t}u=0$ for all $t\ge0$. Therefore, $\Sigma\succ0$.

\paragraph{The Edgeworth correction.}

Note that both $B_1-1/2$ and $B_2-1/2-\delta$ have mean 0, and $B_1-1/2$ is symmetric with skewness 0. Thus, using the independence of $B_1$ and $B_2$ and multilinearity, we have that all the cross terms disappear and
\[
\kappa_3(v^\top\xi_k)=\frac{8}{(5+2\delta)^3}v(z_2)^3\E\sqbk{(B_2-1/2-\delta)^3}
=-\frac{4\delta(1-4\delta^2)}{(5+2\delta)^3}v(z_2)^3.
\]
On the other hand, using the identity \eqref{eqn:P-P_using_xi}, we obtain
\[
\begin{aligned}
\E\sqbk{(v^\top\xi_k)\crbk{v^\top(\widehat P_k-P)\Pi^*\Sigma v}}&=\frac{5+2\delta}{2}\E\sqbk{(v^\top\xi_k)^2}\crbk{\Sigma v(z_3)-\Sigma v(z_2)}\\
&=\frac{v(z_1)^2+(1-4\delta^2)v(z_2)^2}{10+4\delta}\crbk{\Sigma v(z_3)-\Sigma v(z_2)}.
\end{aligned}
\]
where we use
\[
\E\sqbk{(v^\top\xi_k)^2}=\frac{4}{(5+2\delta)^2}\sqbk{\frac14v(z_1)^2+\crbk{\frac14-\delta^2}v(z_2)^2}.
\]

Applying \eqref{eqn:F30} with $\alpha=2$ and $\gamma=1/2$, we compute the tensor $F_{3,0}$ as
\[
\begin{aligned}
F_{3,0}[v^{[3]}]&=8\kappa_3(v^\top\xi_k)+12\E\sqbk{(v^\top\xi_k)\crbk{v^\top(\widehat P_k-P)\Pi^*\Sigma v}}\\
&=-\frac{32\delta(1-4\delta^2)}{(5+2\delta)^3}v(z_2)^3+\frac{6}{5+2\delta}\crbk{v(z_1)^2+(1-4\delta^2)v(z_2)^2}\crbk{\Sigma v(z_3)-\Sigma v(z_2)}.
\end{aligned}
\]
Then, Proposition~\ref{prop:coefficient-recursion} identifies $K_{3,0}$ as the unique solution of $(\cL_3-2I)K_{3,0}=F_{3,0}$, and the first centered Edgeworth polynomial is $P_1(z)=K_{3,0}:H_3(z;\Sigma)/6$.

\paragraph{The bias correction.}
Note that $P\Pi^*$ has $\set{z_2,z_3}$ as the unique recurrent class. So, Corollary \ref{cor:single-recurrent-class} applies and $D=\overline D e$. Moreover, the unique stationary distribution of $P\Pi^*$ is $\eta=\frac1{3+2\delta}(0,2,1+2\delta)$ and $\mrm K^* = e\eta$. 

Write $h_k=\crbk{Q_k(z_1)-Q_k(z_2)}_+$ where $(\cd)_+$ denotes the positive part. Then, the expectation formula \eqref{eqn:D-mean} and the contribution from $m$ in \eqref{eqn:mean-limit} imply 
\begin{equation}\label{eqn:tu_example_ED_m}
\E\overline D=\frac{2}{3+2\delta}\sum_{k=0}^{\infty}\E h_k,\quad \text{and}\quad m = k_0 \mrm K^*(Q_0-Q^*) = e\eta Q^* = -\frac{2+4\delta}{3+2\delta}e.
\end{equation}
Then, $\mu=(\E\overline D)e+m$.

\paragraph{One-dimensional visualization.} 
To visualize this correction, we consider the projection of the normalized error onto the vector of all ones $e$; i.e., $W_n:=\sqrt n\,e^\top(Q_n-Q^*).$ 
To obtain the distributional expansion, we apply Theorem~\ref{thm:combined} to functions of $e^\top z$. For any $g\in C_b^3(\R)$, let $f(z)=g(e^\top z)$. Then $f(\sqrt n(Q_n-Q^*))=g(W_n)$. Moreover, if $Z\sim N(0,\Sigma)$, then $W:=e^\top Z\sim N(0,\sigma_e^2)$, where $\sigma_e^2=e^\top\Sigma e$. In particular,
\[
\int_{\R^3}g(e^\top z)\phi_\Sigma(z)\,dz=\int_\R g(w)\phi_{\sigma_e^2}(w)\,dw. \]

We now consider the centered correction. Recall that $P_1(z)=K_{3,0}:H_3(z;\Sigma)/6$ and define
\[
P_1^e(w):=\frac16K_{3,0}[e^{[3]}]H_3(w;\sigma_e^2)=\frac{K_{3,0}[e^{[3]}]}{6\sigma_e^6}\crbk{w^3-3\sigma_e^2w}.
\]
For $f(z)=g(e^\top z)$, the chain rule gives
$\partial_i\partial_j\partial_k f(z)=e_i e_j e_k g'''(e^\top z)$ and $ \nabla f(z)=g'(e^\top z)e^\top.$
Integrating by parts, we obtain
\[
\begin{aligned}
\int_{\R^3}g(e^\top z)P_1(z)\phi_\Sigma(z)\,dz
&=\frac{K_{3,0}[e^{[3]}]}6\int_{\R^3}g'''(e^\top z)\phi_\Sigma(z)\,dz\\
&=\frac{K_{3,0}[e^{[3]}]}6\int_\R g'''(w)\phi_{\sigma_e^2}(w)\,dw\\
&=\int_\R g(w)P_1^e(w)\phi_{\sigma_e^2}(w)\,dw.
\end{aligned}
\]
Similarly, the bias correction can be rewritten as
\[
\int_{\R^3}g(e^\top z)\mu^\top\Sigma^{-1}z\phi_\Sigma(z)\,dz
=(e^\top\mu)\int_{\R^3}g'(e^\top z)\phi_\Sigma(z)\,dz
=\int_\R g(w)\frac{e^\top\mu}{\sigma_e^2}w\phi_{\sigma_e^2}(w)\,dw.
\]

Substituting them into Theorem~\ref{thm:combined}, we obtain the signed density approximation
\begin{equation}
\phi_{\sigma_e^2}(w)\sqbk{1+\frac1{\sqrt n}\crbk{P_1^e(w)+\frac{e^\top\mu}{\sigma_e^2}w}}
\label{eqn:two-state-scalar-expansion}
\end{equation}
with an error that is $o(n^{-1/2})$ uniformly over the scalar $C_b^3$ unit ball.

\begin{figure}[!htbp]
\centering
\includegraphics[width=\textwidth]{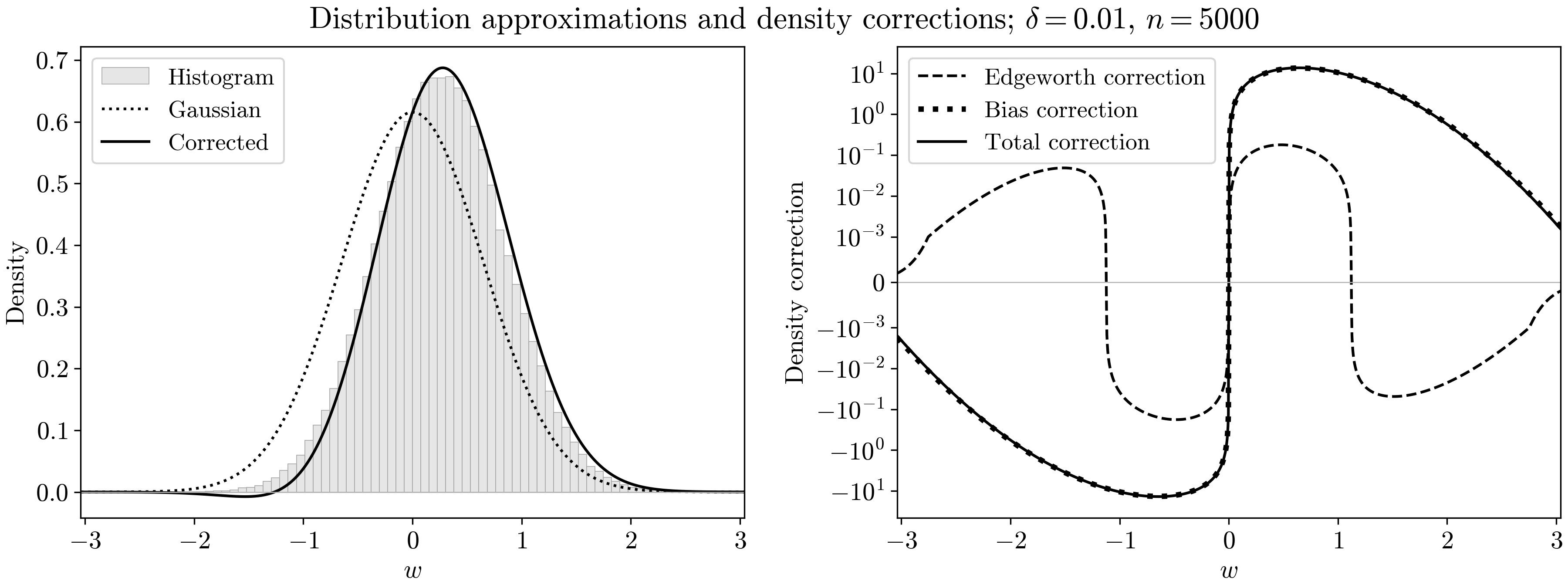}
\caption{Gaussian approximation, Edgeworth and bias corrections, and empirical histogram for $W_n$. }
\label{fig:two-state-density-correction}
\end{figure}

Figure~\ref{fig:two-state-density-correction} visualizes the density corrections for $W_n$ in Theorem \ref{thm:combined}. The left panel shows the empirical Q-learning histogram, the CLT Gaussian density, and the signed density approximation \eqref{eqn:two-state-scalar-expansion} for $W_n$. The right panel plots the two corrections: one arising from the initial and maximization bias $\phi_{\sigma_e^2}(w)(e^\top\mu)w/\sigma_e^2$, and the other from the Edgeworth expansion $\phi_{\sigma_e^2}(w)P_1^e(w)$, along with their combined effect. The vertical axis in the right panel uses a logarithmic scale.

For this simple example, one might expect $n=5000$ to be sufficiently large for the CLT to provide a good approximation. However, we still observe a noticeable deviation between the empirical distribution histogram and the Gaussian density predicted by the CLT, while the first-order corrected density provides a much better approximation. We also observe from the second plot that the $O(n^{-1/2})$ correction is predominantly from the bias correction. Additional calculations using \eqref{eqn:tu_example_ED_m} suggest that this positive bias correction arises primarily from the positive maximization transient, with $\E \overline D \approx 8.532 \pm 0.0146$, while the contribution from $m \approx-0.676$ is negative. Here, $\E \overline D$ is estimated via truncation and Monte Carlo $\E \overline D\approx \frac{1}{M}\sum_{k=1}^M N^b \overline D_{N}$ with $M = 10^7$ and $N = 10^6$. The reported range represents an approximate $95\%$ confidence interval.

\section{Conclusion}

This paper develops an asymptotic characterization of synchronous Q-learning that retains information about its initial nonlinear evolution. Theorem \ref{thm:maximization_transient} identifies a nonnegative maximization transient satisfying $n^bD_n\to D$ almost surely and in $L^p$ for every $1\le p<\infty$. Although maximization residuals stop after a finite random time, their propagated effect persists at order $n^{-b}$. Theorems \ref{thm:edgeworth} and \ref{thm:combined} combine an Edgeworth expansion of the policy-correct fluctuations with this transient through their asymptotic independence. The resulting expansion distinguishes the centered distributional corrections from shifts caused by initialization and maximization, revealing transient mechanisms that Gaussian limits discard and worst-case guarantees can obscure.

Using this characterization to develop more stable Q-learning algorithms is ongoing work. The transient structure and distributional expansion suggest ways to reduce displacement while maintaining the resulting stochastic fluctuations.

Extending the analysis beyond synchronous sampling is another important direction. We believe that, under a uniform ergodicity behavior policy and a suitable step size choice setting, we should obtain a similar theory. Temporal dependence and asynchronous updates nevertheless require substantially harder proofs and coefficient computations for the Edgeworth expansion part. This initial study therefore focuses on the synchronous setting to establish the exact characterization without obscuring its main insight with these additional difficulties.

\bibliographystyle{apalike}
\bibliography{references}

\newpage
\appendix

\section{Proof of Proposition \ref{prop:large-maximization-transient}}\label{section:proof:prop:large-maximization-transient}

\begin{proof}
Since there is only one state, the Q-function depends only on the actions. The mean rewards are $r(1)=0$ and $r(2)=-\delta$. Action $1$ is uniquely optimal, and hence $V^* = 0$,
$Q^*(1)=0$, and $Q^*(2)=-\delta.$

Let $Q_k$ be the Q-learning iterates, define
\[
\Delta_k:=Q_k(2)-Q_k(1),\qquad \widehat\Delta_k:=\Delta_k+\delta.
\]
Since there is only one state, the continuation term $\gamma V_k$ is the same in the updates for $Q_k(1)$ and $Q_k(2)$, and therefore cancels when they are subtracted. Since $\alpha_{k+1}=2/(k+2)$, we have
\[
\Delta_{k+1}=\frac{k}{k+2}\Delta_k+\frac{2}{k+2}\crbk{\varepsilon_{k+1}-\delta}.
\]
Thus, 
\[
\widehat\Delta_{k+1}=\frac{k}{k+2}\widehat\Delta_k+\frac{2}{k+2}\varepsilon_{k+1}.
\]
Since $\widehat\Delta_1=\varepsilon_1$, induction gives
\[
\widehat\Delta_k=\frac{2}{k(k+1)}\sum_{j=1}^k j\varepsilon_j.
\]
In particular, $\E \widehat \Delta_k =0$ for $k\ge1$.

Because action $1$ is optimal,
\[
V_k-\Pi^*Q_k=\crbk{\Delta_k}_+=\crbk{\widehat\Delta_k-\delta}_+.
\]
In this one-state model, $\mrm K^*Pz=ze$ for every scalar $z$. The expectation formula in Theorem~\ref{thm:maximization_transient} reduces to
\[
\E D=d_\delta e,\qquad d_\delta=\sum_{k=0}^\infty\E\crbk{\widehat\Delta_k-\delta}_+.
\]

Write $\sigma_k^2 :=\var\widehat \Delta_k=\E\sqbk{\abs{\widehat\Delta_k}^2}$. A direct calculation gives
\[
\sigma_k^2=\frac{4}{k^2(k+1)^2}\sum_{j=1}^k j^2=\frac{2(2k+1)}{3k(k+1)}\ge\frac{4}{3(k+1)}.
\]
Since $\widehat\Delta_k$ is a Rademacher sum, by symmetry and the Khintchine inequality, we have
\[
\E\crbk{\widehat\Delta_k}_+=\frac12\E\abs{\widehat\Delta_k}\ge\frac{\sigma_k}{2\sqrt2}\ge\frac{1}{\sqrt{6(k+1)}}.
\]
Using $(z-\delta)_+\ge (z)_+-\delta$, we obtain
\[
\E\crbk{\widehat\Delta_k-\delta}_+\ge\frac{1}{\sqrt{6(k+1)}}-\delta.
\]
We truncate the sum for $d_\delta$ to $N=\floor{1/(6\delta^2)}$ and obtain
\[\begin{aligned}
d_\delta &\geq\sum_{k=1}^N\E\crbk{\widehat\Delta_k-\delta}_+\ge\sum_{k=1}^N\crbk{\frac{1}{\sqrt{6(k+1)}}-\delta}\\
&\stackrel{(i)}{\ge}\frac{2}{\sqrt6}\crbk{\sqrt{N+2}-\sqrt2}-N\delta \ge\frac{1}{6\delta}-\frac{2}{\sqrt3}\\
&\ge\frac{1}{10\delta},
\end{aligned}\]
where $(i)$ follows by comparison with the integral of $x^{-1/2}$, and the last inequality uses $\delta\le1/20$. This proves the claim.
\end{proof}

\section{Cumulant Identities}
Recall the definition of cumulant in \eqref{eqn:cumulant_def}. Then, the following properties hold.

\begin{lemma}
\label{lemma:cumulant_identities}
Let $r\ge1$, and let $W_1,\ds,W_r$ be real random variables with finite $r$th absolute mixed moments. Then their joint cumulant is symmetric and multilinear; it is constant-shift invariant when $r\ge2$. Moreover, 
\begin{itemize}
\item 
If $\set{W_i:i=1,\ds,r}$ are essentially bounded, then
\begin{equation}
|\kappa(W_1,\ds,W_r)| \le\crbk{\sum_{\pi\in\mfk{P}_r}(|\pi|-1)!} \prod_{i=1}^r\norm{ W_i}_{L^\infty}.
\label{eqn:cumulant-product-bound}
\end{equation}

\item Let $\cG$ be a sub-$\sigma$-field. Then the following law of total cumulance holds:
\begin{equation}
\kappa(W_1,\ds,W_r) =\sum_{\pi\in\mfk{P}_r} \kappa\crbk{\kappa(W_i:i\in B\mid\cG):B\in\pi}.
\label{eqn:mixed-total-cumulance}
\end{equation}
\end{itemize}
Here the conditional cumulants are defined the same way as the regular cumulants in \eqref{eqn:cumulant_def}, but with $\E$ replaced by the conditional expectation. 
\end{lemma}

The product bound \eqref{eqn:cumulant-product-bound} follows directly from the definition, while the law of total cumulance in \eqref{eqn:mixed-total-cumulance} is established in \citet{brillinger1969cumulance}.

\section{Proof of Results in Section \ref{section:coefficients}}
This section presents the proofs of Proposition \ref{prop:coefficient-recursion} and Theorem \ref{thm:edgeworth}. The proof of Theorem \ref{thm:edgeworth} builds on the expansion and coefficient-matching procedure developed in the proof of Proposition \ref{prop:coefficient-recursion}.

\subsection{Proof of Proposition~\ref{prop:coefficient-recursion}}
\label{section:proof:prop:coefficient-recursion}

We first show that, given any symmetric tensor $F_{r,p}$, the coefficient equation \eqref{eqn:coefficient-equation} has a unique solution under the stated assumptions. We then inductively identify the $F_{r,p}$ using \eqref{eqn:discrepancy-coefficient}, and prove that the resulting coefficients preserve the induction hypothesis. Meanwhile, we identify the cancellations needed in the proof of Theorem~\ref{thm:edgeworth}. 

\subsubsection{Integral Representation and Uniqueness of the Solution}
For a pair $(r,p)$, the condition $r-1+2p\le2b$ gives $p\le b-(r-1)/2$. As $b>1/2$,
\begin{equation}
rb-(r-1+p)\ge(r-1)\crbk{b-\frac12}>0.
\label{eqn:coefficient-gap}
\end{equation}
Since $P\Pi^*$ is stochastic, we have $\norm{P\Pi^\top u}_{1}\leq \sum_{j}|u_j|\sum_{i}(P\Pi^*)_{j,i} = \norm{u}_1$. Thus, 
\[
\norm{e^{-\alpha t(I-\gamma P\Pi^*)^\top}u}_1\le e^{-\alpha t}\sum_{k=0}^\infty\frac{(\alpha\gamma t)^k}{k!}\norm{[(P\Pi^*)^\top]^k u}_1 \le e^{-bt}\norm{u}_1.
\]
Recall the formula for $e^{-t\cL_r}T$ in \eqref{eqn:cL_exponential}. Applying \eqref{eqn:tensor_norm_identity} to the arguments $u_1,\ds, u_r$ therefore gives
\begin{equation}
\norminf{e^{-t\cL_r}T}\le e^{-rbt}\norminf{T}.
\label{eqn:tensor-semigroup-bound}
\end{equation}

For any symmetric rank-$r$ tensor $F_{r,p}$, \eqref{eqn:coefficient-gap} and \eqref{eqn:tensor-semigroup-bound} show that the integral in \eqref{eqn:coefficient-integral} converges. Moreover, $K_{r,p}$ defined \eqref{eqn:coefficient-integral} is symmetric because the same matrix acts on every tensor argument, and it satisfies
\[
\norminf{K_{r,p}}\le\norminf{F_{r,p}}\int_0^\infty e^{-[rb-(r-1+p)]t}\,dt=\frac{\norminf{F_{r,p}}}{rb-(r-1+p)}.
\]

We observe that 
\[
-\frac{d}{dt}\crbk{e^{(r-1+p)t}e^{-t\cL_r}F_{r,p}} = (\cL_r-(r-1+p))e^{(r-1+p)t}e^{-t\cL_r}F_{r,p}.
\]
Thus, the fundamental theorem of calculus implies that
\[
(\cL_r-(r-1+p)I)K_{r,p}
=-\int_0^\infty\frac{d}{dt}\crbk{e^{(r-1+p)t}e^{-t\cL_r}F_{r,p}}\,dt
=F_{r,p},
\]
where the $+\infty$ boundary term vanishes by \eqref{eqn:tensor-semigroup-bound}. Thus, $K_{r,p}$ given in \eqref{eqn:coefficient-integral} is a solution to \eqref{eqn:coefficient-equation}.

To verify the uniqueness, we show that  $\cL_r-(r-1+p)I$ is invertible by checking that the kernel of $\cL_r-(r-1+p)I$ is $\set{0}$. If $T\in\ker(\cL_r-(r-1+p)I)$, then $e^{(r-1+p)t}e^{-t\cL_r}T=T$ for every $t\ge0$.
Thus, by \eqref{eqn:tensor-semigroup-bound} and \eqref{eqn:coefficient-gap}
\[
\norminf{e^{(r-1+p)t}e^{-t\cL_r}T} \leq e^{(r-1+p-rb)t}\norminf{T}\ra 0. 
\]
Hence $T=0$ and $\cL_r-(r-1+p)I$ is invertible. Thus, $K_{r,p}$ is the unique solution.

\subsubsection{Inductive Construction of the Coefficients and Validation}
Starting from $\coeff(0)=\varnothing$, we prove that the construction uniquely determines the coefficients at every stage and maintains the induction hypothesis \eqref{eqn:coefficient-induction} that 
\[
[t^s]\Delta_n^j(tv)=o\crbk{n^{-(s+q)}}
\]
for every $(s,q)$ such that $K_{s,q}\in \coeff(j)$. At $j=0$, there are no such pairs, so the assertion holds trivially. Suppose it holds at stage $j$, and fix a new pair $(r,p)$ with $r-1+p=j+1$ and $r-1+2p\leq 2b$ and $v\in i\R^{|\Z|}$.

\paragraph{Existence of $F_{r,p}$.}
Recall  from \eqref{eqn:def_Deltanj} that for scaler $t$, 
\begin{equation}\label{eqn:tu_Delta_tv}
\Delta_n^j(tv)=\log\E\exp\set{\underbrace{\alpha_{n+1}tv^\top\crbk{\xi_{n+1}+\varepsilon_{n+1}}+\Psi_n^j(tC_{n+1}^\top v)}_{=U_n^j(t)}}-\Psi_{n+1}^j(tv)
\end{equation}
where $\Psi_n^j$, defined in \eqref{eqn:known-cumulant-polynomial}, is formed from the coefficients in $\coeff(j)$. For fixed $v$, $\Psi_n^j(tC_{n+1}^\top v)$ is a polynomial in $t$, with no constant or linear term and degree at most $j+1$. Since $F_{r,p}$ is from the degree $r$ term, we retain terms of $U_n^j(t)$ through degree $r$:
\begin{equation}
U_n^j(t)=\sum_{s=1}^r\frac{t^s}{s!}c_s^j + O_{L^\infty}(|t|^{r+1})\label{eqn:U_polynomial}
\end{equation} where the coefficients of $t^s/s!$ are
\begin{equation}\label{eqn:tu_coeff_csj}
\begin{aligned}
c_1^j&:=\alpha_{n+1}v^\top\sqbk{\xi_{n+1}+\varepsilon_{n+1}},\\
c_s^j&:=\sum_{\substack{q\ge0\\s-1+q\le j\\s-1+2q\le2b}}n^{-(s-1+q)}K_{s,q}[(C_{n+1}^\top v)^{[s]}],\qquad 2\le s\le r.
\end{aligned}
\end{equation}

Since $(C_{n+1},\xi_{n+1},\varepsilon_{n+1})$ are bounded and ${U_n^j(t)} = O_{L^\infty}(t)$, the cumulant expansion therefore gives
\[
\log\E e^{U_n^j(t)}=\sum_{k=1}^r\frac{1}{k!}\kappa_k\crbk{U_n^j(t)}+O(|t|^{r+1}).
\]
We extract the coefficient for $t^r$. 
To obtain the power $t^r$ from the order-$k$ cumulant, we select indices $s_1,\ds,s_k$ so that $\sum_{i=1}^ks_i = r$. Then, the product of the corresponding terms in the polynomial \eqref{eqn:U_polynomial} is $O(t^r)$. Thus, by multilinearity of the cumulant

\[
[t^r]\log\E e^{U_n^j(t)}
=\sum_{k=1}^r\frac{1}{k!}[t^r]\kappa_k\crbk{U_n^j(t)}=\sum_{k=1}^r\frac{1}{k!}\sum_{\substack{s_1,\ds,s_k\ge1\\s_1+\cdots+s_k=r}}\frac{\kappa(c_{s_1}^j,\ds,c_{s_k}^j)}{s_1!\cdots s_k!}.
\]

Notice that when $k = 1$, the corresponding term in the sum is simply $\E c_r^j/r!$ with $r\ge2$. Moreover, by \eqref{eqn:tu_coeff_csj} and the choice of $(r,p)$ so that $r-1+p = j+1$  and $r-1+2p\leq 2b$, we have
\[
\E c_r^j = \sum_{q=0}^{p-1}n^{-(r-1+q)}\E K_{r,q}[(C_{n+1}^\top v)^{[r]}]
\]
Thus, separating the $k = 1$ case and subtracting $[t^r]\Psi_{n+1}^j(tv)$ yields
\begin{equation}
\begin{aligned}
[t^r]\Delta_n^j(tv)
&=\frac{1}{r!}\sum_{q=0}^{p-1}\sqbk{n^{-(r-1+q)}\E K_{r,q}[(C_{n+1}^\top v)^{[r]}]-(n+1)^{-(r-1+q)}K_{r,q}[v^{[r]}]}\\
&\quad +\sum_{k=2}^r\frac{1}{k!}\sum_{\substack{s_1,\ds,s_k\ge1\\s_1+\cdots+s_k=r}}\frac{\kappa(c_{s_1}^j,\ds,c_{s_k}^j)}{s_1!\cdots s_k!}
\end{aligned}
\label{eqn:expanded-discrepancy}
\end{equation}
where the contribution of $-[t^r]\Psi_{n+1}^j(tv)$ only comes from terms with coefficients $\set{K_{r,q}:r-1+q\leq j}$; i.e. $q=0,\ds, p-1$. 

We separately consider the two terms in \eqref{eqn:expanded-discrepancy}. Recall that $C_{n+1}=I +\alpha_{n+1}(\gamma\widehat P_{n+1}\Pi^*-I) = I+O(n^{-1})$. So, using multilinearity in the tensor arguments, we have $K_{s,q}[(C_{n+1}^\top v)^{[s]}]=K_{s,q}[v^{[s]}]+O_{L^\infty}(n^{-1}).$
Together with $(n+1)^{-(r-1+q)}=n^{-(r-1+q)}+O(n^{-(r+q)})$, this implies
\[
n^{-(r-1+q)}\E K_{r,q}[(C_{n+1}^\top v)^{[r]}]-(n+1)^{-(r-1+q)}K_{r,q}[v^{[r]}]=O(n^{-(r+q)}).
\]
Thus the first sum in \eqref{eqn:expanded-discrepancy} starts no earlier than $n^{-r}$ and hence is $O(n^{-r})$.

For the second sum, the cumulants have $k\ge2$ arguments. So, by Lemma \ref{lemma:cumulant_identities}, $\kappa$ is shift-invariant, which allows us to center the cumulant argument. Therefore, for $s\ge 2$,
\[
c_s^j =\sum_{\substack{q\ge0\\s-1+q\le j\\s-1+2q\le2b}}n^{-(s-1+q)}K_{s,q}[(C_{n+1}^\top v)^{[s]}] = \sum_{\substack{q\ge0\\s-1+q\le j\\s-1+2q\le2b}}n^{-(s-1+q)}(K_{s,q}[v^{[s]}] + O_{L^\infty}(n^{-1})).
\]
The deterministic term $K_{s,q}[v^{[s]}]$ disappears upon centering, so for $s\ge2$,
\begin{equation}
\norm{c_s^j-\E c_s^j}_{L^\infty}\le C\sum_{\substack{q\ge0\\s-1+q\le j\\s-1+2q\le2b}}n^{-(s-1+q)}n^{-1}\le Cn^{-s}.\label{eqn:tu_centered_csj_Linfty_bd}
\end{equation}

On the other hand, we notice that if the cumulant in includes $c_1^j =\alpha_{n+1}v^\top\xi_{n+1}+\alpha_{n+1}v^\top\varepsilon_{n+1}$ where, by the remark following \eqref{eqn:YmD-recursions}, $\varepsilon_{n+1} = \gamma(\widehat P_{n+1}-P)\Pi^*m_n = O(n^{-b})$ can be of non-integer order. So, we isolate this contribution by separating $c_1^j$. As $\xi_{n+1}$ is bounded, we have $\norm{\alpha_{n+1}v^\top\xi_{n+1} }_{L^\infty}=O(n^{-1})$ and $ \norm{\alpha_{n+1}v^\top\varepsilon_{n+1}}_{L^\infty} = O(n^{-{(b+1)}})$. In particular, $\norm{c_{1}^j-\E c_{1}^j}_{L^\infty} = O(n^{-1})$. 

Therefore, the cumulant product bound in Lemma \ref{lemma:cumulant_identities} gives
\begin{equation}\label{eqn:tu_cumulant_order_bd}
\abs{\kappa(c_{s_1}^j,\ds,c_{s_k}^j)}\le C\prod_{i=1}^k\norm{c_{s_i}^j-\E c_{s_i}^j}_{L^\infty}\le Cn^{-\sum_i s_i}=Cn^{-r}.
\end{equation}

Moreover, for the terms in the second sum in \eqref{eqn:expanded-discrepancy} that have $c_1^j$ in the argument of $\kappa$, we expand each $c_1^j=\alpha_{n+1}v^\top\xi_{n+1}+\alpha_{n+1}v^\top\varepsilon_{n+1}$ and use multilinearity to separate. Each separated cumulant containing $\varepsilon_{n+1} = O(n^{-b})$ gains at least one additional factor $n^{-b}$. Thus, these terms are
\[
O(n^{-(b+1)}n^{-(r-1)})= O(n^{-(r+b)})=o(n^{-(r+p)}),
\]
since $r-1+2p\leq 2b$ which implies that $p\le b-(r-1)/2<b$. So, the non-integer $b$ of the order of $\varepsilon_{n+1}$ will not contribute to the $O(n^{-(r+p)})$ term. 

Thus, after these $\varepsilon_{n+1}$-related cumulant terms are excluded, the remaining terms in $[t^r]\Delta_n^j(tv)$ have an expansion in integer powers of $n^{-1}$, as the step size and the entries of $C_{n+1}$ admit such an expansion, and the cumulants are multilinear. In addition, the above estimates for the two sums in \eqref{eqn:expanded-discrepancy} show that no term appears before order $n^{-r}$. We may therefore write
\begin{equation}\label{eqn:tu_form_of_trDelta}
[t^r]\Delta_n^j(tv)=\sum_{k = r}^{r+p} n^{-k}T_{r,k}[v^{[r]}] + o(n^{-(r+p)})
\end{equation}
for some symmetric rank-$r$ tensors $\set{T_{r,k}:k=r\ds,r+p}$,
where the $o(n^{-(r+p)})$ term collect the cumulants containing $\varepsilon_{n+1}$ and higher order integer expansion. 

To conclude the existence of $F_{r,p}$, we separately consider the following cases. If $p\ge1$, the induction hypothesis at $(r,p-1)$ gives $[t^r]\Delta_n^j(tv)=o(n^{-(r+p-1)}),$
which forces all coefficients $\set{T_{r,k}:k=r,\ds,r+p-1}$ to vanish. If $p=0$, there are no preceding coefficients and $r = r+p$. Thus, in either case, the coefficient at the first possible order $r+p$ defines $F_{r,p} := r!T_{r,r+p}$, and hence
\begin{equation}\label{eqn:tu_vanish_terms}
[t^r]\Delta_n^j(tv)=\frac{n^{-(r+p)}}{r!}F_{r,p}[v^{[r]}]+o(n^{-(r+p)}),\qquad r+p-1=j+1.
\end{equation}
This establishes \eqref{eqn:discrepancy-coefficient}.

\paragraph{Induction step.}
For each pair $(r,p)$ with $r-1+p=j+1$ and $r-1+2p\le2b$, let $K_{r,p}$ be the unique symmetric solution of \eqref{eqn:coefficient-equation} for the tensor $F_{r,p}$ determined above. By definition, 
\[
\Psi_n^{j+1}(v)=\Psi_n^j(v)+\sum_{\substack{s\ge2,\ q\ge0\\s-1+q=j+1\\s-1+2q\le2b}}\frac{n^{-(j+1)}}{s!}K_{s,q}[v^{[s]}],
\]
and form $\Delta_n^{j+1}$. 

We analyze $[t^r]\crbk{\Delta_n^{j+1}(tv)-\Delta_n^j(tv)}$ by comparing the two sums in \eqref{eqn:expanded-discrepancy} at stages $j+1$ and $j$, beginning with the second sum having $k\ge2$ cumulant arguments. 

By \eqref{eqn:tu_coeff_csj}, $c_1^{j+1}=c_1^j$, while for $s\ge2$,
\[
c_s^{j+1}-c_s^j=\sum_{\substack{q\ge0\\s-1+q=j+1\\s-1+2q\le2b}}n^{-(j+1)}K_{s,q}[(C_{n+1}^\top v)^{[s]}].
\]
Since $C_{n+1}=I+O(n^{-1})$, the leading term of each summand is the deterministic quantity $n^{-(j+1)}K_{s,q}[v^{[s]}]$. Subtracting the expectation therefore gives
\[
\norm{(c_s^{j+1}-c_s^j)-\E(c_s^{j+1}-c_s^j)}_{L^\infty}\le Cn^{-(j+1)}n^{-1}=Cn^{-(j+2)}.
\]
Moreover, by \eqref{eqn:tu_centered_csj_Linfty_bd}, for $k=j,j+1$ each $\norm{c_{s_l}^k - \E c_{s_l}^k}_{L^\infty} = O_{L^\infty}(n^{-1})$ or smaller by the same argument as in \eqref{eqn:tu_coeff_csj} and the later discussion of the special case when $s_l = 1$. 

With these identities, we use multilinearity and create a telescoping sum to obtain
\begin{align*}
\abs{\kappa(c_{s_1}^{j+1},\ds,c_{s_k}^{j+1})-\kappa(c_{s_1}^j,\ds,c_{s_k}^j)}&=\abs{\sum_{i=1}^k\kappa(c_{s_1}^{j+1},\ds,c_{s_{i-1}}^{j+1},c_{s_i}^{j+1}-c_{s_i}^j,c_{s_{i+1}}^j,\ds,c_{s_k}^j)}\\
&\leq \sum_{i=1}^k Cn^{-{(j+2)}}\crbk{\prod_{l<i} \norm{c_{s_l}^{j+1} - \E c_{s_l}^{j+1}}_{L^\infty}}\crbk{\prod_{l>i} \norm{c_{s_l}^{j} - \E c_{s_l}^{j}}_{L^\infty}}\\
&\le Cn^{-(j+2)}n^{-(k-1)}\\
&=O(n^{-(j+3)})
\end{align*}
where we used $k\ge2$ so that the cumulants are invariant to centering. 
Since the second sum in \eqref{eqn:expanded-discrepancy} involves a finite sum of the cumulants, the difference contribution to $[t^r]\crbk{\Delta_n^{j+1}(tv)-\Delta_n^j(tv)}$ is still $O(n^{-(j+3)})$.

On the other hand, for any fixed $r$ with $p=j-r+2$ (hence $r-1+p = j+1$) such that $r-1+2p\le2b$, the only additional degree-$r$ term of $\Delta_n^{j+1}(tv)$ compared to $\Delta_n^{j}(tv)$ in the first sum in \eqref{eqn:expanded-discrepancy} is the one involving $K_{r,p}$. Consequently, combining the two contributions from the two sums, we get
\begin{equation}\label{eqn:tu_tr_delta_diff}
[t^r]\crbk{\Delta_n^{j+1}(tv)-\Delta_n^j(tv)}
=\frac1{r!}\Big\{n^{-(j+1)}\E K_{r,p}[(C_{n+1}^\top v)^{[r]}]-(n+1)^{-(j+1)}K_{r,p}[v^{[r]}]\Big\}+O(n^{-(j+3)}).
\end{equation}

To analyze \eqref{eqn:tu_tr_delta_diff}, we note that by the definition of $C_{n+1}$ and multilinearity
\begin{align*}
\E K_{r,p}[(C_{n+1}^\top v)^{[r]}]
&=K_{r,p}[v^{[r]}]+r\alpha_{n+1}\E K_{r,p}[v^{[r-1]},(\gamma\widehat P_{n+1}\Pi^*-I)^\top v]+O(\alpha_{n+1}^2\norminf{K_{r,p}}\norm{v}_1^r)\\
&=K_{r,p}[v^{[r]}]-n^{-1}(\cL_rK_{r,p})[v^{[r]}]+O(n^{-2}\norminf{K_{r,p}}\norm{v}_1^r),
\end{align*}
where the last equality uses $\E\widehat P_{n+1}=P$, $\alpha_{n+1}=\alpha/n+O(n^{-2})$, and the definition of $\cL_r$. On the other hand, note that $(n+1)^{-(j+1)}=n^{-(j+1)}-(j+1)n^{-(j+2)}+O(n^{-(j+3)})$. 

Applying these expansions to \eqref{eqn:tu_tr_delta_diff} then gives
\begin{equation}
\begin{aligned}
[t^r]\crbk{\Delta_n^{j+1}(tv)-\Delta_n^j(tv)}
&=-\frac{n^{-(j+2)}}{r!}\sqbk{(\cL_r-(j+1)I)K_{r,p}}[v^{[r]}]+O(n^{-(j+3)})\\
&=-\frac{n^{-(j+2)}}{r!}F_{r,p}[v^{[r]}]+O(n^{-(j+3)}),
\end{aligned}
\label{eqn:tu_Delta_rp_diff}
\end{equation}
where the common leading term $n^{-(j+1)}K_{r,p}[v^{[r]}]$ cancels, and the second equality follows from \eqref{eqn:coefficient-equation}. Therefore, adding \eqref{eqn:discrepancy-coefficient} to \eqref{eqn:tu_Delta_rp_diff}, the order-$n^{-(j+2)}$ terms cancel, yielding
\begin{equation}
[t^r]\Delta_n^{j+1}(tv)=o(n^{-(j+2)})=o(n^{-(r+p)})\label{eqn:tu_new_rp_j+1}
\end{equation}
for every new pair $(r,p)$ so that $K_{r,p}\in \coeff(j+1)\setminus\coeff(j)$.

Finally, for any previously computed coefficient $K_{s,q}\in\coeff(j)$, the same calculation shows that the first sum changes by $O(n^{-(j+2)})$, while the second sum changes by $O(n^{-(j+3)})$. Hence
\[
[t^s]\crbk{\Delta_n^{j+1}(tv)-\Delta_n^j(tv)}=O(n^{-(j+2)})=o(n^{-(s+q)}),
\]
since $s+q\le j+1$. We then apply the induction hypothesis at $j$ to obtain 
\[
[t^s]\Delta_n^{j+1}(tv)=o\crbk{n^{-(s+q)}}. 
\]
Together with \eqref{eqn:tu_new_rp_j+1}, this establishes the induction hypothesis \eqref{eqn:coefficient-induction} for $j+1$, completing the proof of Proposition \ref{prop:coefficient-recursion}.

\subsection{Proof of Theorem~\ref{thm:edgeworth}}
\label{section:edgeworth-proof}

 Let $J=\floor{2b-1}$ and $\mu_n$ be the law of $\sqrt nY_n$. We consider signed measures $\set{\nu_n:n\ge 1}$ with signed density
\begin{equation}
\psi_n(z)=\phi_\Sigma(z)\sqbk{1+\sum_{j=1}^J n^{-j/2}P_j(z)};
\label{eqn:approx-density}
\end{equation}
i.e., $\nu_n$ is the targeted order $J$ Edgeworth expansion approximation to $\mu_n$ we want to establish. 

We note that $\nu_n$ satisfies
\begin{equation}
\nu_n(\R^{|\Z|}) = \int_{\R^{|\Z|}}\psi_n(z)dz=1\quad\text{and}\quad \int_{\R^{|\Z|}} z\psi_n(z)dz=0,
\label{eqn:comparison-mass-mean}
\end{equation}
where the second integral is computed coordinate by coordinate. Indeed, recall that each $\phi_\Sigma P_j$ is a linear combination of entries of $\phi_\Sigma H_r$ with $r\ge2$. By the definition of the Hermite tensors and integration by parts,
\begin{align*}
\int_{\R^{|\Z|}}\phi_\Sigma(z)H_{r,i_1\ds i_r}(z;\Sigma)\,dz
&=(-1)^r\int_{\R^{|\Z|}}\partial_{i_1}\cdots\partial_{i_r}\phi_\Sigma(z)\,dz=0,\\
\int_{\R^{|\Z|}} z_\ell\phi_\Sigma(z)H_{r,i_1\ds i_r}(z;\Sigma)\,dz
&=(-1)^{r+1}\1{\{\ell=i_1\}}\int_{\R^{|\Z|}}\partial_{i_2}\cdots\partial_{i_r}\phi_\Sigma(z)\,dz=0.
\end{align*}
The boundary terms vanish because Gaussian derivatives are polynomials times the Gaussian density, which vanish at infinity fast enough. The second identity also uses $r\ge 2$ so that the integral involves at least one Gaussian derivative. This implies \eqref{eqn:comparison-mass-mean}.

We prove Theorem \ref{thm:edgeworth} by showing that $\nu_n$ approximates $\mu_n$  up to a uniform error over $\set{f:\norm f_{C_b^2} \leq 1}$. To achieve this, we first note that the dynamics for $\sqrt n Y_n$, viewed as an SA recursion, implies a corresponding update rule for its law $\mu_{n}$ of the form $\mu_{n+1} = \mu_n\cK_{n+1}$ for some operator $\cK_{n+1}$ to be defined next. Since $\nu_n$ is the proposed Edgeworth expansion approximation for $\mu_n$, the defect from one-step update $\nu_{n+1}-\nu_n\cK_{n+1}$ should be small. 

\paragraph{One-step defect of the proposed Edgeworth expansion.} We first write down the dynamics of $\sqrt n Y_n$. Multiplying \eqref{eqn:YmD-recursions} by $\sqrt{n+1}$ gives
\begin{align*}
\sqrt{n+1}Y_{n+1}
&=\sqrt{\frac{n+1}{n}}C_{n+1}(\sqrt nY_n)+\sqrt{n+1}\alpha_{n+1}(\xi_{n+1}+\varepsilon_{n+1})
\end{align*}
where, for notational simplicity, we write $M_{n+1}:=\sqrt{\frac{n+1}{n}}C_{n+1}.$

Note that $\set{(\xi_{n+1},\widehat P_{n+1},C_{n+1}):n\ge0}$ is an independent sequence. Hence, $\set{\sqrt n Y_n:n\ge1}$  is a (time-inhomogeneous) Markov chain. Thus, the law $\mu_{n+1}$ of $\sqrt{n+1}Y_{n+1}$ can be written as $\mu_{n+1}=\mu_n\cK_{n+1}$, where
\begin{equation}\label{eqn:tu_def_cK}
\cK_{n+1}f(z):=\E f\crbk{M_{n+1}z+\sqrt{n+1}\alpha_{n+1}(\xi_{n+1}+\varepsilon_{n+1})}
\end{equation}
for bounded measurable $f$. Since we want to show that $\nu_n$ approximates $\mu_n$, we apply this transition and consider the one-step defect $\nu_n\cK_{n+1}-\nu_{n+1}$. We establish that this defect is indeed of a higher order in the following key lemma.

\begin{lemma}
\label{lemma:one-step}
For all sufficiently large $n$, $\nu_n\cK_{n+1}-\nu_{n+1}$ has a signed density $r_n$ with $\int r_n(z)\,dz=0$ and $\int z r_n(z)\,dz=0$ such that
\begin{equation}
\int(1+|z|^2)|r_n(z)|\,dz\le Cn^{-(J+3)/2}.
\label{eqn:local-error}
\end{equation}
\end{lemma}

The proof of Lemma \ref{lemma:one-step} is long and technical, leveraging a more refined version of the arguments in Proposition \ref{prop:coefficient-recursion}. To streamline the presentation, the proof is deferred to Section \ref{section:proof:lemma:one-step}. 

To simplify notation, we identify the density $r_n$ with the measure and write $r_n(dz) = r_n(z)dz$.

\paragraph{Propagating the one-step defect.}
Equipped with Lemma \ref{lemma:one-step} that bounds the one-step defect incurred by using $\nu_n$ to approximate $\mu_n$, we propagate this error to obtain a bound on $|\mu_n f-\nu_n f|$.

Fix $N\ge1$ large enough that Lemma~\ref{lemma:one-step} applies for every $k\ge N$, and fix $f\in C_b^2(\R^{|\Z|})$. We define for all $N\le k<n,$
\[
f_{k,n}:=\cK_{k+1}\cdots\cK_n f\quad \text{and}\quad f_{n,n}:=f. 
\]
Recalling that $\cK_{n}$ is the transition kernel of the Markov chain $\set{\sqrt n Y_n:n\ge 1}$, $ f_{k,n}(z)$ can be viewed the expected terminal value of the reward function $f$ at time $n$ when the normalized recursion starts from $z$ at time $k$. In particular, it satisfies the backward induction $f_{k,n}=\cK_{k+1}f_{k+1,n}$. 

On the other hand, for $k\ge N$, $\nu_k\cK_{k+1} = \nu_{k+1}+ r_{k}$. Thus, the difference of the two laws satisfies the following recursion
\[
\mu_{k+1}-\nu_{k+1}=(\mu_k-\nu_k)\cK_{k+1}+r_k.
\]
We iterate this from $N$ to $n$ and obtain
\begin{equation}
\mu_nf -\nu_nf=\int f_{N,n}\,d(\mu_N-\nu_N)+\sum_{k=N}^{n-1}\int f_{k+1,n}(z)r_k(z)\,dz.
\label{eqn:tu_mu_nf-nu_nf_propagated}
\end{equation}
We bound the two terms separately, using the fact that both $\mu_N-\nu_N$ and $r_n$ have zero total mass and mean by the definition of $\mu_n$, the centering in $Y_n$, \eqref{eqn:comparison-mass-mean}, and lemma \ref{lemma:one-step}.

We consider a finite signed measure $\rho$ with zero mass and mean and a finite second absolute moment; i.e.$\int_{\R^{|\Z|}}d\rho = 0$, $\int_{\R^{|\Z|}} z\rho(dz) = 0\in\R^{|\Z|}$, and $\int |z|^2|\rho|(dz)<\infty$ where $\abs{\rho}$ denotes the total variation measure. Since $\rho$ has finite mean and mass, subtracting the zero and first order Taylor expansion will not change the integral. Thus, the integral remainder formula for Taylor series implies
\begin{equation}
\begin{aligned}
\abs{\int_{\R^{|\Z|}} f\,d\rho}
&=\abs{\int_{\R^{|\Z|}}\{f(z)-f(0)-\nabla f(0) z\}\,\rho(dz)}\\
&\le\frac12\norm{\nabla^2f}_\infty\int_{\R^{|\Z|}}|z|^2\,|\rho|(dz),
\end{aligned}
\label{eqn:second-order-comparison}
\end{equation}
where we recall that $\norminf{\nabla^2f}$ first take matrix 2-norm then sup norm in $z$. 

Thus, to bound the two terms in \eqref{eqn:tu_mu_nf-nu_nf_propagated} using \eqref{eqn:second-order-comparison}, we need to control the Hessians of the functions $\set{f_{k,n}:N\leq k\leq n}$.

By backward induction, we have
\[
\begin{aligned}
f_{k,n}(z)
&=\E f_{k+1,n}(M_{k+1}z+\sqrt{k+1}\alpha_{k+1}(\xi_{k+1}+\varepsilon_{k+1}))\\
&=\E f_{k+2,n}(M_{k+2}(M_{k+1}z+\sqrt{k+1}\alpha_{k+1}(\xi_{k+1}+\varepsilon_{k+1}))+\sqrt{k+2}\alpha_{k+2}(\xi_{k+2}+\varepsilon_{k+2}))\\
&=\ds \\
&= \E f\crbk{M_n\cdots M_{k+1}z+\sum_{\ell=k+1}^n M_n\cdots M_{\ell+1}(\sqrt{\ell}\alpha_\ell(\xi_\ell+\varepsilon_\ell))}.
\end{aligned}
\]
where the second equality uses the independence of the sampled reward-transition sequence. 

Note that the sum part does not depend on $z$, and the linear part is
\[
M_n\cdots M_{k+1}=\left(\prod_{\ell=k+1}^n\sqrt{\frac{\ell}{\ell-1}}\right)C_n\cdots C_{k+1}=\sqrt{\frac nk}\,C_n\cdots C_{k+1}.
\]
Recall \eqref{eqn:product-bound}. Then, applying chain rule to $f$ yields an integrable bound for the first and second derivatives of the integrand, uniformly in a neighborhood of $z$. Dominated convergence therefore gives
\begin{equation}
\begin{aligned}
\norm{\nabla^2f_{k,n}}_\infty
&\le\norm{\nabla^2f}_\infty\E\left\|\sqrt{\frac nk}\,C_n\cdots C_{k+1}\right\|_2^2\\
&\le C\frac nk\left(\frac{k+1+k_0}{n+k_0}\right)^{2b}\norm{\nabla^2f}_\infty\\
&\le C(k/n)^{2b-1}\norm{\nabla^2f}_\infty.
\end{aligned}
\label{eqn:tu_hessian-contraction}
\end{equation}
where the norm in the expectation is the matrix 2-norm. Since $|\Z|$ is finite, the equivalence of norms and \eqref{eqn:product-bound} yields \eqref{eqn:tu_hessian-contraction} with a factor $C$ that could depend on $N$ and $|\Z|$ but not $f,\;k$, and $n$. 

Applying \eqref{eqn:second-order-comparison} to each term of \eqref{eqn:tu_mu_nf-nu_nf_propagated}, and then using \eqref{eqn:tu_hessian-contraction} and \eqref{eqn:local-error}, yields
\[
\begin{aligned}
\abs{\mu_nf -\nu_n f}
&\le\frac12\norm{\nabla^2f_{N,n}}_\infty\int|z|^2\,|\mu_N-\nu_N|(dz)+\frac12\sum_{k=N}^{n-1}\norm{\nabla^2f_{k+1,n}}_\infty\int|z|^2|r_k(z)|\,dz\\
&\le C\norm{\nabla^2f}_\infty\left\{n^{-(2b-1)}+\sum_{k=N}^{n-1}\left(\frac{k+1}{n}\right)^{2b-1}k^{-(J+3)/2}\right\}.
\end{aligned}
\]
Here, we note that the first two moments are finite because $\sqrt N Y_N$ is almost surely bounded by a constant independent of $n$, while the density $\psi_N$ of $\nu_N$ is a polynomial times a Gaussian density and therefore has sufficiently fast tail decay. 

For $a,c>0$, comparing the sum with the integral of powers, 
\[
n^{-a}\sum_{k=N}^{n-1}k^{a-c-1}\le C\begin{cases}
n^{-c},&a>c,\\
n^{-a}\log(en),&a=c,\\
n^{-a},&a<c.
\end{cases}
\]
Recall that $\mu_n$ is the law of ${\sqrt{n}Y_n}$ and $\psi_n$ is the density of $\nu_n$. Then, with $a=2b-1$, $c=(J+1)/2$, and $\eta=\min\{2b-1,(J+1)/2\}$, we have
\begin{equation}
\abs{\E f(\sqrt nY_n)-\int f(z)\psi_n(z)\,dz}\le C\norm{\nabla^2f}_\infty n^{-\eta}\log(en),
\label{eqn:smooth-error}
\end{equation}
for all $f\in C_b^2(\R^{|\Z|})$. Since $b>1/2$ and $J+1>2b-1$, both exponents defining $\eta$ exceed $b-1/2$. Hence $n^{-\eta}\log(en)=o(n^{1/2-b})$, and taking the supremum over $\norm{f}_{C_b^2}\le1$ proves \eqref{eqn:policy-edgeworth}.

\subsubsection{Proof of Lemma \ref{lemma:one-step}}
\label{section:proof:lemma:one-step}

Recall the definition of $\cK_{n+1}$ in \eqref{eqn:tu_def_cK}. This transition kernel contains the shift $\sqrt{n+1}\alpha_{n+1}\varepsilon_{n+1}$, which includes the deterministic $m_n=O(n^{-b})$. To prove the defect bound, we separate the contribution of this term.

Define
\[
V_{n+1}:=\sqrt{n+1}\alpha_{n+1}\xi_{n+1},\qquad G_{n+1}(z):=M_{n+1}z+V_{n+1},\qquad W_{n+1}:=\sqrt{n+1}\alpha_{n+1}\varepsilon_{n+1}.
\]
Thus the normalized update is $\sqrt{n+1}Y_{n+1}=G_{n+1}(\sqrt nY_n)+W_{n+1}$. Since $C_{n+1}=I+O_{L^\infty}(n^{-1})$ and $\sqrt{(n+1)/n}=1+O(n^{-1})$, we have $M_{n+1}=I+O_{L^\infty}(n^{-1})$. Hence there is a non-random $N\ge1$ such that $M_{n+1}$ is invertible almost surely for all $n\ge N$, with positive determinant. It follows that $G_{n+1}$ is also invertible. Moreover, $m_n$ is deterministic and $\varepsilon_{n+1}=\gamma(\widehat P_{n+1}-P)\Pi^*m_n$, so $\E W_{n+1}=0$ and $|W_{n+1}|\le Cn^{-1/2}|m_n|$.

For $n\ge N$, let $\hat p_{n+1}$ be the signed density obtained by applying $G_{n+1}$ to $\psi_n$. Thus, for every bounded measurable $f$,
\begin{equation}
\int_{\R^{|\Z|}}f(y)\hat p_{n+1}(y)\,dy=\int_{\R^{|\Z|}}f(G_{n+1}(z))\psi_n(z)\,dz.
\label{eqn:updated-density-action}
\end{equation}
This density exists because $G_{n+1}$ is an invertible affine map. Its total variation is $\int_{\R^{|\Z|}}|\psi_n(z)|\,dz\le C$, so Fubini's theorem and a translation give
\begin{equation}
\nu_n\cK_{n+1}f=\E\int_{\R^{|\Z|}}f(G_{n+1}(z)+W_{n+1})\psi_n(z)\,dz=\int_{\R^{|\Z|}}f(y)\E\hat p_{n+1}(y-W_{n+1})\,dy.
\label{eqn:nu_cK_f}
\end{equation}
Consequently, the defect $\nu_n\cK_{n+1}-\nu_{n+1}$ has density
\begin{equation}
r_n(y)=\crbk{\E\hat p_{n+1}(y)-\psi_{n+1}(y)}+\E\sqbk{\hat p_{n+1}(y-W_{n+1})-\hat p_{n+1}(y)}.
\label{eqn:full-density-error}
\end{equation}
We bound the two differences separately. For the first difference, we apply Proposition~\ref{prop:coefficient-recursion} to show that the Fourier transform has vanishing lower-order coefficients. We then invert the Fourier transform to get a density bound. The second difference will be bounded by a Taylor expansion in $W_{n+1}$.

\paragraph{Bounding the Fourier transform of $\E\hat p_{n+1}(y)-\psi_{n+1}(y)$. }
Every coefficient $K_{r,p}$ used to construct $\psi_n$ satisfies $r-1+p\le r-1+2p\le2b$, so its coefficient belongs to $\coeff(J+1)$, where $J+1=\floor{2b}$. Define the residual error
\begin{equation}\label{eqn:tu_def_cR}
\cR_n(v):=e^{\Sigma[v^{[2]}]/2}\sqbk{\exp\set{\sum_{j=1}^J n^{-j/2}\ell_j(v)}-1-\sum_{j=1}^J n^{-j/2}q_j(v)}.
\end{equation}
By \eqref{eqn:Psi_J_1} and \eqref{eqn:correction-transform},
\[
\int_{\R^{|\Z|}}e^{v^\top z}\psi_n(z)\,dz=\exp\set{\Psi_n^{J+1}(\sqrt n\,v)}-\cR_n(v).
\]
for $v\in i\R^{|\Z|}$ within a neighborhood of the origin. Applying \eqref{eqn:updated-density-action} with $f(y)=e^{v^\top y}$, we have
\[
\begin{aligned}
\int_{\R^{|\Z|}}e^{v^\top y}\hat p_{n+1}(y)\,dy
&=\int_{\R^{|\Z|}}e^{v^\top G_{n+1}(z)}\psi_n(z)\,dz\\
&=e^{v^\top V_{n+1}}\int_{\R^{|\Z|}}e^{(M_{n+1}^\top v)^\top z}\psi_n(z)\,dz\\
&=e^{v^\top V_{n+1}}\sqbk{\exp\set{\Psi_n^{J+1}(\sqrt n\,M_{n+1}^\top v)}-\cR_n(M_{n+1}^\top v)}.
\end{aligned}
\]

To proceed, we define the following error terms
\[
\begin{aligned}
I_1(v)&:=\log\E\exp\set{v^\top V_{n+1}+\Psi_n^{J+1}(\sqrt n\,M_{n+1}^\top v)-\Psi_{n+1}^{J+1}(\sqrt{n+1}\,v)},\\
I_2(v)&:=\cR_{n+1}(v)-\E\sqbk{e^{v^\top V_{n+1}}\cR_n(M_{n+1}^\top v)}.
\end{aligned}
\]
Taking expectations in the preceding identity and subtracting the Fourier transform of $\psi_{n+1}$ yields
\begin{equation}
\int_{\R^{|\Z|}}e^{v^\top y}\sqbk{\E\hat p_{n+1}(y)-\psi_{n+1}(y)}\,dy=\exp\set{\Psi_{n+1}^{J+1}(\sqrt{n+1}\,v)}\crbk{e^{I_1(v)}-1}+I_2(v).
\label{eqn:one-step-transform-decomposition}
\end{equation}
Thus, to bound the left-hand side of \eqref{eqn:one-step-transform-decomposition}, we show that both $I_1$ and $I_2$ are $O(n^{-(J+3)/2})$, uniformly for imaginary $v$ in a neighborhood of zero.

We first give a brief discussion on bounding $I_1$ using the structure given by Proposition~\ref{prop:coefficient-recursion}. Since $\sqrt n\,M_{n+1}=\sqrt{n+1}\,C_{n+1}$, $I_1(v)$ is the discrepancy $\Delta_n^{J+1}(\sqrt{n+1}\,v)$ with $\varepsilon_{n+1}$ omitted. The argument preceding \eqref{eqn:tu_form_of_trDelta} in the proof of Proposition~\ref{prop:coefficient-recursion} shows that this omission changes each unnormalized degree-$r$ coefficient by $O(n^{-(r+b)})$. Therefore, for every $K_{r,q}\in\coeff(J+1)$ and each fixed imaginary $v$ near zero, \eqref{eqn:coefficient-induction} gives
\begin{equation}
\begin{aligned}
[t^r]I_1(tv)
&=(n+1)^{r/2}\crbk{[t^r]\Delta_n^{J+1}(tv)+O(n^{-(r+b)})}\\
&=(n+1)^{r/2}\crbk{o(n^{-(r+q)})+O(n^{-(r+b)})}\\
&=o(n^{-(r+2q)/2}),
\end{aligned}
\label{eqn:normalized-discrepancy}
\end{equation}
where $q\le b-(r-1)/2<b$ by the definition of $\coeff(J+1)$. The factor $(n+1)^{r/2}$ follows from homogeneity in tensor arguments at degree-$r$ as we extract the normalization $\sqrt{n+1}$.

Then the natural step is to sum \eqref{eqn:normalized-discrepancy} over $r$ to get a bound on $I_1$. However, there are several caveats. First, \eqref{eqn:normalized-discrepancy} is picking up only the $t^r$-coefficients, not taking into account for other $n$ dependence. Second, it is currently only shown to hold only for fixed $v$, not uniformly in $|v|<\delta$. Third, and more importantly, even though the first two issues are addressed, choosing the largest possible $q$ such that $K_{r,q}\in \coeff(J+1)$, we have that $[t^r]I_1(tv) = o(n^{-(J+1)/2})$ or $o(n^{-(J+2)/2})$ depending on if $J+2-r$ is even. This cannot justify a $O(n^{-(J+3)/2})$ expansion. 

The idea to resolve these issues is to realize that $I_1(v)$ is the log expectation of an exponential of a polynomial with random but bounded coefficients. Moreover the polynomial coefficients are smooth in $h = n^{-1/2}$ as well. This upgrades the pointwise bound to a uniform bound in $v$ and little-$o$ to big-$O$ with one higher power of $n^{-1/2}$. Furthermore, to resolve the $o(n^{-(J+1)/2})$ case, we need to note a special structure of the power-of-$n^{-1/2}$ beyond just integer power. We isolate the proof of this bound by introducing the following lemma for readability. 

\begin{lemma}\label{lemma:I1_bd}
There exists $\delta > 0$ such that
\begin{equation}
I_1(v)=O(n^{-(J+3)/2})
\label{eqn:normalized-cancellation}
\end{equation}
uniformly for imaginary $v\in i\R^{|\Z|}$ with $|v|\le\delta$. 
\end{lemma}

\begin{proof}
To pass from these fixed-degree estimates to a uniform bound on $I_1$, we express all its time dependence through a real parameter $h$ near zero. Use the same bounded random pair $(\widehat P_1,\xi_1)$ for every $h$ and define
\[
\begin{aligned}
M(h)&:=\sqrt{1+h^2}\sqbk{I+\frac{\alpha h^2}{1+(k_0+1)h^2}(\gamma\widehat P_1\Pi^*-I)},\\
V(h)&:=\frac{\alpha h\sqrt{1+h^2}}{1+(k_0+1)h^2}\xi_1.
\end{aligned}
\]
At $h=n^{-1/2}$, these have the same joint law as $(M_{n+1},V_{n+1})$. For real $h$ sufficiently close to zero, these functions are $C^\infty$, with derivatives of every fixed order uniformly bounded. 

Substituting \eqref{eqn:Psi_J_1} into the definition of $I_1$, and using $(n+1)^{-j/2}=h^j(1+h^2)^{-j/2}$ at $h=n^{-1/2}$, gives $I_1(v)=\log\E e^{U(n^{-1/2},v)},$
where
\[
\begin{aligned}
U(h,v):=v^\top V(h)+\frac12v^\top\sqbk{M(h)\Sigma M(h)^\top-\Sigma}v+\sum_{j=1}^J h^j\sqbk{\ell_j(M(h)^\top v)-(1+h^2)^{-j/2}\ell_j(v)}.
\end{aligned}
\]
We note that $U(0,v)=0$. Thus all dependence on $n$ is now explicit through $h$.

Since the $\ell_j$ are polynomials in $v$ and $(M(h),V(h))$ are $C^\infty$, $U(h,v)$ is jointly $C^\infty$ in $(h,v)$. Thus, for fixed sufficiently small $h_0,\delta>0$,$U$ and its mixed derivatives of every fixed order are uniformly bounded w.p.1 on $|h|\le h_0$, $|v|\le\delta$. This justifies the interchange of derivative and expectation and thus $\log\E e^{U(h,v)}$ is jointly $C^\infty$. 

We then argue that $U$ is a random polynomial in $v$ and each degree-$r$ term can be factored as $h^r$ times a polynomial in $v$ whose coefficients are smooth functions of $h^2$, with derivatives uniformly bounded near $h=0$; i.e.
\begin{equation}\label{eqn:tu_U_factorization}
[v_{i_1}v_{i_2}\ds v_{i_r}]U(h,v)= h^r R_{i_1,\ds,i_r}(h^2)
\end{equation}
where $R_{i_1,\ds,i_r}$ is a random function and its derivatives of any fixed order are uniformly bounded for $|h|<h_0$. 
Indeed, $M(h)$ is a smooth function of $h^2$ with $M(0)=I$, while $V(h)$ is $h$ times such a function. Thus the linear and quadratic terms of $U$ have the required factors $h$ and $h^2$, respectively. For the sum term, observe that
\[
\ell_j(M(h)^\top v)-(1+h^2)^{-j/2}\ell_j(v)=\sqbk{\ell_j(M(h)^\top v)-\ell_j(v)}+\sqbk{1-(1+h^2)^{-j/2}}\ell_j(v).
\]
We consider the degree-$r$ component of the two differences. By the definition of $\ell_j$ in \eqref{eqn:def_ell_j}, the degree-$r$ component of $\ell_j$ is $K_{r,p}[v^{[r]}]/r!$, where $j=r-2+2p$. 

For the first difference, writing $e_i$ for the $i$th coordinate vector,
\[
[v_{i_1}\cdots v_{i_r}]\sqbk{\ell_j(M(h)^\top v)-\ell_j(v)}= C\sqbk{K_{r,p}[M(h)^\top e_{i_1},\ldots,M(h)^\top e_{i_r}]-K_{r,p}[e_{i_1},\ldots,e_{i_r}]}.
\]
where $C$ only depends on $r$ and the coordinates $i_1,\ds ,i_r$. 
Since $M(h)=I+h^2A(h^2)$ for some smooth matrix-valued random function $A$, multilinearity of $K_{r,p}$ shows that the right-hand side is $h^2$ times a smooth function of $h^2$, with uniformly bounded derivatives. For the second difference, note that $1-(1+h^2)^{-j/2} = O(h^2)$ is $C^\infty$ with bounded derivatives for $|h|\leq h_0$. Thus, after multiplication by $h^j$, the contribution of this degree-$r$ coefficient to $U$ therefore has the form $h^{j+2}R'_{i_1,\ldots,i_r}(h^2)
=h^r[(h^2)^pR'_{i_1,\ldots,i_r}(h^2)],$
where $R'_{i_1,\ldots,i_r}$ is smooth with uniformly bounded derivatives. This proves the required factorization in \eqref{eqn:tu_U_factorization}.

As a consequence, we have that uniformly for $0\le t\le1$, $|v|\leq\delta$,  $\norm{\partial_t^r U(h,tv)}_{L^\infty}\le C|h|^r$ for all $1\le r\le J+3$.
Since the exchange of expectation and derivatives is justified, we have \begin{equation}\label{eqn:tu_derivative_log_E_exp_bd_in_h}
    \abs{\partial_t^r\log\E e^{U(h,tv)}}\le C|h|^r. 
\end{equation}

With these properties, we get back to analyzing $I_1$. To simplify notation, define
\[
T_r(h,v):=\frac1{r!}\left.\partial_t^r\log\E e^{U(h,tv)}\right|_{t=0},
\]
where $T_r$ is $C^\infty$, inheriting the smoothness from $\log\E e^{U(h,v)}$. Then, $\log\E e^{U(h,tv)}$ is also $C^\infty$ in $(h,t)$ and we apply Taylor's theorem at $t=0$ and evaluated at $t=1, h=n^{-1/2}$ to obtain 
\begin{equation}
I_1(v)=\log\E e^{U(n^{-1/2},1v)}=\sum_{r=2}^{J+2}T_r(n^{-1/2},v)+O(n^{-(J+3)/2}),
\label{eqn:finite-logarithmic-expansion}
\end{equation}
where the remainder is uniform in $|v|\le \delta$. Here, we note that the constant and linear coefficients vanish because $U(h,0)=0$ and $\E V(h)=0$.

It remains to bound the early $T_r$ terms using \eqref{eqn:normalized-discrepancy}. The factorization established in \eqref{eqn:tu_U_factorization} and the uniform boundedness of the derivatives that enable the exchange of integral and derivatives then implies
\[
T_r(h,v)=h^r b_r(h^2,v),
\]
where $b_r(s,v)$ is smooth in $s$ near zero, with derivatives of every fixed order bounded uniformly for $|v|\le\delta$. 

For $2\le r\le J+2$, set $p_r:=\floor{(J+2-r)/2}$. Since $r-1+p_r\le r-1+2p_r\le J+1\le2b$, we have $K_{r,p_r}\in\coeff(J+1)$. Therefore, \eqref{eqn:normalized-discrepancy} implies
\[
T_r(n^{-1/2},v)=[t^r]I_1(tv)=o(n^{-(r+2p_r)/2}).
\]
This forces $b_r(n^{-1},v)=o(n^{-p_r})$. 

On the other hand, Taylor's theorem applied to $b_r$ in $s$ gives
\[
b_r(s,v)=\sum_{q=0}^{p_r}\frac{s^q}{q!}\partial_s^q b_r(0,v)+O(|s|^{p_r+1}),
\]
Thus, for each fixed $|v|\leq \delta$, by setting $s=n^{-1}$ we see that the expansion of $b_r$ up to and including $p_r$th derivative must vanish, otherwise it would contradict $b_r(n^{-1},v)=o(n^{-p_r})$. Thus, 
\[
|T_r(h,v)|\le C|h|^{r+2p_r+2}\le C|h|^{J+3},
\]
where the last inequality uses the definition of $p_r$. Substituting these bounds on $T_r$ into \eqref{eqn:finite-logarithmic-expansion} proves the claim.

\end{proof}

For $I_2$, we further decompose it as follows:
\begin{equation}
I_2(v)=\crbk{\cR_{n+1}(v)-\cR_n(v)}+\cR_n(v)\crbk{1-\E e^{v^\top V_{n+1}}}-\E\sqbk{e^{v^\top V_{n+1}}\crbk{\cR_n(M_{n+1}^\top v)-\cR_n(v)}}.
\label{eqn:tu_I2_decomp}
\end{equation}
To bound these terms, fix sufficiently small $\delta>0$ and consider imaginary $v$ with $|v|\le\delta$. By the Bell polynomial construction \eqref{eqn:bell}, the analytic function
\[
\exp\set{\sum_{j=1}^J h^j\ell_j(v)}-1-\sum_{j=1}^J h^jq_j(v)
\]
has vanishing Taylor coefficients through order $J$ in $h$. It is therefore $O(h^{J+1})$, its $h$-derivative is $O(h^J)$, and its $v$-gradient is $O(h^{J+1})$, uniformly on a larger complex ball, say $|v|\le2\delta$. Then, we extend the definition of $\cR_n$ in \eqref{eqn:tu_def_cR} to real $t\ge1$ by replacing $n$ with $t$. Using $h=t^{-1/2}$ and thus $dh/dt=-\frac12t^{-3/2}$, 
\[
\begin{aligned}
|\cR_t(v)|+\norm{\nabla_v\cR_t(v)}&\le Ct^{-(J+1)/2},\\
|\partial_t\cR_t(v)|&\le Ct^{-3/2}(t^{-1/2})^J=Ct^{-(J+3)/2}
\end{aligned}
\]uniformly on $|v|\le2\delta$. 
Thus the time increment term satisfies
\begin{equation}
|\cR_{n+1}(v)-\cR_n(v)|=\abs{\int_n^{n+1}\partial_t\cR_t(v)\,dt}\le C\int_n^{n+1}t^{-(J+3)/2}\,dt\le Cn^{-(J+3)/2}.
\label{eqn:tu_I2_time_increment}
\end{equation}

For the third term in \eqref{eqn:tu_I2_decomp}, $\norm{M_{n+1}-I}_{L^\infty}=O(n^{-1})$ implies that, for all sufficiently large $n$, for all $0\le s\le1$
\[
|v+s(M_{n+1}^\top-I)v|\le\crbk{1+\norm{M_{n+1}-I}}|v|\le2\delta
\]
w.p.1. Thus, integrating the gradient along this segment, we have
\begin{equation}
\begin{aligned}
\norm{\cR_n(M_{n+1}^\top v)-\cR_n(v)}_{L^\infty}
&\le\sup_{|w|\le2\delta}\norm{\nabla_w\cR_n(w)}\,\norm{(M_{n+1}^\top-I)v}_{L^\infty}\\
&\le Cn^{-(J+1)/2}n^{-1}=Cn^{-(J+3)/2}.
\end{aligned}
\label{eqn:tu_I2_RM_R}
\end{equation}

Finally, $\E V_{n+1}=0$ and $|V_{n+1}|\le Cn^{-1/2}$ imply
\[
\begin{aligned}
|\cR_n(v)|\abs{\E e^{v^\top V_{n+1}}-1}
&\le Cn^{-(J+1)/2}\abs{\E\sqbk{e^{v^\top V_{n+1}}-1-v^\top V_{n+1}}}\\
&\le Cn^{-(J+1)/2}\E|v^\top V_{n+1}|^2\le Cn^{-(J+3)/2}.
\end{aligned}
\]
Combining these bounds in \eqref{eqn:tu_I2_decomp}, using \eqref{eqn:tu_I2_time_increment}, \eqref{eqn:tu_I2_RM_R}, and $|e^{v^\top V_{n+1}}|=1$, we see that $I_2(v)=O(n^{-(J+3)/2})$ uniformly for imaginary $v$ with $|v|\le\delta$. 

By \eqref{eqn:Psi_J_1}, the deterministic exponential $\exp\set{\Psi_{n+1}^{J+1}(\sqrt{n+1}\,v)}$ in \eqref{eqn:one-step-transform-decomposition} is bounded on $|v|\le\delta$. Thus, combining \eqref{eqn:normalized-cancellation}, the bound on $I_2$, and $|e^h-1|\le |h|e^{|h|}$ therefore yields
\begin{equation}
\int_{\R^{|\Z|}}e^{v^\top y}\sqbk{\E\hat p_{n+1}(y)-\psi_{n+1}(y)}\,dy=O(n^{-(J+3)/2}),
\label{eqn:one-step-transform}
\end{equation}
for $v\in i\R^{|\Z|}$ with $|v|\le\delta$.
\paragraph{Obtaining the weighted density bound by Fourier inversion.}
We proceed to get the weighted density bound by inverting the Fourier transform obtained above. However, the local bound \eqref{eqn:one-step-transform}, which hold only for $|v|\leq\delta$, is not by itself sufficient for Fourier inversion. To invert, we use the Gaussian-times-polynomial form of the transform to obtain an integrable remainder bound at all frequencies. This provides the uniformity needed.

We employ a Taylor expansion of the Fourier transform of the density error in both $v$ and $n^{-1/2}$. To do this, we again replace the discrete $n^{-1/2}$ with real $h$ near zero and recall the definition of $(M(h),V(h))$. Then, for $v\in i\R^{|\Z|}$, we consider the following quantity
\begin{equation}
\begin{aligned}
\mathcal F(h,v)
&:=\E\sqbk{\exp\set{v^\top V(h)+\frac{1}{2}v^\top M(h)\Sigma M(h)^\top v}\crbk{1+\sum_{j=1}^J h^jq_j(M(h)^\top v)}}\\
&\quad-e^{v^\top\Sigma v/2}\crbk{1+\sum_{j=1}^J h^j(1+h^2)^{-j/2}q_j(v)}.
\end{aligned}
\label{eqn:tu_continuous_fourier_defect}
\end{equation}

Set $h=n^{-1/2}$. Then, by \eqref{eqn:correction-transform} and \eqref{eqn:updated-density-action} with $f(y)=e^{v^\top y}$, we have
\[
\begin{aligned}
\int_{\R^{|\Z|}}e^{v^\top y}\E\hat p_{n+1}(y)\,dy
&=\E\sqbk{e^{v^\top V_{n+1}}\int_{\R^{|\Z|}}e^{(M_{n+1}^\top v)^\top z}\psi_n(z)\,dz}\\
&=\E\sqbk{\exp\set{v^\top V(h)+\frac12v^\top M(h)\Sigma M(h)^\top v}\crbk{1+\sum_{j=1}^J h^jq_j(M(h)^\top v)}}.
\end{aligned}
\]
Moreover, since $(n+1)^{-j/2}=h^j(1+h^2)^{-j/2}$,
\[
\int_{\R^{|\Z|}}e^{v^\top y}\psi_{n+1}(y)\,dy=e^{v^\top\Sigma v/2}\sqbk{1+\sum_{j=1}^J h^j(1+h^2)^{-j/2}q_j(v)}.
\]
Therefore, subtracting these two identities gives exactly the two terms defining $\mathcal F(h,v)$; i.e. for $v\in i\R^{|\Z|}$, 
\begin{equation}
\mathcal F(n^{-1/2},v)=\int_{\R^{|\Z|}}e^{v^\top y}\sqbk{\E\hat p_{n+1}(y)-\psi_{n+1}(y)}\,dy.
\label{eqn:tu_exact_fourier_defect}
\end{equation}

Since $M(0)=I$, $\Sigma\succ0$, and $v$ is imaginary, there are deterministic $h_0,c>0$ such that $v^\top M(h)\Sigma M(h)^\top v\leq- c|v|^2$  uniformly for $|h|\le h_0$ and $v\in i\R^{|\Z|}$ almost surely. Moreover, $M(h)$ and $V(h)$, seen as a random function in $h$, have all their derivatives of any fixed order uniformly bounded for $|h|\le h_0$. Also, recall that  $q_j$ are polynomials of finite degree in $v$. Thus, $\cF(h,v)$ is $C^\infty$ on $[-h_0,h_0]\times i\R^{|\Z|}$ and differentiating \eqref{eqn:tu_continuous_fourier_defect} in $h$ and $v$ and evaluate within $|h|\le h_0$ and $v\in i\R^{|\Z|}$ only introduces polynomial growth factors in $v$. So, for every nonnegative integer $k$ and multi-index $\beta$, there is $C>0$ and a finite integer $L\ge 0$ that can depend on $\beta$ and $k$ but not $h$ and $v$ such that
\begin{equation}\label{eqn:tu_partial_cF_bd}
|\partial_h^k\partial_v^\beta\mathcal F(h,v)|\le C(1+|v|)^Le^{-c|v|^2}
\end{equation}
for all $|h|\le h_0$ and $v\in i \R^{|\Z|}$. Here we note that the uniform boundedness and the Gaussian-times-polynomial form of the derivative justifies the interchange of differentiation and expectation.

We now Taylor expand $\cF$ and use the identity \eqref{eqn:one-step-transform} to eliminate the low-order Taylor coefficients in $h$. By \eqref{eqn:one-step-transform} and \eqref{eqn:tu_exact_fourier_defect}, $\mathcal F(n^{-1/2},v)=O(n^{-(J+3)/2})$ for $|v|\le\delta$. By \eqref{eqn:tu_partial_cF_bd}, we can apply Taylor's formula with integral remainder at $h=0$ to $\cF(\cd,v)$, which then implies

\[\partial_h^k\mathcal F(0,v)=0,\quad\forall\  0\le k\le J+2 \text{ and } |v|\le\delta.
\]
Indeed, if there is a nonzero derivative with $k\le J+2$ then, the Taylor expansion of $\cF(\cd,v)$ evaluated at $n^{-1/2}$ would have a term of $\Omega(n^{-k/2})$, contradicting the bound $\mathcal F(n^{-1/2},v)=O(n^{-(J+3)/2})$.

We next extend this property to every $v\in i\R^{|\Z|}$. The key is that, from the definition of $\cF$, $\partial_h^k\mathcal F(0,v)$ is a Gaussian factor in $v$ multiplied by a polynomial in $v$. Thus, if it is $0$ within a neighborhood of the origin, it is identically 0. To verify this, we first observe  
from \eqref{eqn:tu_continuous_fourier_defect} that the other term except 
\[
\exp\set{v^\top V(h)+\frac12v^\top M(h)\Sigma M(h)^\top v}=e^{v^\top\Sigma v/2}e^{g(h,v)},
\]
where
\[
g(h,v):=v^\top V(h)+\frac12v^\top\sqbk{M(h)\Sigma M(h)^\top-\Sigma}v
\]
is already of the desired form. So, it suffices to show that $h$-derivatives of $e^{g(h,v)}$ at $h=0$ is a polynomial in $v$. Since $M(0)=I$ and $V(0)=0$, we have $g(0,v)=0$. Moreover, for every $k\ge1$,
\[
\partial_h^k g(0,v)=v^\top V^{(k)}(0)+\frac12v^\top\left.\frac{d^k}{dh^k}\sqbk{M(h)\Sigma M(h)^\top}\right|_{h=0}v,
\]
which is a polynomial of degree at most two in $v$. On the other hand, since $e^{g(0,v)} = 1$, by chain rule, the derivative of $\del_h^k e^{g(h,v)}|_{h=0}$ is the sum of products of the derivatives $\del_h^kg(0,v)$, hence a polynomial in $v$. For example,
\[
\left.\partial_h e^{g(h,v)}\right|_{h=0}=\partial_h g(0,v),\quad
\left.\partial_h^2 e^{g(h,v)}\right|_{h=0}=\partial_h^2 g(0,v)+\crbk{\partial_h g(0,v)}^2.
\]

Therefore, we conclude that $\partial_h^k\mathcal F(0,v)$ is a positive Gaussian factor $e^{v^\top \Sigma v/2}$ multiplied by a polynomial in $v$. Since a polynomial vanishing on a real open ball is identically zero, we have
\begin{equation}\label{eqn:tu_cF_h_derivative_0}
\partial_h^k\mathcal F(0,v)=0,\quad\forall\  0\le k\le J+2 \text{ and } v\in i\R^{\Z}
\end{equation}
and all $v$-derivatives of $\partial_h^k\mathcal F(0,v)$ also vanish.

By \eqref{eqn:tu_cF_h_derivative_0}, for any multi-index $\beta$, we Taylor expand $\partial_v^\beta\mathcal F(h,v)$ at $h=0$ to order $J+2$ and realize that all the terms through $h^{J+2}$ vanish. Thus, the integral remainder formula implies that
\begin{equation}\label{eqn:tu_global_fourier_remainder}
\begin{aligned}
\abs{\partial_v^\beta\mathcal F(h,v)}&\leq\frac{1}{(J+2)!}\int_0^h(h-\rho)^{J+2}\abs{\partial_h^{J+3}\partial_v^\beta\mathcal F(\rho,v)}d\rho\\
&=\frac{|h|^{J+3}}{(J+2)!}\int_0^1(1-s)^{J+2}\abs{\partial_h^{J+3}\partial_v^\beta\mathcal F(sh,v)}ds\\
&\leq C|h|^{J+3} (1+|v|)^Le^{-c|v|^2}
\end{aligned}
\end{equation}
for all $|h|\le h_0$ and $v\in i\R^{|\Z|}$. Here, the second line follows from a change of variable $\rho=sh$, the last inequality applies \eqref{eqn:tu_partial_cF_bd}, and $(C,c,L)$ depends only on $\beta$, not on $h$ or $v$.  

Choose an integer $\ell>(|\Z|+2)/2$ and let $\Delta_v=\sum_{j=1}^{|\Z|}\partial_{v_j}^2$. By \eqref{eqn:tu_exact_fourier_defect}, Fourier inversion and the differentiation identity for Fourier transforms imply
\[
(1+|y|^2)^\ell\sqbk{\E\hat p_{n+1}(y)-\psi_{n+1}(y)}
=\frac{1}{(2\pi i)^{|\Z|}}\int_{i\R^{|\Z|}}e^{-v^\top y}(1+\Delta_v)^\ell\mathcal F(n^{-1/2},v)\,dv.
\]
To see this, we apply integration by parts and note that the boundary term vanishes due to the Gaussian-polynomial  bounds on the derivatives \eqref{eqn:tu_partial_cF_bd}. Moreover, applying \eqref{eqn:tu_global_fourier_remainder} to the derivatives of $\cF$, we have
\[
\begin{aligned}
(1+|y|^2)^\ell\abs{\E\hat p_{n+1}(y)-\psi_{n+1}(y)}&\le\frac{1}{(2\pi)^{|\Z|}}\int_{\R^{|\Z|}}\abs{(1+\Delta_v)^\ell\mathcal F(n^{-1/2},iu)}\,du\\
&\le Cn^{-(J+3)/2}\int_{\R^{|\Z|}}(1+|u|)^Le^{-c|u|^2}\,du\\
&\le Cn^{-(J+3)/2}.
\end{aligned}
\]
Consequently,
\begin{equation}
\begin{aligned}
\int_{\R^{|\Z|}}(1+|y|^2)\abs{\E\hat p_{n+1}(y)-\psi_{n+1}(y)}\,dy
&\le Cn^{-(J+3)/2}\int_{\R^{|\Z|}}(1+|y|^2)^{1-\ell}\,dy\\
&\le Cn^{-(J+3)/2},
\end{aligned}
\label{eqn:matched-density-error}
\end{equation}
where the last integral is finite by the choice of $\ell$.

\paragraph{Bounding $\E\sqbk{\hat p_{n+1}(y-W_{n+1})-\hat p_{n+1}(y)}$.}
It remains to bound the shift contribution in \eqref{eqn:full-density-error}. Taylor's formula with a change of variable as in \eqref{eqn:tu_global_fourier_remainder} implies
\begin{equation}\label{eqn:tu_pyW-py}
\begin{aligned}
\hat p_{n+1}(y-W_{n+1})-\hat p_{n+1}(y)=-\nabla\hat p_{n+1}(y)W_{n+1}+\int_0^1(1-s)\nabla^2\hat p_{n+1}(y-sW_{n+1})[W_{n+1}^{[2]}]\,ds.
\end{aligned}
\end{equation}
Since $\E W_{n+1}=0$, we have
\[
\E\sqbk{\nabla\hat p_{n+1}(y)W_{n+1}}=\E\sqbk{\crbk{\nabla\hat p_{n+1}(y)-\nabla\psi_n(y)}W_{n+1}}.
\]

Recall the definition in \eqref{eqn:updated-density-action} and that for $n\ge N$ sufficiently large, $M_{n+1}$ and hence $G_{n+1}$ is invertible. Thus, applying the change of variable formula gives
\begin{equation}
\hat p_{n+1}(y)=\frac{\psi_n(G_{n+1}^{-1}(y))}{\det M_{n+1}}.
\label{eqn:updated-density}
\end{equation}
Differentiating \eqref{eqn:updated-density} and adding and subtracting $\nabla\psi_n(G_{n+1}^{-1}(y))$, we obtain
\begin{equation}\label{eqn:tu_grad_p-grad_psi}
\begin{aligned}
\nabla\hat p_{n+1}(y)-\nabla\psi_n(y)
=\nabla\psi_n(G_{n+1}^{-1}(y))\sqbk{\frac{M_{n+1}^{-\top}}{\det M_{n+1}}-I}+\sqbk{\nabla\psi_n(G_{n+1}^{-1}(y))-\nabla\psi_n(y)}.
\end{aligned}
\end{equation}

Since $M_{n+1}=I+O_{L^\infty}(n^{-1})$, the matrix factor in the preceding decomposition is $O_{L^\infty}(n^{-1})$, and
\begin{equation}
|G_{n+1}^{-1}(y)-y|\le Cn^{-1}|y|+Cn^{-1/2}\le Cn^{-1/2}(1+|y|).\label{eqn:tu_Gy_y_bound}
\end{equation}
By \eqref{eqn:approx-density} and the Gaussian-times-polynomial derivative property used earlier, the first two derivatives of $\psi_n$ have a common bound $C(1+|y|)^Le^{-c|y|^2}$ for sufficiently large $n$ when evaluated along the line segment $\set{(1-s)y+sG_{n+1}^{-1}(y):s\in[0,1]}$, or at $G_{n+1}^{-1}(y-sW_{n+1})$. 

Thus, integrating the Hessian along the segment, we have
\[
\begin{aligned}
\abs{\nabla\psi_n(G_{n+1}^{-1}(y))-\nabla\psi_n(y)}&\leq\int_0^1\abs{\crbk{G_{n+1}^{-1}(y)-y}^\top \nabla^2
\psi_n((1-s)y+sG_{n+1}^{-1}(y))}ds\\
&\leq Cn^{-1/2}(1+|y|)^Le^{-c|y|^2}.
\end{aligned}
\]
Here, we use \eqref{eqn:tu_Gy_y_bound} and the derivative bounds to control the norm of the Hessian-displacement product by $Cn^{-1/2}(1+|y|)^{L+1}e^{-c|y|^2}$. Combining this and \eqref{eqn:tu_grad_p-grad_psi} with the matrix factor bound, and again use the fact that the gradient of $\psi_n$ is a polynomial-Gaussian product, we get that
\[
|\nabla\hat p_{n+1}(y)-\nabla\psi_n(y)|\le Cn^{-1/2}(1+|y|)^Le^{-c|y|^2}
\]
for some possibly different $L$. 

Similarly, differentiating \eqref{eqn:updated-density} twice yields
\[
\nabla^2\hat p_{n+1}(y-sW_{n+1})=\frac{M_{n+1}^{-\top}\nabla^2\psi_n(G_{n+1}^{-1}(y-sW_{n+1}))M_{n+1}^{-1}}{\det M_{n+1}}.
\]
The matrix factors are uniformly bounded, so its norm is at most $C(1+|y|)^Le^{-c|y|^2}$, uniformly in $s\in[0,1]$. Multiplying these pointwise bounds by $1+|y|^2$ and integrating therefore gives, almost surely with deterministic constants independent of $n$,
\begin{equation}\label{eqn:tu_1+y_two_terms}
\begin{aligned}
\int_{\R^{|\Z|}}(1+|y|^2)|\nabla\hat p_{n+1}(y)-\nabla\psi_n(y)|\,dy&\le Cn^{-1/2},\\
\sup_{0\le s\le1}\int_{\R^{|\Z|}}(1+|y|^2)\norm{\nabla^2\hat p_{n+1}(y-sW_{n+1})}\,dy&\le C.
\end{aligned}
\end{equation}
where $\norm{\cd}$ denotes the matrix 2-norm.

Recall the definition $|W_{n+1}| = \sqrt{n+1}\alpha_{n+1}|\varepsilon_{n+1}|\le Cn^{-1/2}|m_n|$ and $m_n=O(n^{-b})$. Then, the weighted integral of \eqref{eqn:tu_pyW-py} can be bounded by
\[
\begin{aligned}
&\int_{\R^{|\Z|}}(1+|y|^2)\abs{\E\sqbk{\hat p_{n+1}(y-W_{n+1})-\hat p_{n+1}(y)}}\,dy\\
&\leq \E\sqbk{\abs{W_{n+1}}\int_{\R^{|\Z|}}(1+|y|^2)\abs{\nabla\hat p_{n+1}(y)-\nabla\psi_n(y)}dy} \\
&+ \int_0^1(1-s) \E\sqbk{ |W_{n+1}|^2\sup_{0\le \rho\le1}\int_{\R^{|\Z|}}(1+|y|^2)\norm{\nabla^2\hat p_{n+1}(y-\rho W_{n+1})}dy}ds.\end{aligned}
\]
These resulting two terms can be further bounded using \eqref{eqn:tu_1+y_two_terms} as
\begin{equation}
\begin{aligned}
\int_{\R^{|\Z|}}(1+|y|^2)\abs{\E\sqbk{\hat p_{n+1}(y-W_{n+1})-\hat p_{n+1}(y)}}\,dy
&\le Cn^{-1}|m_n|+Cn^{-1}|m_n|^2\\
&\le Cn^{-1-b}\\
&\le Cn^{-(J+3)/2},
\end{aligned}
\label{eqn:W_density_error}
\end{equation}
where the last inequality uses $J\le2b-1$. 

\paragraph{Concluding Lemma \ref{lemma:one-step}.}
Applying \eqref{eqn:matched-density-error} and \eqref{eqn:W_density_error} to \eqref{eqn:full-density-error} proves \eqref{eqn:local-error}. 

In addition, we verify that $r_n$ has zero mass and mean. Since $\cK_{n+1}$ is a Markov kernel, it preserves the mass of $\nu_n$. So, $$\int_{\R^{|\Z|}} r_n(z)dz = \nu_n\cK_{n+1}(\R^{|\Z|})-\nu_{n+1}(\R^{|\Z|})= 0.$$
Moreover, applying \eqref{eqn:nu_cK_f} with $f(y) = y$, we obtain
\[
\begin{aligned}
\nu_n\cK_{n+1} f
&=\E\int_{\R^{|\Z|}}\crbk{G_{n+1}(z)+W_{n+1}}\psi_n(z)\,dz\\
&=\E\int_{\R^{|\Z|}}\crbk{M_{n+1}z+V_{n+1}+W_{n+1}}\psi_n(z)\,dz\\
&=\E\sqbk{V_{n+1}+W_{n+1}}\\
&=0,
\end{aligned}
\]
where the third equality uses \eqref{eqn:comparison-mass-mean} and the last equality employs $\E\xi_{n+1}=0$ and $\E W_{n+1}=0$. Since $\nu_n$ also has mean zero by \eqref{eqn:comparison-mass-mean}, $r_n$ has zero mean as well.

\end{document}